\documentclass[draft=false,bibliography=totoc,index=totoc]{scrartcl}

\usepackage{amsmath}
\usepackage{amssymb}
\usepackage{amsthm}
\usepackage{cite}
\usepackage[T1]{fontenc}
\usepackage{hyperref}
\usepackage{lmodern}
\usepackage{mathscinet}
\usepackage{mathtools}
\usepackage{microtype}
\usepackage{pmat}
\usepackage[automark]{scrpage2}

\usepackage{141k}

 \newtheorem{thm}{Theorem}[section]
 
 \newtheorem{lem}[thm]{Lemma}
 \newtheorem{prop}[thm]{Proposition}
\theoremstyle{definition}
 \newtheorem{defn}[thm]{Definition}
 
 \newtheorem*{prob}{Problem}
\theoremstyle{remark}
 \newtheorem{rem}[thm]{Remark}

\mathtoolsset{mathic=true}

\numberwithin{equation}{section}
\title{On a Simultaneous Approach to the Even
and Odd Truncated Matricial \Stieltjes{} Moment Problem}
\subtitle{II. An $\ug$-Schur-Stieltjes-type algorithm for sequences of holomorphic matrix-valued functions}
\author{Bernd Fritzsche \and Bernd Kirstein \and Conrad M\"adler}

\begin{document}
\maketitle

\begin{abstract}
 The main goal of this paper is to achieve a simultaneous treatment of the even and odd truncated matricial Stieltjes moment problems in the most general case. These results are generalizations of results of Chen and Hu~\zitas{MR1807884,MR2038751} which considered the particular case \(\ug=0\). Our approach is based on Schur analysis methods. More precisely, we use two interrelated versions of Schur-type algorithms, namely an algebraic one and a function-theoretic one. The algebraic version was worked out in a former paper of the authors. It is an algorithm which is applied to finite or infinite sequences of complex matrices. The construction and investigation of the function-theoretic version of our Schur-type algorithm is a central theme of this paper. This algorithm will be applied to relevant subclasses of holomorphic matrix-valued functions of the Stieltjes class. Using recent results on the holomorphicity of the Moore-Penrose inverse of matrix-valued Stieltjes functions, we obtain a complete description of the solution set of the moment problem under consideration in the most general situation.
\end{abstract}

\section{Introduction}
 This paper is a continuation of~\zita{114arxiv}. The main goal is to achieve a simultaneous treatment of the even and odd cases of a truncated matricial Stieltjes moment problem. We are guided by our investigations on the truncated matricial Hamburger moment problem in~\zita{MR3380267}, where we simultaneously treated the even and odd cases. Our approach in~\zita{MR3380267} is based on Schur analysis methods the origin of which goes back to the fundamental memoir~\zita{Nev22} by R.~Nevanlinna. We used two interrelated versions of Schur-type algorithms, namely an algebraic one and a function-theoretic one. The algebraic version had been worked out in the paper~\zita{MR3014199}, whereas the function-theoretic version and the interplay of both versions where investigated in~\zita{MR3380267}. Against to this background, we want to treat now simultaneously the even and odd cases for an analogous truncated matricial Stieltjes problem. More precisely, our approach will be based again on Schur analysis methods, namely on the interplay between two versions of Schur-type algorithms, an algebraic one and a function-theoretic one. The investigation of the algebraic version of our Schur-type algorithm was the main content of the paper~\zita{114arxiv}. The first central theme of this paper is to work out the function-theoretic version of our Schur-type algorithm. A careful analysis of the interplay between both versions will lead us finally to a complete solution of the original truncated matricial Stieltjes problem in the most general case. This is the main achievement of this paper.
 
 Though, roughly speaking, we repeat the basic strategy used in~\zitas{MR3014199,MR3380267} the concrete constructions  are rather different. What concerns the algebraic versions of our Schur-type algorithm, these differences where described in great detail in~\zita{114arxiv}. Roughly speaking, the main difference can be characterized as follows: In the Hamburger case the algebraic version of our Schur-type algorithm is a two-step algorithm, whereas in the Stieltjes case we are lead to a one-step algorithm. This observation reflects the nature of finite \tHnnd{} sequences from \(\Cqq\) on the one side and finite \taSnnd{} sequences from \(\Cqq\) on the other side. Using matricial versions of Hamburger-Nevanlinna theorems in~\zita{MR3380267}, we reformulated the truncated matricial Hamburger moment problem under study into an equivalent problem of prescribing an asymptotic expression for matrix-valued functions which are holomorphic in the open upper half plane. Guided by these experiences, we prove now Stieltjes-type versions of Hamburger-Nevanlinna-type theorems, which enable us to reformulate the original truncated matricial Stieltjes moment problem under study into an equivalent problem of prescribing an asymptotic expansion for a special class of matrix-valued functions which are holomorphic in \(\Cs\) where \(\ug\) is a given real number. More precisely, this special class is the class \(\SFqa\) of all \tlSFq{}, which is introduced in \rdefn{D1242}. We frequently apply results from our former paper~\zita{142}, which contains a detailed treatment of the class \(\SFqa\) and some important subclasses. The construction and investigation of a Schur-type algorithm for \tqqa{matrix-valued} functions holomorphic in \(\Cs\) which preserves special subclasses of \(\SFqa\) is a central theme of this paper. The idea how to build the algorithm in the matrix case was inspired by some constructions in the papers~\zitas{MR2038751,MR1807884} by Chen and Hu. An essential point of our approach is an intensive use of the interplay between the function-theoretic and algebraic versions of our matricial Schur-type algorithms. Both algorithms are formulated in terms of Moore-Penrose inverses of matrices. What concerns the function-theoretic version, it can be said that its effectiveness is mostly caused by recent results from~\zitas{142} on the holomorphicity of the Moore-Penrose inverse of special classes of holomorphic matrix-valued functions.

 In order to describe more concretely the central topics studied in this paper, we introduce some notation. Throughout this paper, let $p$ and $q$ be positive integers. Let \symba{\N}{n}, \symba{\NO}{n}, \symba{\Z}{z}, \symba{\R}{r}, and  \symba{\C}{c} be the set of all positive integers, the set of all \tnn{} integers, the set of all integers, the set of all real numbers, and the set of all complex numbers, respectively. For every choice of $\rho, \kappa \in \R \cup\{-\infty,\infp\}$, let $\mn{\rho}{\kappa}\defeq \setaa{k\in\Z}{\rho\leq k \leq \kappa}$\index{z@$\mn{\rho}{\kappa}$}. We will write \symba{\Cpq}{c}, \symba{\CHq}{c}, \symba{\Cggq}{c}, and \symba{\Cgq}{c} for the set of all complex \tpqa{matrices}, the set of all \tH{} complex \tqqa{matrices}, the set of all \tnnH{} complex \tqqa{matrices}, and the set of all \tpH{} complex \tqqa{matrices}, respectively. 

 We will use \symba{\BorR}{b} to denote the $\sigma$\nobreakdash-algebra of all Borel subsets of $\R$. Let $\Omega \in \BorR\setminus\set{\emptyset}$. Then, let $\BorO  \defeq \BorR\cap \Omega$\index{b@$\BorO $}. Furthermore, we will write $\MggqO$ to designate the set of all \tnnH{} \tqqa{measures} defined on $\BorO$, \ie{}, the set of all $\sigma$\nobreakdash-additive mappings $\mu\colon\BorO\to \Cggq$. We will use the integration theory with respect to \tnnH{} \tqqa{measures}, which was worked out independently by I.~S.~Kats~\zita{MR0080280} and M.~Rosenberg~\zita{MR0163346}. Some features of this theory are sketched in \rapp{A1608}. For all $\kappa\in \NOinf$, we will use \symba{\MggquO{\kappa}}{m} to denote the set of all $\sigma \in \MggqO$ such that the integral\index{s@$\suo{j}{\sigma}$}
 \begin{equation}\label{F1*1}
  \suo{j}{\sigma}
  \defeq \int_\Omega x^j \sigma(\dif x)
 \end{equation}
 exists for all \(j\in\mn{0}{\kappa}\).

\breml{R1*1}
 If $k,\ell\in\NO$ with $k<\ell$, then the inclusion $\MggquO{\ell}\subseteq\MggquO{k}$ holds true.
\erem

 The central problem in this paper is the truncated version of the following power moment problem of Stieltjes-type\index{m@\mproblem{\rhl}{\kappa}{=}}
\begin{prob}[\mproblem{\rhl}{\kappa}{=}]
 Let \(\ug\in\R\), let $\kappa\in \NOinf $ and let $\seqska $ be a sequence of complex \tqqa{matrices}. Describe the set $\MggqKksg$\index{m@$\MggqKksg$} of all $\sigma \in \MggquK{\kappa}$ for which  $\suo{j}{\sigma}=\su{j}$ is fulfilled for all $j \in \mn{0}{\kappa}$.
\end{prob}
 We note that there is a further matricial version of the truncated Stieltjes moment problem:\index{m@\mproblem{\rhl}{m}{\leq}}
\begin{prob}[\mproblem{\rhl}{m}{\leq}]
 Let \(\ug\in\R\), let $m\in\NO$ and let $\seqs{m}$ be a sequence of complex \tqqa{matrices}. Describe the set $\MggqKskg{m}$\index{m@$\MggqKskg{m}$} of all $\sigma\in\MggquK{m}$ for which  $\suo{j}{\sigma}=\su{j}$ is satisfied for each $j\in\mn{0}{m-1}$, whereas the matrix $\su{m}-\suo{m}{\sigma}$ is \tnnH{}.
\end{prob}
 A detailed treatment of the history of these two moment problems is contained in the introduction to the paper~\zita{MR2570113}.

 In order to state a necessary and sufficient condition for the solvability of each of the above formulated moment problems, we have to recall the notion of two types of sequences of matrices. If $n\in \NO$ and if $\seqs{2n}$ is a sequence of complex \tqqa{matrices}, then $\seqs{2n}$ is called \notii{\tHnnd{}} if the block Hankel matrix\index{h@$\Hu{n}$}
\[
 \Hu{n}
 \defeq \matauuo{\su{j+k}}{j,k}{0}{n}
\]
 is \tnnH{}. Let $\seqsinf $ be a sequence of complex \tqqa{matrices}. Then $\seqsinf $ is called \notii{\tHnnd{}} if $\seqs{2n}$ is \tHnnd{} for all $n\in\NO$. For all $\kappa\in\NOinf $, we will write \symba{\Hggqu{2\kappa}}{h} for the set of all \tHnnd{} sequences $\seqs{2\kappa}$ of complex \tqqa{matrices}. Furthermore, for all $n \in \NO$, let \symba{\Hggequ{2n}}{h} be the set of all sequences $\seqs{2n}$ of complex \tqqa{matrices} for which there exist complex \tqqa{matrices} $\su{2n+1}$ and $\su{2n+2}$ such that $\seqs{2(n+1)} \in\Hggqu{2(n+1)}$, whereas \symba{\Hggequ{2n+1}}{h} stands for the set of all sequences $\seqs{2n+1}$ of complex \tqqa{matrices} for which there exist some $\su{2n+2}\in\Cqq$ such that $\seqs{2(n+1)} \in\Hggqu{2(n+1)}$. For each $m\in\NO$, the elements of the set $\Hggequ{m}$ are called \emph{\tHnnde{} sequences}. For technical reason, we set $\Hggeqinf\defeq\Hggqinf$\index{h@$\Hggeqinf$}.
 
 Besides the just introduced classes of sequences of complex \tqqa{matrices} we need analogous classes of sequences of complex \tqqa{matrices} which take into account the influence of the prescribed number \(\ug\in\R\): Let \(\seqs{\kappa}\) be a sequence of complex \tpqa{matrices}. Then, for all \(n\in\NO\) with \(2n+1\leq\kappa\), we introduce the block Hankel matrix\index{k@$\Ku{n}$}
\[
 \Ku{n}
 \defeq\matauuo{\su{j+k+1}}{j,k}{0}{n}.
\]
 Let \(\ug\in\R\). Now we will introduce several classes of finite or infinite sequences of complex \tqqa{matrices} which are characterized by the sequences \(\seqs{\kappa}\) and \(\seq{-\ug\su{j}+\su{j+1}}{j}{0}{\kappa-1}\). Let \(\Kggqua{0}\defeq\Hggqu{0}\)\index{k@$\Kggqua{0}$}, and, for all \(n\in\N\), let \symba{\Kggqua{2n}}{k} be the set of all sequences \(\seqs{2n}\) of complex \tqqa{matrices} for which the block Hankel  matrices \(\Hu{n}\) and \(-\ug\Hu{n-1}+\Ku{n-1}\) are both \tnnH{}, \ie{},
\beql{Kgg2n}
 \Kggqua{2n}
 =\setaa*{\seqs{2n}\in\Hggqu{2n}}{\seq{-\ug\su{j}+\su{j+1}}{j}{0}{2(n-1)}\in\Hggqu{2(n-1)}}.
\eeq
 Furthermore, for all \(n\in\NO\), let \symb{\Kggqua{2n+1}} be the set of all sequences \(\seqs{2n+1}\) of complex \tqqa{matrices} for which the block Hankel  matrices \(\Hu{n}\) and \(-\ug\Hu{n}+\Ku{n}\) are both \tnnH{}, \ie{},
\beql{Kgg2n+1}
 \Kggqua{2n+1}
 \defeq\setaa*{\seqs{2n+1}\text{ from }\Cqq}{\set*{\seqs{2n},\seq{-\ug\su{j}+\su{j+1}}{j}{0}{2n}}\subseteq\Hggqu{2n}}.
\eeq
 Let \symba{\Kggqinfa}{k} be the set of all sequences \(\seqsinf\) of complex \tqqa{matrices} such that \(\seqs{m}\in\Kggqua{m}\) for all \(m\in\NO\). Formulas \eqref{Kgg2n} and \eqref{Kgg2n+1} show that the sets \(\Kggqua{2n}\) and \(\Kggqua{2n+1}\) are determined by two conditions. The condition \(\seqs{2n}\in\Hggqu{2n}\) ensures that a particular Hamburger  moment problem associated with the sequence \(\seqs{2n}\) is solvable (see, \eg{}~\zitaa{MR2570113}{\cthm{4.16}}). The second condition \(\seq{-\ug\su{j}+\su{j+1}}{j}{0}{2(n-1)}\in\Hggqu{2(n-1)}\) (resp.\ \(\seq{-\ug\su{j}+\su{j+1}}{j}{0}{2n}\in\Hggqu{2n}\)) controls that the original sequences \(\seqs{2n}\) and \(\seqs{2n+1}\) are well adapted to the interval \(\rhl \). Let \(m\in\NO\). Then, let \symba{\Kggequa{m}}{k} be the set of all sequences \(\seqs{m}\) of complex \tqqa{matrices} for which there exists a complex \tqqa{matrix} \(\su{m+1}\) such that \(\seqs{m+1}\) belongs to \(\Kggqua{m+1}\). Further, let \(\Kggeqinfa\defeq\Kggqinfa\)\index{k@$\Kggeqinfa$}. We call a sequence \(\seqska\) of complex \tqqa{matrices} \notii{\taSnnd{}} if it belongs to \(\Kggqka\) and \notii{\taSnnde{}} if it belongs to \(\Kggeqka\). 
 
 Now we can characterize the situations in which the mentioned problems have a solution:
 
\bthmnl{\zitaa{MR2735313}{\cthmss{1.3}{1.4}}}{T1*3+2}
 Let \(\ug\in\R\), let \(m\in\NO\), and let \(\seqs{m}\) be a sequence of complex \tqqa{matrices}. Then:
 \benui
  \il{T1*3+2.a} \(\Mggqaag{\rhl }{\seqs{m}}\neq\emptyset\) if and only if \(\seqs{m}\in\Kggequa{m}\).
  \il{T1*3+2.b} \(\Mggqaakg{\rhl }{\seqs{m}}\neq\emptyset\) if and only if \(\seqs{m}\in\Kggqua{m}\).
 \eenui
\ethm

%

 In the case \(\ug=0\) there is some important work of Bolotnikov~\zita{MR1362524} and Chen/Hu~\zitas{MR1807884,MR2038751} which influenced our investigations. This will be explained now in more detail. \rthmp{T1*3+2}{T1*3+2.a} was proved in Bolotnikov~\zitaa{MR2038751}{\clem{1.7}}. Bolotnikov also stated a parametrization of the solution set of Problem \mproblem{[0,\infp)}{m}{\leq} in the language of the Stieltjes transforms (see~\zitaa{MR2038751}{\cthmss{5.3}{6.6}}. Chen and Hu stated a parametrization of the solution set of Problem \mproblem{[0,\infp)}{m}{=} in the language of the Stieltjes transforms (see~\zitaa{MR2038751}{\cthm{4.1}}).
 
 The main goal of our paper here is to give via Stieltjes transformation a parametrization of the solution set of Problem \mproblem{\rhl}{m}{=} (see \rthmss{T12*1}{T12*7}).

\section{The class $\SFqa$}\label{S*2}
 The use of several classes of holomorphic matrix-valued functions is one of the special features of this paper. In this section, we summarize some basic facts about the class of \tlSF{s} of order \(q\), which are mostly taken from our former paper~\zita{142}. If $A\in\Cqq$, then let $\re A\defeq\frac{1}{2}(A+A^\ad)$\index{r@$\re A$} and $\im A\defeq\frac{1}{2\iu}(A-A^\ad)$\index{i@$\im A$} be the real part and the imaginary part of $A$, respectively. Let \(\ohe\defeq\setaa{z\in\C}{\im z>0}\)\index{p@\(\ohe\)} be the open upper half plane of \(\C\). The first class of functions which plays an essential role in this paper, is the  following.

 
\bdefnl{D1242}
 Let \(\ug\in\R\) and let \(F\colon\Cs\to\Cqq\). Then \(F\) is called a \notii{\tlSFq{}} if \(F\) satisfies the following three conditions:
 \bAeqi{0}
  \il{D1242.I} \(F\) is holomorphic in \(\Cs\).
  \il{D1242.II} For all \(w\in\ohe\), the matrix \(\im\ek{F(w)}\) is \tnnH{}.
  \il{D1242.III} For all \(w\in\lhl\), the matrix \(F(w)\) is \tnnH{}.
 \eAeqi
 We denote by \(\SFqa\)\index{s@\(\SFqa\)} the set of all \tlSF{s} of order \(q\).
\edefn
 For a comprehensive survey on the class $\SFqa$, we refer the reader to~\zita{142}. The functions belonging to the class \(\SFqa\) admit an important integral representation. To state this, we introduce some terminology: If $\mu$ is a \tnnH{} \tqqa{measure} on a measurable space $(\Omega,\mathfrak{A})$, then we will use \symba{\LaaaK{\Omega}{\mathfrak{A}}{\mu}}{l} with \(\mathbb{K}\in\set{\R,\C}\) to denote the space of all Borel-measurable functions $f\colon\Omega\to\mathbb{K}$ for which the integral $\int_\Omega f\dif\mu$ exists. In preparing the desired integral representation, we observe that, for all $\mu\in\MggqK$ and each $z\in\Cs$, the function $\hu{z}\colon\rhl\to\C$\index{h@$\hu{z}$} defined by $\hua{z}{t}\defeq(1+t-\ug)/(t-z)$ belongs to $\LaaaC{\rhl}{\BorK}{\mu}$.

\bthmnl{\cf{}~\zitaa{142}{\cthm{3.6}}}{T2*1}
 Let \(\ug \in\R\) and let \(F\colon\Cs  \to \Cqq \). Then:
 \benui
 \il{T2*1.a} If \(F \in \SFqa \), then there are a unique matrix \(\gamma \in \Cggq\) and a unique \tnnH{} measure \(\mu \in\Mggqrhl \) such that
 \beql{F2*2}
  F(z)
  = \gamma +\int_{\rhl} \frac{1+ t-\ug}{t-z} \mu (\dif t)
 \eeq
 holds true for each \(z\in\Cs \).
 \il{T2*1.b} If there are a matrix \(\gamma\in\Cggq\) and a \tnnH{} measure \(\mu \in\Mggqrhl \) such that \(F\) can be represented via \eqref{F2*2} for each \(z\in \Cs \), then \(F\) belongs to the class \(\SFqa \).
 \eenui
\ethm

 For all $F\in\SFqa$, we will write \symba{(\gammau{F},\muu{F})}{g} for the unique pair \((\gamma,\mu)\in\Cggq\times\Mggqrhl \) for which the representation \eqref{F2*2} holds true for all \(z\in \Cs \).

In the special case that \(q=1\) and \(\ug=0\) \rthm{T2*1} can be found in Krein/Nudelman~\zitaa{MR0458081}{Appendix}.

\bremnl{\zitaa{142}{\crem{3.12}}}{R0832}
 Let \(\ug\in\R\) and let \(F\in\SFqa\). Further let \(A\in\Cqq\) be such that \(\gammaF+A\) is \tnnH{}. Then \(G\defeq F+A\in\SFqa\) and \((\gammaG,\muG)=(\gammaF+A,\muF)\).
\erem

\bpropnl{\zitaa{142}{\cprop{3.13}}}{P2*5}
 Let \(\ug\in\R\) and let \(F\in\SFqa\). Then \(\lim_{y\to\infp}F(\iu y)=\gammau{F}\).
\eprop

 We are particularly interested in the structure of the values of functions belonging to $\SFqa$. For this reason, we introduce some notation: For all $A\in\Cpq$, let \symba{\nul{A}}{n} be the null space of $A$ and \symba{\ran{A}}{r} be the column space of $A$.

\bpropnl{\cf{}~\zitaa{142}{\cprop{3.15}}}{P2*7}
 Let \(\ug \in\R\) and let \(F\in\SFqa \). For all \(z\in \Cs \), then
 \begin{align*}
  \Ran{F(z)}&=\ran{\gammau{F}}+\Ran{\muuA{F}{\rhl }}&
  &\text{and}&
  \Nul{F(z)}&=\nul{\gammau{F}}\cap\Nul{\muuA{F}{\rhl }}.
 \end{align*}
\eprop
 \rprop{P2*7} contains essential information on the class $\SFqa$. It indicates that, for an arbitrary function $F$ belonging to $\SFqa$, the null space $\nul{F(z)}$ and the column space $\ran{F(z)}$ are independent of the concrete point $z\in\Cs$ and, furthermore, in which way these linear subspaces of $\Cq$ are determined by the pair $(\gammaF,\muF)$ of $F$.

 In the sequel, we will sometimes meet situations where interrelations of the null space (resp.\ column space) of a function $F\in\SFqa$ to the null space (resp.\ column space) of a given matrix $A\in\Cpq$ are of interest. More precisely, we will frequently apply the following auxilliary result. In preparing its formulation we note that for a matrix \(A\in\Cpq\) we denote by \symb{A^\mpi} its Moore-Penrose inverse. This means \(A^\mpi\) is the unique matrix \(X\) from \(\Cqp\) which satisfies the four equations \(AXA=A\), \(XAX=X\), \((AX)^\ad=AX\) and \((XA)^\ad=XA\).

\blemnl{\cf{}~\zitaa{142}{\clem{3.18}}}{L2*8}
 Let \(\ug \in\R\), let \(A\in\Cpq\), and let \(F\in\SFqa\). Then the following statements are equivalent:
 \baeqi{0}
  \il{L2*8.i} \(\nul{A}\subseteq\nul{F(z)}\) for all \(z\in\Cs\).
  \il{L2*8.ii} There is a \(z_0\in\Cs\) such that \(\nul{A}\subseteq\nul{F(z_0)}\).
  \il{L2*8.iii} \(\nul{A}\subseteq\nul{\gammau{F}}\cap\nul{\muua{F}{\rhl }}\).
  \il{L2*8.iv} \(FA^\mpi A=F\).
  \il{L2*8.v} \(\ran{F(z)}\subseteq\ran{A^\ad}\) for all \(z\in\Cs\).
  \il{L2*8.vi} There is a \(z_0\in\Cs\) such that \(\ran{F(z_0)}\subseteq\ran{A^\ad}\).
  \il{L2*8.vii} \(\ran{\gammau{F}}+\ran{\muua{F}{\rhl}}\subseteq\ran{A^\ad}\).
  \il{L2*8.viii} \(A^\mpi AF=F\).
 \eaeqi
\elem

 A generic application of \rlem{L2*8} will be concerned with situations where the matrix $A$ even belongs to $\Cggq$.

 In our subsequent considerations we will very often use the Moore-Penrose inverse of functions belonging to the class $\SFqa$. In this connection, the following result turns out to be of central importance.

\bpropnl{\zitaa{142}{\cthm{6.3}}}{P2*9}
 Let \(\ug\in\R\) and let \(F\in\SFqa\). Then \(G\colon\Cs\to\Cqq\) defined by \(G(z)\defeq-(z-\ug)^\inv\ek{F(z)}^\mpi\) belongs to the class \(\SFqa\) as well.
\eprop

\section{On some subclasses of $\SFqa$}\label{S*3}
 An essential feature of our subsequent considerations is the use of several subclasses of $\SFqa$. In this section, we summarize basic facts about these subclasses under the special orientation of this paper. These subclasses are characterized by growth properties on the positive imaginary axis. It should be mentioned that the scalar versions of the function classes were introduced and studied in the paper~\zita{KK74}. Next we recall the Euclidean norm of a matrix. For each \(B\in\Cqq\), let \symba{\tr B}{t} be the trace of \(B\). If $A\in \Cpq$, then \(\normE{A}\defeq\sqrt{\tr\rk{A^\ad A}}\)\index{\(\normE{A}\)} is the Euclidean norm of \(A\). An important subclass of the class $\SFqa$ is the set
 \beql{F3*3}
  \SFuqa{0}
  \defeq\setaa*{F\in\SFqa}{\sup_{y\in[1,\infp)}y\normE*{F(\iu y)}<\infp}.
 \eeq
 Now we want to characterize the class \(\SFuqa{0}\) via an integral representation.
 
\begin{lem}[\zitaa{142}{\clem{A.4}}]\label{63A}
 Let \(\Omega\) be a \tne{} closed subset of \(\R\) and let \(\sigma\in\MggqO \). Then:
 \benui
  \il{63A.a} For each \(w\in\C\setminus \Omega\), the function \(g_w \colon\Omega\to\C\)\index{g@\(g_w\)} defined by \(g_w (t)\defeq1/(t-w)\) belongs to \(\LaaaC{\Omega}{\BorO}{\sigma}\).
  \il{63A.b} The matrix-valued function \(S\colon\C\setminus \Omega\to\Cqq \) given by \(S(w)\defeq\int_\Omega g_w \dif\sigma\) satisfies \(\ran{S(z)}=\ran{\sigma (\Omega)}\) and in particular \(\rank S(z) =\rank\sigma (\Omega)\) for each  \(z\in\C\setminus [\inf \Omega, \sup \Omega]\).
 \il{63A.c} \(-\iu\lim_{y\to\infp}yS(\iu y)=\sigma(\Omega)\).
 \eenui
\end{lem}
 Let \(\Omega\) be a \tne{} closed subset of \(\R\) and let \(\sigma\in\Mggqa{\Omega}\). Then, in view of \rlemp{63A}{63A.a}, the matrix-valued function \(\OSt{\sigma}\colon\C\setminus\Omega\to\Cqq\)\index{s@\(\OSt{\sigma}\)} given by
\beql{F3*2}
 \OSta{\sigma}{z}
 \defeq\int_\Omega\frac{1}{t-z}\sigma(\dif t)
\eeq
 is well-defined and called \notii{\tOSto{\Omega}{\sigma}}.
 
 Let \(z\in\C\setminus\Omega\). Then obviously \(\ko z\in\C\setminus\Omega\) and from \eqref{F3*2} it follows
\beql{G0836}
 \OSta{\sigma}{\ko z}
 =\int_\Omega\frac{1}{t-\ko z}\sigma(\dif t)
 =\int_\Omega\ko{\rk*{\frac{1}{t-z}}}\sigma(\dif t)
 =\ek*{\int_\Omega\frac{1}{t-z}\sigma(\dif t)}^\ad
 =\ek*{\OSta{\sigma}{z}}^\ad.
\eeq
 
 Now we consider the case \(\Omega=\rhl\). We want to characterize the set of all \taSt{s} of measures belonging to $\MggqK$.

\bthml{T3*2}
 The mapping $\sigma\mapsto\OSt{\sigma}$ is a bijective correspondence between $\MggqK$ and $\SFuqa{0}$. In particular $\SFuqa{0}=\setaa{\OSt{\sigma}}{\sigma\in\MggqK}$.
\ethm
\bproof
 This is an immediate consequence of~\zitaa{142}{\cthm{5.1}}.
\eproof

 For each $F\in\SFuqa{0}$, the unique measure $\sigma\in\MggqK$ satisfying $\OSt{\sigma}=F$ is called the \notii{\taSmo{F}}. We will also write \symba{\OSm{F}}{s} for \(\sigma\). \rthm{T3*2} indicates that the \taSt{} $\OSt{\sigma}$ of a measure $\sigma\in\MggqK$ is characterized by a particular mild growth on the positive imaginary axis.

\bcorl{C0838}
 Let \(\ug\in\R\) and let \(F\in\SFuqa{0}\). For \(z\in\Cs\) then \(\ko z\in\Cs\) and \(F(\ko z)=\ek{F(z)}^\ad\).
\ecor
\bproof
 Combine \rthm{T3*2} and \eqref{G0836}
\eproof

\bpropnl{\zitaa{142}{\cprop{6.4}}}{P0844}
 Let \(\ug\in\R\), let \(F\in\SFuqa{0}\), and let \(\OSm{F}\) be the \taSmo{F}. Then \(G\colon\Cs \to\Cqq\) defined by \(G(z)\defeq-(z-\ug )^\inv\ek{F(z)}^\mpi\) belongs to \(\SFqa\) and
 \[
  \gammaG
  =\ek*{\OSmA{F}{\rhl}}^\mpi.
 \]
 In particular, if \(F\) is not the constant function with value \(\Oqq\), then \(G\in\SFqa\setminus\SFuqa{0}\).
\eprop

 In view of \rthm{T3*2}, \rprob{\mproblem{\rhl}{\kappa}{=}} can be given a first reformulation as an equivalent problem in the class $\SFuqa{0}$ as follows:

\index{s@\iproblem{\rhl}{\kappa}{=}}
\begin{prob}[\iproblem{\rhl}{\kappa}{=}]
 Let \(\ug\in\R\), let \(\kappa\in\NOinf\), and let \(\seqska\) be a sequence of complex \tqqa{matrices}. Describe the set \symba{\SFqaska}{s} of all \(F\in\SFuqa{0}\) the \taSm{} of which belongs to \(\MggqKksg\).
\end{prob}

 In \rsec{S*6}, we will state a reformulation of the original power moment problem \mproblem{\rhl}{\kappa}{=} as an equivalent problem of finding a prescribed asymptotic expansion in a sector of the open upper half plane $\ohe$. Furthermore, we will see that a detailed analysis of the behavior on the positive imaginary axis of the concrete functions of $F\in\SFqa$ under study is extremely useful. For this reason, we turn now our attention to some subclasses of $\SFqa$ which are described in terms of their growth on the positive imaginary axis. First we consider the set\index{s@$\SFdqa$}
\beql{F3*4}
 \SFdqa
 \defeq\setaa*{F\in\SFqa}{\lim_{y\to\infp}\normE*{F(\iu y)}=0}.
\eeq
 In the following considerations we associate with a function \(F\in\SFqa\) often the unique ordered pair \((\gammaF,\muF)\) given via \rthmp{T2*1}{T2*1.a}.

\bremnl{\zitaa{142}{\ccor{3.14}}}{R3*5}
 Let \(\ug\in\R\). Then
 \[
  \SFdqa
  =\setaa*{F\in\SFqa}{\gammaF=\Oqq}.
 \]
\erem

\breml{R3*6}
 Let \(\ug\in\R\), let $n\in\N$, and let $\seq{p_k}{k}{1}{n}$ be a sequence of positive integers. For all $k\in\mn{1}{n} $, let $F_k\in\SFdua{p_k}$ and let $A_k\in\Coo{p_k}{q}$. In view of \rrem{R3*5} and~\zitaa{142}{\crem{3.11}}, then $\sum_{k=1}^nA_k^\ad F_kA_k\in\SFdqa$.
\erem

 Now we state modifications of \rprop{P2*7} and \rlem{L2*8} for $\SFdqa$.

\bpropl{P3*7}
 Let \(\ug\in\R\), let $F\in\SFdqa$, and let \(\muF\) be given via \rthmp{T2*1}{T2*1.a}. For all $z\in\Cs$, then
 \begin{align*}
  \Ran{F(z)}&=\Ran{\muFA{\rhl}}&
  &\text{and}&
  \Nul{F(z)}&=\Nul{\muFA{\rhl}}.
 \end{align*}
\eprop
\bproof
 Combine \rprop{P2*7} and \rrem{R3*5}.
\eproof

\bleml{L3*8}
 Let $A\in\Cpq$ and let $F\in\SFdqa$. Then the statements
 \baeqi{8}
  \il{L3*8.ix} $\ran{\muFa{\rhl}}\subseteq\ran{A^\ad}$.
 \eaeqi
 and
 \baeqi{9}
  \item $\nul{A}\subseteq\nul{\muFa{\rhl}}$.
 \eaeqi
 are equivalent. Furthermore,~\ref{L3*8.ix} is equivalent to each of the statements~\ref{L2*8.i}--\ref{L2*8.viii} stated in \rlem{L2*8}.
\elem
\bproof
 Combine \rrem{R3*5} and \rlem{L2*8}.
\eproof

 For all \(\ug\in\R\) and all $\kappa\in\Ninf$, we now consider the class\index{s@$\SFkaqa$}
\beql{F3*11}
 \SFkaqa
 \defeq\setaa*{F\in\SFuqa{0}}{\OSm{F}\in\MggquK{\kappa}}.
\eeq

\bleml{L3*22}
 Let \(\ug\in\R\), let $\kappa\in\NOinf$, and let \(F\in\SFkaqa\) with \taSm{} \(\OSm{F}\). For each \(z\in\Cs\), then
 \begin{align*}
  \Ran{F(z)}&=\Ran{\OSmA{F}{\rhl}}&
 &\text{and}&
  \Nul{F(z)}& =\Nul{\OSmA{F}{\rhl}}.
 \end{align*}
\elem
\bproof
 This is an immediate consequence of~\zitaa{142}{\cprop{5.3}}.
\eproof

\breml{R3*23}
 Let \(\ug\in\R\),  let \(\kappa\in\NOinf\), and let \(F\in\SFkaqa\). Then from \eqref{F3*11} and \eqref{F3*3} we see that \(\lim_{y\to\infp}F(\iu y)=\Oqq\).
\erem

\breml{R3*18}
 Let \(\ug\in\R\). In view of the \rremss{R1*1}{R3*23}, then
 \beql{Nr.TMB}
 \SFinfqa
 \subseteq\SFuqa{\ell}
 \subseteq\SFuqa{m}
 \subseteq\SFuqa{0}
 \subseteq\SFdqa
 \subseteq\SFqa
 \eeq
 for all \(\ell,m\in\NO\) with \(\ell\geq m\).
\erem

\section{On the class $\SFdqaa{A}$}\label{S*4}
 In this section, we consider a particular subclass of the class $\SFdqa$, which was introduced in \eqref{F3*4}. We have seen in \rprop{P2*7} that, for an arbitrary function $F\in\SFqa$, the null space of $F(z)$ is independent from the concrete choice of $z\in\Cs$. For the case $F\in\SFdqa$, a complete description of this constant null space was given in \rprop{P3*7}. Against to this background, we single out now a special subclass of $\SFdqa$, which is characterized by the interrelation of the constant null space to the null space of a prescribed matrix $A\in\Cpq$. More precisely, for all $A\in\Cpq$, let\index{s@\(\SFdqaa{A}\)}
\begin{equation}\label{F4*1}
  \SFdqaa{A}
  \defeq\setaa*{F\in\SFdqa}{\nul{A}\subseteq\Nul{\muuA{F}{\rhl}}} 
 \end{equation}
 In our later investigations the role of the matrix $A$ will be taken by matrices which are generated from the sequence of data matrices of the problem under consideration via a Schur-type algorithm.

\breml{R4*1}
 Let \(\ug\in\R\). If $A\in\Cpq$ satisfies $\nul{A}=\set{\Ouu{q}{1}}$, then $\SFdqaa{A}=\SFdqa$. In particular, this situation arises in the case that $p=q$ and $\det A\neq0$ are fulfilled.
\erem

 In the following considerations, we associate with a function \(F\in\SFqa\) often the unique ordered pair \((\gammaF,\muF)\) given via \rthmp{T2*1}{T2*1.a}.
 
\bexaml{E4*2}
 Let \(\ug\in\R\), let \(A\in\Cpq\), and let \(F\colon\Cs\to\Cqq\) be defined by \(F(z)\defeq\Oqq\). In view of~\zitaa{142}{\cexam{3.19}}, then one can easily see that \(F\in\SFdqaa{A}\) and \(\muu{F}\) is the zero measure belonging to \(\MggqK\).
\eexam

\bleml{L4*3}
 Let \(\ug\in\R\), let $A\in\Cpq$ and $F\colon\Cs\to\Cqq$. Then the following statements are equivalent:
 \baeqi{0}
  \item $F\in\SFdqaa{A}$.
  \item $F\in\SFdqa $ and $\nul{A}\subseteq\nul{F(z)}$ for all $z\in\Cs$.
  \item $F\in\SFqa$, $\lim_{y\to\infp}\normE{F(\iu y)}=0$, and $\nul{A}\subseteq\nul{F(z)}$ for all $z\in\Cs$.
  \item $F\in\SFqa$, $\gammaF=\Oqq$, and $\nul{A}\subseteq\nul{\muFa{\rhl}}$.
 \eaeqi
\elem
\bproof
 This follows from \eqref{F4*1}, \rlem{L3*8}, \eqref{F3*4}, and \rrem{R3*5}.
\eproof

\breml{R-C4*7}
 Let \(\ug\in\R\), let \(A\in\Cpq\), and let \(F\in\SFdqa\). In view of \rrem{R3*5} and \rlem{L2*8}, then the following statements are equivalent:
 \baeqi{0}
  \item \(F\in\SFdqaa{A}\).
  \item There is a \(z_0\in\Cs\) such that \(\nul{A}\subseteq\nul{F(z_0)}\).
  \item \(\nul{A}\subseteq\nul{\muua{F}{\rhl }}\).
  \item \(FA^\mpi A=F\).
  \item \(\ran{F(z)}\subseteq\ran{A^\ad}\) for all \(z\in\Cs\).
  \item There is a \(z_0\in\Cs\) such that \(\ran{F(z_0)}\subseteq\ran{A^\ad}\).
  \item \(\ran{\muua{F}{\rhl}}\subseteq\ran{A^\ad}\).
  \item \(A^\mpi AF=F\).
 \eaeqi
\erem

\breml{R4*5}
 Let \(\ug\in\R\), let \(r\in\N\), and let \(A\in\Cpq\) and \(B\in\Coo{r}{q}\) be such that \(\nul{A}\subseteq\nul{B}\). Then \(\SFdqaa{B}\subseteq\SFdqaa{A}\).
\erem

\breml{R-P4*6}
 Let \(\ug\in\R\) and let \(A\in\Cpq\). In view of \rrem{R-C4*7}, \(A^\mpi A=(A^\mpi A)^\ad\) and~\zitaa{142}{\crem{3.20}}, then \(\SFdqaa{A}=\set{A^\mpi AFA^\mpi A\colon F\in\SFdqa}\).
\erem

The following result contains essential information on the structure of the set $\SFdqaa{A}$, where $A\in\Cpq$.

\bpropl{P4*9}
 Let \(\ug\in\R\) and let $A\in\Cpq$. Then:
 \benui
  \il{P4*9.a} If $A=\Opq$, then $\SFdqaa{A}=\set{F}$, where $F\colon\Cs\to\Cqq$ is defined by $F(z)\defeq\Oqq$.
  \il{P4*9.b} Suppose that $r\defeq\rank A$ fulfills $r\geq1$. Let $u_1,u_2,\dotsc,u_r$ be an orthonormal basis of $\ran{A^\ad}$ and let $U\defeq\mat{ u_1,u_2,\dotsc,u_r}$. Then:
  \begin{enuii}
   \il{P4*9.b1} \(\SFdqaa{A}=\setaa{UfU^\ad}{f\in\SFdua{r}}\).
   \il{P4*9.b2} If $f,g\in\SFdua{r}$ are such that $UfU^\ad=UgU^\ad$, then $f=g$.
  \end{enuii}
 \eenui
\eprop
\bproof
 \eqref{P4*9.a} This follows from \rrem{R-P4*6} and \rexam{E4*2}.
 
 \eqref{P4*9.b1} Let $G\in\SFdqaa{A}$. In view of \rrem{R-P4*6}, there exists an $F\in\SFdqa$ such that $G=A^\mpi AFA^\mpi A$. Let $f\defeq U^\ad FU$. Because of \rrem{R3*6}, then $f\in\SFdua{r}$. In view of \rrem{RA*6}, we have $UU^\ad=A^\mpi A$. Thus, $G=UU^\ad FUU^\ad=UfU^\ad$. Hence, \(\SFdqaa{A}\subseteq\setaa{UfU^\ad}{f\in\SFdua{r}}\).
 
 Conversely, let $f\in\SFdua{r}$. In view of \rrem{R3*6}, then \(UfU^\ad\in\SFdqa\). Now we consider an arbitrary $x\in\nul{A}$. In view of the construction of $U$ and the relation $[\nul{A}]^\oc=\ran{A^\ad}$, we get $U^\ad x=\Ouu{r}{1}$. Thus, $x\in\nul{U^\ad}$. Consequently, for each $z\in\Cs$ we get $\nul{A}\subseteq\nul{(UfU^\ad)(z)}$. The application of \rlem{L4*3} yields now \(\setaa{UfU^\ad}{f\in\SFdua{r}}\subseteq\SFdqaa{A}\). This completes the proof of~\eqref{P4*9.b1}.
 
 \eqref{P4*9.b2} In view of \rrem{RA*6}, we have $U^\ad U=\Iu{r}$. Thus $UfU^\ad=UgU^\ad$ implies $f=g$.
\eproof

\section{The class $\SFqas{\kappa}$}\label{S*5}
 In this section, we consider particular subclasses of the class $\SFkaqa$, which was introduced in \eqref{F3*3} for $\kappa=0$ and in \eqref{F3*11} for $\kappa\in\Ninf $. In view of \rthm{T3*2}, for each function $F$ belonging to one of the classes $\SFkaqa$ with some $\kappa\in\NOinf $, we can consider the \taSm{} $\OSm{F}$ of \(F\). Now taking \rrem{R1*1} into account, we turn our attention to subclasses of functions $F\in\SFkaqa$ with prescribed first $\kappa+1$ power moments of the \taSm{} $\OSm{F}$.

 For all \(\ug\in\R\), all $\kappa \in \NOinf $, and each sequence $\seqska$ of complex \tqqa{matrices}, now we consider the class\index{s@$\SFqaska$}
 \begin{equation} \label{F5*3}
  \SFqaska
  \defeq\setaa*{F \in\SFkaqa}{\OSm{F}\in \MggqKSg{\kappa}}.
 \end{equation}
 
\breml{R5*3}
 Let \(\ug\in\R\), let $\kappa\in\NOinf $, and let $\seqska$ be a sequence of complex \tqqa{matrices}. If $\iota\in\NOinf $ with $\iota\leq\kappa$ then it is readily checked that
 \[
  \SFqaS{\kappa}
  \subseteq\SFqaS{\iota}.
 \]
 and
 \[
  \SFqaS{\kappa}
  =\bigcap_{m=0}^\kappa\SFqaS{m}.
 \]
\erem

\breml{R1228}
 Let \(\ug\in\R\), let \(\kappa\in\NOinf\), and let \(\seqska\) be a sequence of complex \tqqa{matrices}. Let \(F\in\SFqaska\) with \taSm{} \(\OSm{F}\). Then
 \[
  \su{0}
  =\suo{0}{\OSm{F}}
  =\int_\rhl t^0\OSma{F}{\dif t}
  =\OSmA{F}{\rhl}.
 \]
\erem

 Now we characterize those sequences for which the sets defined in \eqref{F5*3} are \tne{}.

\bthml{T-P5*4}
 Let \(\ug\in\R\), let \(\kappa\in\NOinf\), and let \(\seqska\) be a sequence of complex \tqqa{matrices}. Then \(\SFqas{\kappa}\neq \emptyset\) if and only if \(\seqska\in\Kggeqka\).
\ethm
\bproof
 Combine \eqref{F5*3} and \eqref{F3*4}  with \rthmss{T1*3+2}{T3*2}.
\eproof

Now we state a useful characterization of the set of functions given in \eqref{F5*3}.


\bthml{T-C5*6}
 Let \(\ug\in\R\), let $\kappa\in\NOinf$, and let $\seqska$ be a sequence of complex \tqqa{matrices}. In view of \eqref{F3*2} then
 \[
  \SFqas{\kappa}
  =\setaa*{\OSt{\sigma}}{\sigma\in\MggqKSg{\kappa}}.
 \]
\ethm
\bproof
 Combine \eqref{F5*3} with \rthm{T3*2}.
\eproof

 \rthm{T-C5*6} shows that $\SFqaska $ coincides with the solution set of Problem~\iproblem{\rhl}{\kappa}{=}, which is via \taSt{} equivalent to the original Problem~\mproblem{\rhl}{\kappa}{=}. Thus, the investigation of the set $\SFqaska $ is a central theme of our further considerations. The next result contains essential information on the functions belonging to this set.

\bpropl{P5*7}
 Let $\kappa\in\NOinf $, let $\seqska \in\Kggeqka$, and let $F\in\SFqas{\kappa}$. Then:
 \benui
  \il{P5*7.a} $\ran{F(z)}=\ran{\su{0}}$ and $\nul{F(z)}=\nul{\su{0}}$ for all $z\in\Cs$.
  \il{P5*7.b} $[F(z)][F(z)]^\mpi=\su{0}\su{0}^\mpi$ and $[F(z)]^\mpi[F(z)]=\su{0}^\mpi\su{0}$ for all $z\in\Cs$.
  \il{P5*7.c} The function $F$ belongs to the class $\SFqa$ with 
  \begin{align}\label{P5*7.A}
   \ran{\su{0}}&=\ran{\gammaF}+\Ran{\muFA{\rhl}},&
   \nul{\su{0}}&=\nul{\gammaF}\cap\Nul{\muFA{\rhl}}.
  \end{align}
 \eenui
\eprop
\bproof
 \eqref{P5*7.a} In view of the choice of $F$, we get
 \beql{P5*7.1}
  F
  \in\SFuqa{0}
 \eeq
 and $\OSm{F}\in\MggqKksg$. Thus, we have $\OSma{F}{\rhl}=\su{0}$. Combining this with \eqref{P5*7.1}, we obtain from \rlem{L3*22} all assertions of~\eqref{P5*7.a}.
 
 \eqref{P5*7.b} The assertions of~\eqref{P5*7.b} follow from~\eqref{P5*7.a} by application of \rrem{RA*3}.
 
 \eqref{P5*7.c} From \eqref{P5*7.1} and \eqref{F3*3} we get $F\in\SFqa$. Now the assertions of~\eqref{P5*7.c} follow by combination of \rprop{P2*7} with~\eqref{P5*7.a}.
\eproof

 The next result establishes a connection to the class $\SFdqaa{\su{0}}$ introduced in \rsec{S*4}.

\bleml{L5*8}
 Let \(\ug\in\R\), let \(\kappa\in\NOinf \) and let \(\seqska \) be a sequence of complex \tqqa{matrices}. Then \(\SFqas{\kappa}\subseteq\SFdqaa{\su{0}}\).
\elem
\bproof
 Let \(F\in\SFqas{\kappa}\). Then \(F\in\SFuqa{0}\) and the \taSm{} \(\OSm{F}\) of \(F\) belongs to \(\MggqKksg\). From \eqref{Nr.TMB} we get then \(F\in\SFdqa\) and, in view of \eqref{F1*1}, we have \(\su{0}=\suo{0}{\OSm{F}}=\OSma{F}{\rhl}\). Taking into account \rlem{L3*22}, we conclude then \(\nul{\su{0}}=\nul{\OSma{F}{\rhl}}=\nul{F(z)}\) for all \(z\in\Cs\). From \rrem{R-C4*7} we obtain then \(F\in\SFdqaa{\su{0}}\).
\eproof

\begin{rem}\label{R1530}
 Let \(\ug\in\R\) and let \(F\colon\Cs\to\Cqq\) be a matrix-valued function. Then:
 \benui
  \il{R1530.a} If \(F\in\SFuqa{0}\), then \(F\in\SFqas{0}\) with \(s_0\defeq\OSma{F}{\rhl}\).
  \il{R1530.b} Let \(\kappa\in\NOinf\) and let \(\seqska\in\Kggeqka\). If \(F\in\SFqas{\kappa}\), then \(F\in\SFuqa{0}\) and \(\OSma{F}{\rhl}=s_0\).
 \eenui
\end{rem}

\section{On \Hamburger{}-\Nevanlinna{}-type results for $\SFqas{\kappa}$}\label{S*6}
 This section should be compared with~\zitaa{MR3380267}{\csec{6}}, where we investigated matricial generalizations of the classical Hamburger-Nevanlinna theorems for the class \(\RFq\) of \tqqa{matrix-valued} functions which are holomorphic in \(\ohe\) and have \tnnH{} imaginary values at each point. Using these results, we were able to reformulate the truncated matricial Hamburger moment problem in~\zita{MR3380267} into a problem of determining all functions from special subclasses of \(\RFq\) which have prescribed asymptotic expansion in some sector of \(\ohe\).
 
 The main goal of this section can be described as follows. Let $m\in\NO$ and let $\seqs{m}$ be a sequence of complex \tqqa{matrices}. Then we are going to reformulate the moment problem~\mproblem{\rhl}{m}{=} into a problem of determining all functions from special subclasses of \(\SFqa\) which have a prescribed asymptotic expansion in some sector of \(\ohe\). This aim can be realized by application of a Stieltjes-type version of a matricial Hamburger-Nevanlinna theorem. Before formulating the result, we introduce some notation. For each \(\Omega\in\BorR\setminus\set{\emptyset}\), each \(\kappa\in\NOinf\), and each sequence \(\seqska\) of complex \tqqa{matrices} let \(\Mggqaag{\Omega}{\seqska}\)\index{m@\(\Mggqaag{\Omega}{\seqska}\)} be the set of all \(\sigma\in\Mgguqa{\kappa}{\Omega}\) such that \(\suo{j}{\sigma}=\su{j}\) is fulfilled for each \(j\in\mn{0}{m}\). For all $r\in (0,\infp)$ and each $\delta \in(0,\pi/2 ]$, let\index{s@$\ohesecuu{r}{\delta}$}  
\[
 \ohesecuu{r}{\delta}
 \defeq\setaa*{z \in \C}{\abs{z}\geq r\text{ and }\delta \leq \arg z \leq \pi -\delta}.
\]
 Taking \eqref{F3*2} into account, now we can formulate the following Stieltjes-type version of the  Hamburger-Nevanlinna theorem:

\bthml{T6*1}
 Let \(\Omega\) be a \tne{} closed subset of \(\R\), let $m\in\NO$, and let $\seqs{m}$ be a sequence of complex \tqqa{matrices}.
 \benui
  \il{T6*1.a} Let \(\sigma\in\Mggouaag{q}{m}{\Omega}{\seqs{m}}\) and let \(\OSt{\sigma}\) be the \tOSto{\Omega}{\sigma}. For all \(\ell\in\mn{0}{m}\) and all \(\delta\in (0,\pi/2]\), then
  \[
   \lim_{r\to \infp} \sup_{z\in\ohesecuu{r}{\delta}}\normE*{z^{\ell+1} \ek*{ \OSta{\sigma}{z} +\sum_{j=0}^\ell z^{-(j+1)} s_j}}
   = 0. 
  \]
  \il{T6*1.b} Suppose that \(\inf\set{\abs{\inf \Omega}, \abs{\sup \Omega}}<\infp\) and that the matrices \(\su{0},\su{1},\dotsc,\su{m}\) are \tH{}. Let \(\sigma\in\Mggqa{\Omega}\) be such that the \tOSt{\Omega} \(\OSt{\sigma}\) of \(\sigma\) fulfills
  \beql{T6*1.A}
   \lim_{y\to\infp}\normE*{(\iu y)^{m+1}\ek*{\OSta{\sigma}{\iu y}+\sum_{j=0}^m(\iu y)^{-(j+1)}\su{j}}}
   =0.
  \eeq
  Then  \(\sigma\) belongs to \(\Mggouaag{q}{m}{\Omega}{\seqs{m}}\).
 \eenui
\ethm

 \rthm{T6*1} should be compared with~\zitaa{MR3380267}{\cthm{6.1}}. A closer look shows that \rthm{T6*1} simultaneously works for all \(m\in\NO\) whereas in~\zitaa{MR3380267}{\cthm{6.1}} it is assumed that \(m\) is even.

 \rPart{T6*1.b} of \rthm{T6*1} will be often applied in the following for \(\Omega=\rhl\). It contains a sufficient condition which implies that a function $F\in\SFuqa{0}$ is the \taSt{} of a solution $\sigma$ of Problem~\mproblem{\rhl}{m}{=}.
 
 The proof of \rpart{T6*1.a} of \rthm{T6*1} can be done analogously to the proof in the case that \(\Omega=\R\) and \(m=2n\) with some \(n\in\NO\) hold (see, \eg{}~\zita{MR703593},~\zitaa{MR1624548}{\clem{2.1}}). We omit the details. Now we are going to prove \rpart{T6*1.b} of \rthm{T6*1}. In order to do this, we need some auxiliary results:

\begin{rem}\label{63C}
 Let \(m\in\NO\), let \(\mathcal{G}\) be a subset of \(\C\) such that \(\setaa{\iu y}{y\in (0,\infp)} \subseteq \mathcal{G}\), let \(F\colon\mathcal{G}\to \Cpq \) be a matrix-valued function, and let \(\seqs{m}\) be a sequence of complex \tpqa{matrices} such that
 \begin{equation}\label{Nr.L2}
  \lim_{y\to\infp}(\iu y)^{\ell+1} \ek*{F (\iu y) + \sum_{j=0}^\ell (\iu y)^{-(j+1)} s_j }
  =\Opq
 \end{equation}
 holds for \(\ell=m\). Then one can easily check by induction that \eqref{Nr.L2} is valid for each \(\ell\in\mn{0}{m}\).
\end{rem}

\begin{lem}\label{64A} 
 Let \(\Omega\) be a \tne{} closed subset of \(\R\) such  that \(\inf\set{\abs{\inf\Omega},\abs{\sup \Omega}}<\infp\), let \(m\in \N\), let \(\seqs{m-1}\) be a sequence of complex \tqqa{matrices}, and let \(s_m\) be a \tH{} \tqqa{matrix}. Let \(\sigma\in\Mggouaag{q}{m-1}{\Omega}{\seqs{m-1}}\) be such that the \tOSt{\Omega} \(\OSt{\sigma}\) of \(\sigma\) fulfills
 \begin{equation}\label{Nr.SIY}
  \lim_{y\to\infp} (\iu y)^{m+1} \ek*{\OSta{\sigma}{\iu y} +\sum_{j=0}^m(\iu y)^{-(j+1)} s_j}
  =\Oqq.
 \end{equation}
 Then \(\sigma\) belongs to \(\Mggouaag{q}{m}{\Omega}{\seqs{m}}\).
\end{lem}
\begin{proof}
 Because of \eqref{Nr.SIY} we have
\begin{equation}\label{Nr.SIZ}
 s_m
 = -\lim_{n\to\infi} (\iu n)^{m+1} \ek*{ \OSta{\sigma}{\iu n} +\sum_{j=0}^{m-1} (\iu n)^{-(j+1)} s_j}.
\end{equation}
 For each \(n\in\N\), let \(f_n\colon\Omega\to\C\) and \(g_n\colon\Omega\to\C\) be defined by
\begin{align*}
f_n (t)&\defeq  \frac{n^2 t^m}{t^2 + n^2} -\iu \frac{nt^{m+1}}{t^2 + n^2}&
&\text{and}&
g_n (t)&\defeq  \frac{n^2 t^m}{t^2 + n^2}.
\end{align*}
For each \(n\in\N\) we get in view of \(\sigma\in\Mggouaag{q}{m-1}{\Omega}{\seqs{m-1}}\) then
\[
\begin{split}
 & -(\iu n)^{m+1} \ek*{ \OSta{\sigma}{\iu n} +\sum_{j=0}^{m-1} (\iu n)^{-(j+1)} s_j}\\
 & = - (\iu n)^{m+1} \int_{\Omega} \frac{1}{t-\iu n} \sigma (\dif t) - \sum_{j=0}^{m-1} (\iu n)^{m-j}\int_\Omega t^j \sigma (\dif t)\\
 & = \int_\Omega \frac{1}{t-\iu n} \ek*{ -(\iu n)^{m+1} - \sum_{j=0}^{m-1} (\iu n)^{m-j} t^j (t-\iu n)} \sigma (\dif t)\\
 & = \int_\Omega \frac{-\iu n t^m}{t - \iu n} \sigma (\dif t)
 =\int_\Omega \frac{-\iu n t^{m+1} - (\iu n)^2 t^m}{t^2 + n^2}\sigma (\dif t)
 =\int_\Omega f_n\dif\sigma
\end{split}
\]
 and in particular \(f_n\in\LaaaC{\Omega}{\BorO}{\sigma}\). Hence, for each \(n\in\N\), because of \(g_n = \re  f_n\) and \rrem{R1647} we obtain 
\begin{equation}\label{Nr.GG}
 g_n
 \in\LaaaC{\Omega}{\BorO}{\sigma}
\end{equation}
 and
\begin{equation}\label{Nr.L7}
 \re  \rk*{ -(\iu n)^{m+1} \ek*{ \OSta{\sigma}{\iu n} +\sum_{j=0}^{m-1} (\iu n)^{-(j+1)} s_j} }
 =\int_\Omega g_n \dif\sigma.
\end{equation}
 Using \eqref{Nr.L7}, \eqref{Nr.SIZ} and \(\su{m}^\ad=\su{m}\) we get
 \begin{equation}\label{Nr.D77}\begin{split}
  \lim_{n\to\infi} \int_\Omega g_n \dif\sigma
  &=\lim_{n\to\infi}\re\rk*{ -(\iu n)^{m+1} \ek*{ \OSta{\sigma}{\iu n} +\sum_{j=0}^{m-1} (\iu n)^{-(j+1)} s_j} }\\
  &=\re\rk*{ -\lim_{n\to\infi}(\iu n)^{m+1} \ek*{ \OSta{\sigma}{\iu n} +\sum_{j=0}^{m-1} (\iu n)^{-(j+1)} s_j} }
  =\re s_m 
  =\su{m}.
 \end{split}\end{equation}
 Now we are going to prove that the function \(g\colon\Omega\to\C\) given by \(g(t) \defeq  t^m\) belongs to \(\LaaaC{\Omega}{\BorO}{\sigma}\) and that \(\lim_{n\to\infi} \int_\Omega g_n \dif\sigma = \int_\Omega g\dif\sigma\). Obviously, for each \(t\in\Omega\) we have
\begin{equation}\label{Nr.D15}
 \lim_{n\to\infi} g_n (t) 
 = t^m = g (t). 
\end{equation}

 First we consider the case
\begin{equation}\label{Nr.D8}
\inf \set*{\abs{\inf \Omega},\abs{\sup \Omega}} 
=\abs{\inf \Omega}.  
\end{equation}
 Then \(\inf \Omega\in\R\) and \(\ug \defeq  \min \set{0,\inf \Omega}\) is a non-positive real number. Obviously,
\begin{align*}
 [\ug,0] \cup (0,\infp) &= \rhl ,&
 [\ug, 0] \cap (0,\infp)&=\emptyset&
 &\text{and}&
 \Omega&\subseteq \rhl.
\end{align*}
 Thus the sets \(\cA \defeq  \Omega\cap [\ug, 0]\) and \(\cB \defeq  \Omega\cap (0,\infp)\) are disjoint and fulfill \(\cA\cup\cB =\Omega\). Consequently,
 \begin{equation}\label{64A.1}
  1_\cA + 1_\cB
  = 1_\Omega
 \end{equation}
 where, for each \(\mathcal{M}\in\BorO \), the notation \(1_\mathcal{M}\) denotes the indicator function of \(\mathcal{M}\) (defined on \(\Omega\)). For each \(t\in\Omega\), from \eqref{Nr.D15} we conclude that 
\begin{align}\label{Nr.D16}
\lim_{n\to\infi} (1_\cA g_n) (t)&= (1_\cA g) (t)&
&\text{and}&
\lim_{n\to\infi} (1_\cB g_n) (t)&= (1_\cB g) (t)
\end{align}
 hold. For each \(n\in\N\) and each \(t\in\Omega\) the estimate \(\abs{(1_\cA g_n) (t)} \leq 1_\cA  (t) \cdot\abs{t}^m\le\abs{\ug}^m 1_\cA  (t)\) is true. Since \(\abs{\ug}^m 1_\cA\) belongs to \(\LaaaC{\Omega}{\BorO}{\sigma}\) application of \rprop{P1102} yields then \(1_\cA g \in\LaaaC{\Omega}{\BorO}{\sigma}\), \(1_\cA g_n \in\LaaaC{\Omega}{\BorO}{\sigma}\) for all \(n\in\N\), and 
\begin{equation}\label{Nr.D20}
 \lim_{n\to\infi} \int_\Omega 1_\cA g_n \dif\sigma
 =\int_\Omega 1_\cA g \dif\sigma.
\end{equation}
 We consider an arbitrary \(x\in\Cq\). For each \(n\in\N\), from \eqref{Nr.GG} we see that \(1_\cB g_n\) belongs to \(\LaaaR{\Omega}{\BorO}{x^\ad\sigma x}\). Moreover, one can easily check that \(0\leq 1_\cB g_n\leq 1_\cB g_{n+1}\) holds for every \(n\in\N\). Because of \(1_\cA + 1_\cB = 1_\Omega\), we have \(g_n -1_\cA g_n = 1_\cB g_n\) for each \(n\in\N\). Using this, \eqref{Nr.D77}, and \eqref{Nr.D20}, we conclude that
\begin{equation}\label{Nr.D26}
 x^\ad\rk*{ s_m - \int_\Omega 1_\cA g \dif\sigma} x
 =  \lim_{n\to\infi} \int_\Omega 1_\cB g_n \dif(x^\ad\sigma x).
\end{equation}
 Hence, in view of \eqref{Nr.D16}, \eqref{Nr.D26}, and the monotone convergence theorem, we get \(1_\cB g\in \LaaaR{\Omega}{\BorO}{x^\ad\sigma x}\) and \(\lim_{n\to\infi} \int_\Omega 1_\cB g_n  \dif(x^\ad\sigma x)= \int_\Omega 1_\cB g \dif(x^\ad\sigma x)\). Consequently, \(1_\cB g\in \LaaaR{\Omega}{\BorO}{\sigma}\) and
\begin{equation}\label{Nr.D29}
 \lim_{n\to\infi} \int_\Omega 1_\cB g_n  \dif\sigma
 = \int_\Omega 1_\cB g \dif\sigma
\end{equation}
 hold. Since \eqref{64A.1} the identity \(g =  1_\cA g + 1_\cB g\) is valid and since the functions \(1_\cA g\) and \(1_\cB g\) both belong to \(\LaaaR{\Omega}{\BorO}{\sigma}\) we get
\begin{equation}\label{Nr.GL}
 g
 \in  \LaaaR{\Omega}{\BorO}{\sigma}
\end{equation}
 and
\[
 \int_\Omega g \dif\sigma
 = \int_\Omega 1_\cA g  \dif\sigma  + \int_\Omega 1_\cB g  \dif\sigma.
\]
 From \eqref{Nr.D20}, \eqref{Nr.D29}, \eqref{64A.1} and \eqref{Nr.D77} we infer then
 \begin{equation}\label{64A.2}\begin{split}
  \int_\Omega t^m\sigma(\dif t)
  =\int_\Omega g\dif\sigma
  &=\lim_{n\to\infi}\int_\Omega1_\cA g_n\dif\sigma+\lim_{n\to\infi}\int_\Omega1_\cB g_n\dif\sigma\\
  &=\lim_{n\to\infi}\int_\Omega\rk{1_\cA g_n+1_\cB g_n}\dif\sigma
  =\lim_{n\to\infi}\int_\Omega g_n\dif\sigma
  =\su{m}.
 \end{split}\end{equation}

 Now we consider the case \(\inf\set{\abs{\inf \Omega},\abs{\sup \Omega}}=\abs{\sup \Omega}\). Then we have \(\sup \Omega\in\R\) and \(\beta \defeq  \max\set{0, \sup \Omega}\) is a \tnn{} real number. Obviously, \([0,\beta]\cup (-\infty, 0) = (-\infty, \beta]\), \([0,\beta] \cap (-\infty, 0) =\emptyset\), and \(\Omega \subseteq (-\infty, \beta]\). Thus the sets \(\cC \defeq  \Omega\cap [0,\beta]\) and \(\cD \defeq \Omega\cap (-\infty, 0)\) are disjoint and fulfill \(\cC \cup\cD =\Omega\). Analogous to the case \eqref{Nr.D8} one can prove then that \eqref{Nr.GL} and \eqref{64A.2} hold. Thus \eqref{Nr.GL} and \eqref{64A.2} are proved for each case. Because of \(\sigma \in\Mggouaag{q}{m-1}{\Omega}{\seqs{m-1}}\) and \eqref{64A.2}, we conclude then that \(\sigma\) belongs to \(\Mggouaag{q}{m}{\Omega}{\seqs{m}}\).
\end{proof}

\begin{proof}[Proof of \rpart{T6*1.b} of \rthm{T6*1}]
 Using \rlemp{63A}{63A.c} one can easily check that
 \[
 \suo{0}{\sigma}
 = \sigma (\Omega) 
 =  -\iu\lim_{y\to\infp} y  \OSta{\sigma}{\iu y}
 =  -\ek*{ \lim_{y\to\infp} \iu y \OSta{\sigma}{\iu y}+ s_0} + s_0.
 \]
 In view of \eqref{T6*1.A}, the application of \rrem{63C} yields then
 \[
  \suo{0}{\sigma} 
  = -\lim_{y\to\infp} \iu y  \ek*{ \OSta{\sigma}{\iu y}+ (\iu y)^\inv  s_0} + s_0 
  = s_0.
 \]
 Hence \(\sigma\) belongs to \(\Mggouaag{q}{0}{\Omega}{\seqs{0}}\). Consequently, there is an \(\ell\in\mn{0}{m}\) such that \(\sigma\) belongs to \(\Mggouaag{q}{\ell}{\Omega}{\seqs{\ell}}\). We consider the case \(\ell<m\). Then \eqref{T6*1.A} and \rrem{63C} provide us
 \[
   \lim_{y\to\infp} (\iu y)^{\ell +2} \ek*{ \OSta{\sigma}{\iu y}+ \sum_{j=0}^{\ell +1} (\iu y)^{-(j+1)} s_j} 
   =\Oqq.
 \]
 Application of \rlem{64A} yields then \(\sigma\in\Mggouaag{q}{\ell+1}{\Omega}{\seqs{\ell+1}}\). Consequently, \(\sigma\in\Mggouaag{q}{m}{\Omega}{\seqs{m}}\) is inductively proved.
\end{proof}

 \rthm{T6*1} was inspired by Benaych-Georges~\zitaa{MR2225442}{\cthm{1.3}}. This result was obtained in the context of free probability in the scalar case for the so-called \(R\)\nobreakdash-transform. The \(R\)\nobreakdash-transform realizes a bijective correspondence between the set of probability measures on \(\R\) and some set of holomorphic functions. The \(R\)\nobreakdash-transform is different from our \tOSt{\Omega} in \eqref{F3*2}. For this reason, even in the scalar case \rthm{T6*1} is different from~\zitaa{MR2225442}{\cthm{1.3}}. However in a wider sense both results have some common features. In particular, the assumption of the semiboundedness of the support of the measure in \rpart{T6*1.b} of \rthm{T6*1} was suggested by~\zitaa{MR2225442}{\cthm{1.3(b)}}.
 
 Now we give corollaries of the Stieltjes-type matricial \Hamburger{}-\Nevanlinna{} \rthm{T6*1} for the particular case \(\Omega=\rhl\).
\bcorl{C6*3}
 Let \(\ug\in\R\), let \(m\in\NO\), let \(\seqs{m}\) be a sequence of complex \tqqa{matrices}, and let \(F\in\SFqas{m}\). For all \(\ell \in\mn{0}{m}\) and all \(\delta\in(0,\pi/2]\), then
 \[
  \lim_{r\to\infp}\sup_{z\in\ohesecuu{r}{\delta}}\normE*{z^{\ell +1}\ek*{F(z)+\sum_{j=0}^\ell z^{-(j+1)}\su{j}}}
  =0.
 \]
\ecor
\bproof
 Let \(\sigma\) be the \taSmo{F} and let \(\OSt{\sigma}\) be the \taSto{\sigma}. Then \(\OSt{\sigma}=F\) and, in view of \eqref{F5*3}, furthermore \(\sigma\in\MggqKsg{m}\). The application of \rthmp{T6*1}{T6*1.a} completes the proof.
\eproof

\bcorl{C-T6*4}
 Let \(\ug\in\R\), let \(m\in\NO\), let \(\seqs{m}\) be a sequence of \tH{} \tqqa{matrices}, and let \(F\in\SFqa\) be such that
 \beql{C-T6*4.1}
  \lim_{y\to\infp}(\iu y)^{m+1}\ek*{ F(\iu y)+\sum_{j=0}^m(\iu y)^{-(j+1)}\su{j}}
  =\Oqq.
 \eeq
 Then \(F\) belongs to the class \(\SFqas{m}\).
\ecor
\bproof
 \rrem{63C} shows that \eqref{C-T6*4.1} holds for \(m=0\). This implies that \(\sup_{y\in[1,\infp)}y\normE{F(\iu y)}<\infp\) is valid, \ie{}, \(F\) belongs to \(\SFuqa{0}\). In particular, \(\OSt{\sigma}=F\), where \(\sigma\) is the \taSmo{F} and \(\OSt{\sigma}\) is the \taSto{\sigma}. Using \rthmp{T6*1}{T6*1.b} and \eqref{F5*3}, the proof is finished.
\eproof

 Let $m\in\NO$ and let $\seqs{m}$ be a sequence of Hermitian \tqqa{matrices}. Then our preceding considerations lead us to the desired description of the set $\SFqas{m}$.

\bthml{T-P6*5}
 Let \(\ug\in\R\), let \(m\in\NO\), and let \(\seqs{m}\) be a sequence of \tH{} \tqqa{matrices}. Then 
 \begin{multline*}
  \SFqas{m}\\
  =\setaa*{F\in\SFqa}{\lim_{y\to\infp}(\iu y)^{m+1}\ek*{ F(\iu y)+\sum_{j=0}^m(\iu y)^{-(j+1)}\su{j}}=\Oqq}.
 \end{multline*}
\ethm
\bproof
 Combine \rcorss{C6*3}{C-T6*4}.
\eproof

 \rthm{T-P6*5} is the main result of this section. In combination with \rthm{T-C5*6}, we recognize now that \rprob{\mproblem{\rhl}{m}{=}} via \taSt{} can be reformulated into an equivalent problem of prescribing the asymptotic behavior on the positive imaginary axis of functions belonging to \(\SFqa\), with the aid of the given sequence \(\seqs{m}\) from \(\Cqq\).
 
 \rthm{T-P6*5} should be compared with~\zitaa{MR3380267}{\cthm{6.5}}. A closer look shows that \rthm{T-P6*5}  works for all \(m\in\NO\), whereas in~\zitaa{MR3380267}{\cthm{6.5}} it is assumed that \(m\) is even. In~\zita{MR3380267}, the case of an odd \(m\in\NO\) requires a separate treatment (see~\zitaa{MR3380267}{\cthm{6.6}, \cprop{6.7}}).

\section{On a Schur-type algorithm for sequences of complex ${p\times q}$~matrices}\label{S*7}
 In this section, we recall some essential facts on the Schur-type algorithm for sequences from $\Cpq$, which was introduced and investigated in~\zita{114arxiv}. The elementary step of this algorithm is based on the use of the reciprocal sequence of a finite or infinite sequence from $\Cpq$. For this reason, we first remember the definition of the reciprocal sequence:

 Let $\kappa \in \NOinf $ and let $\seqska $ be a sequence of complex \tpqa{matrices}. Then the sequence \symb{\seq{\su{j}^\rez }{j}{0}{\kappa}} given by $\su{0}^\rez \defeq \su{0}^\mpi$ and, for all $k \in \mn{1}{\kappa}$, recursively by
\[
 \su{k}^\rez 
 \defeq - \su{0}^\mpi\sum_{j=0}^{k-1}\su{k-j}\su{j}^\rez ,
\]
 is called the \notii{reciprocal sequence corresponding to  $\seqska $}. For a detailed treatment of the concept of reciprocal sequences, we refer the reader to~\zita{MR3014197}. Furthermore, we recall the notion of \tspat{ation} for sequences of complex matrices (\cf{}~\zitaa{114arxiv}{\cdefn{4.1}}):

 Let \(\ug\in\C\), let \(\kappa\in\NOinf \), and let \(\seqska \) be a sequence of complex \tpqa{matrices}. Let \(\su{-1}\defeq\Oqq\). Then the sequence \symb{\seq{\su{j}^\spa}{j}{0}{\kappa}} given by
 \beql{F-D1455}
  \su{j}^\spa
  \defeq-\ug\su{j-1}+\su{j}
 \eeq
 is called the \notii{\tspato{\seqska }}. In order to prepare the basic construction, we use the \tsraute{} corresponding to the \tspat{} of a sequence of complex matrices:
 
 Let \(\ug\in\C\), let \(\kappa\in\NOinf \), and let \(\seqska \) be a sequence of complex \tpqa{matrices} with \tspat{} \(\seq{u_j}{j}{0}{\kappa}\). Then the sequence \symb{\seq{\su{j}^\reza }{j}{0}{\kappa}} is given by
\[
 \su{j}^\reza 
 \defeq u_j^\rez,
\]
 \ie{}, the sequence \(\seq{\su{j}^\reza }{j}{0}{\kappa}\) is defined to be the \tsrautea{the \tspato{\seqska }}. Now we explain the elementary step of the Schur-type algorithm under consideration (\cf{}~\zitaa{114arxiv}{\cdefn{7.1}}):
 
 Let \(\ug\in\C\), let \(\kappa\in\Ninf \), and let \(\seqska \) be a sequence of complex \tpqa{matrices}. Then the sequence \(\seq{t_j}{j}{0}{\kappa-1}\) given by
 \beql{F*sta}
  t_j
  \defeq-\su{0}\su{j+1}^\reza \su{0}
 \eeq
 is called the \notii{\taScht{}} (or short \notii{\saScht{}}) \emph{of \(\seqska $}.

\bremnl{\cf{}~\zitaa{114arxiv}{\crem{7.3}}}{R7*1}
 Let \(\ug\in\C\), let \(\kappa\in\Ninf\), and let \(\seqska\) be a sequence of complex \tpqa{matrices} with \saScht{} \(\seq{t_j}{j}{0}{\kappa-1}\). For all \(m\in\mn{1}{\kappa}\), then the sequence \(\seq{t_j}{j}{0}{m-1}\) depends only on the matrices \(\su{0},\su{1},\dotsc,\su{m}\) and, hence, it coincides with the \saSchto{\seqs{m}}.
\erem

 The repeated application of the \taScht{} generates in a natural way a corresponding algorithm for (finite or infinite) sequences of complex \tpqa{matrices} (\cf{}~\zitaa{114arxiv}{\cdefn{8.1}}):

 Let \(\ug\in\C\), let \(\kappa\in\NOinf \), and let \(\seqska \) be a sequence of complex \tpqa{matrices}. The sequence \(\seq{\su{j}^\sta{0}}{j}{0}{\kappa}\) given by \(\su{j}^\sta{0}\defeq\su{j}\) is called the \notii{\sthaSchto{0}{\seqska}}. In the case \(\kappa\geq1\), for all \(k\in\mn{1}{\kappa}\), the \notii{\sthaScht{k} \(\seq{\su{j}^\sta{k}}{j}{0}{\kappa-k}\) of \(\seqska\)}\index{$\seq{\su{j}^\sta{k}}{j}{0}{\kappa-k}$} is recursively defined to be the \taScht{} of the \sthaScht{(k-1)} \(\seq{\su{j}^\sta{k-1}}{j}{0}{\kappa-(k-1)}\) of \(\seqska\).

 One of the central properties of the just introduced Schur-type algorithm is that it preserves the \taSnnd{} extendability of sequences of matrices. This is the content of the following result.

\bthmnl{\cf{}~\zitaa{114arxiv}{\cthm{8.10}}}{T-P7*2}
 Let \(\ug\in\R\), let \(\kappa\in\NOinf\), let \(\seqska\in\Kggqka \), and let \(k\in\mn{0}{\kappa}\). Then the \sthaScht{k} \(\seq{\su{j}^\sta{k}}{j}{0}{\kappa-k}\) of \(\seqska\) belongs to \(\Kggequa{\kappa-k}\).
\ethm

 In our considerations below, the special parametrization introduced in~\zita{MR3014201} will play an essential role, the so-called \taSp{}. For the convenience of the reader, we recall this notion. Let \(\ug\in\R\), let $\kappa\in\NOinf$, and let $\seqska$ be a sequence of complex \tpqa{matrices}. For every choice of integers $\ell,m\in\NO$ with $\ell\leq m\leq\kappa$, let\index{y@$\yuu{\ell}{m}$}\index{z@$\zuu{\ell}{m}$}
\begin{align*}
 \yuu{\ell}{m}&\defeq
 \bMat
  \su{\ell}\\
  \su{\ell+1}\\
  \vdots\\
  \su{m}
 \eMat&
 &\text{and}&
 \zuu{\ell}{m}&\defeq\mat{\su{\ell},\su{\ell+1},\dotsc,\su{m}}.
\end{align*}
 For all $n\in\NO$ with $2n\leq\kappa$, let\index{h@$\Hu{n}$}
\[
 \Hu{n}
 \defeq\matauuo{\su{j+k}}{j,k}{0}{n}
\]
 and, for all $n\in\NO$ with $2n+1\leq\kappa$, let\index{k@$\Ku{n}$}
\[
 \Ku{n}
 \defeq\matauuo{\su{j+k+1}}{j,k}{0}{n}.
\]

\begin{defn}[\cf{}~\zitaa{MR3014201}{\cdefn{4.2}}]\label{D7*3}
 Let \(\ug\in\R\), let \(\kappa\in\NOinf \), and let \(\seqska \) be a sequence of complex \tpqa{matrices}. Then the sequence \(\seq{\Spu{j}}{j}{0}{\kappa}\)\index{q@$\seq{\Spu{j}}{j}{0}{\kappa}$} given by \(\Spu{0}\defeq\su{0}\) and
 \[
  \Spu{2k}
  \defeq\su{2k}-\zuu{k}{2k-1}\Hu{k-1}^\mpi\yuu{k}{2k-1}
 \]
 for all \(k\in\N\) with \(2k\leq\kappa\) and by \(\Spu{1}\defeq-\ug\su{0}+\su{1}\) and
 \[
  \Spu{2k+1}
  \defeq(-\ug\su{2k}+\su{2k+1})-(-\ug\zuu{k}{2k-1}+\zuu{k+1}{2k})(-\ug\Hu{k-1}+\Ku{k-1})^\mpi(-\ug\yuu{k}{2k-1}+\yuu{k+1}{2k})
 \]
 for all \(k\in\N\) with \(2k+1\leq\kappa\) is called the \notii{\taSpo{\seqska}}.
\end{defn}

 In~\zitas{MR3014201}, several important classes of sequences of complex \tqqa{matrices} were characterized in terms of their \taSp{}. From the view of this paper, the class $\Kggeqka$ of \taSnnde{} sequences is of extreme importance (see \rthmss{T1*3+2}{T-P5*4}). In the case of a sequence $\seqska\in\Kggeqka$, the \taSp{} can be generated by the above constructed Schur-type algorithm. This is the content of the following theorem.

\bthmnl{\zitaa{114arxiv}{\cthm{9.15}}}{T7*5}
 Let \(\ug\in\R\), let \(\kappa\in\NOinf\), and let \(\seqska\in\Kggeqka\). Then \(\seq{\su{0}^\sta{j}}{j}{0}{\kappa}\) is exactly the \taSpo{\seqska}.
\ethm

 An essential step in the further considerations of this paper can be described as follows. Let \(\ug\in\R\), let \(\kappa\in\Ninf\) and let $\seqska \in\Kggeqka$ with \taScht{} $\seq{t_j}{j}{0}{\kappa-1}$. In view of \rthm{T-P7*2}, we get then $\seq{t_j}{j}{0}{\kappa-1}\in\Kggequa{\kappa-1}$. Thus, \rthm{T-P5*4} yields that both sets $\SFqas{\kappa}$ and $\SFuqaa{\kappa-1}{\seq{t_j}{j}{0}{\kappa-1}}$ are \tne{}. Then a central aspect of our strategy is based on the construction of a special bijective mapping $\HTu{\su{0}}$ with the property
\[
 \HTuA{\su{0}}{\SFqaS{\kappa}}
 =\SFuqAA{\kappa-1}{\seq{t_j}{j}{0}{\kappa-1}}.
\]
 This mapping $\HTu{\su{0}}$ will help us to realize the basic step for our construction of a special Schur-type algorithm in the class $\SFqa$, which stands in a bijective correspondence to the above described Schur-type algorithm for \taSnnde{} sequences.
 
 Against to this background we will use the inverse of the \saScht{}.
 
\begin{defn}[\zitaa{114arxiv}{\cdefn{10.1}}]\label{114.D1712}
 Let \(\ug\in\C\), let \(\kappa\in\NOinf \), let \(\seqska\) be a sequence of complex \tpqa{matrices}, and let \(A\) be a complex \tpqa{matrix}. The sequence \(\seq{r_j}{j}{0}{\kappa+1}\) with \(r_0\defeq A\) recursively defined by
 \[
  r_j
  \defeq\ug^jA+\sum_{\ell=1}^j\ug^{j-\ell}AA^\mpi\ek*{\sum_{k=0}^{\ell-1}\su{\ell-k-1}A^\mpi r_k^\spa}
 \]
 is called the \notii{\siaSchto{\seqska}{A}}.
\end{defn}

\section{Further considerations on the class of \taSnnde{} sequences and some subclass}\label{S1459}
 In order to verify several interrelations between the algebraic and function theoretic versions of our Schur-type algorithm, we will need certain properties of \taSnnde{} sequences of complex \tqqa{matrices}. Now we give a short summary of this material:

\bleml{L7*6}
 Let \(\ug\in\R\), let $\kappa\in\NOinf$ and let $\seqska\in\Kggeqka$. Then \(\su{j}\in\CHq\) for all \(j\in\mn{0}{\kappa}\) and $\su{2k}\in\Cggq$ for all $k\in\NO$ with \(2k\leq\kappa\).
\elem
\bproof
 This is an immediate consequence of~\zitaa{MR3014201}{\clem{2.9}}.
\eproof

\bpropnl{\cf{}~\zitaa{MR3014201}{\cthm{4.12(b)}}}{P-L7*7}
 Let \(\ug\in\R\), let \(\kappa\in\NOinf\), and let \(\seqska\) be a sequence of complex \tqqa{matrices} with \taSp{} \(\seqQ{\kappa}\). Then \(\seqska\in\Kggqka\) if and only if \(\Spu{j}\in\Cggq\) for all \(j\in\mn{0}{\kappa}\) and \(\ran{\Spu{j+1}}\subseteq\ran{\Spu{j}}\) for all \(j\in\mn{0}{\kappa-2}\).
\eprop

Now we recall a class of sequences of complex matrices which, as the considerations in~\zita{114arxiv} have shown, turned out to be extremely important in the framework of the above introduced Schur algorithm:

\begin{defn}[\zitaa{MR3014197}{\cdefn{4.3}}]\label{D7*8}
 Let $\kappa\in\NOinf$ and $\seqska $ be a sequence of complex \tpqa{matrices}. We then say that $\seqska $ is \notii{dominated by its first term} (or, simply, that it is \notii{\tftd{}}) if \(\nul{\su{0}}\subseteq\bigcap_{j=0}^{\kappa}\nul{\su{j}}\) and \(\bigcup_{j=0}^{\kappa}\ran{\su{j}}\subseteq\ran{\su{0}}\). The set of all \tftd{} sequences $\seqska $ of complex \tpqa{matrices} will be denoted by \symba{\Dpqu{\kappa}}{d}.
\end{defn}

 For a comprehensive investigation of \tftd{} sequences, we refer the reader to the paper~\zita{MR3014197}. From the view of the Schur-type algorithm for sequences of matrices, the following result proved to be of central importance:
 
\begin{prop}[\zitaa{114arxiv}{\cprop{3.8(a)}}]\label{P-L7*9}
 Let \(\ug\in\R\) and let \(\kappa\in\NOinf \). Then \(\Kggeqka\subseteq\Dqqu{\kappa}\).
\end{prop}

 Now we turn our attention to further important subclasses of the class of all \taSnnd{} sequences. Let \(\ug\in\R\) and let $m\in\NO$. A sequence $\seqs{m}$ of complex \tqqa{matrices} is called \notii{\taSpd{}} if in the case \(m=0\) the block Hankel matrix \(\Hu{0}\) is \tpH{}, in the case \(m=2n\) with some \(n\in\N\) the block Hankel matrices $\Hu{n}$ and $-\ug\Hu{n-1}+\Ku{n-1}$ are \tpH{}, and in the case \(m=2n+1\) with some \(n\in\NO\) the block Hankel matrices $\Hu{n}$ and $-\ug\Hu{n}+\Ku{n}$ are \tpH{}, respectively. A sequence $\seqsinf$ of complex \tqqa{matrices} is called \notii{\taSpd{}} if for all $m\in\NO$ the sequence $\seqs{m}$ is \taSpd{}. For all $\kappa\in\NOinf$, we will write $\Kgqka$\index{k@$\Kgqka$} for the set of all \taSpd{} sequences $\seqska$ of complex \tqqa{matrices}.

\bpropnl{\zitaa{MR3014201}{\cprop{2.20}}}{P-R7*10}
 Let \(\ug\in\R\) and let \(\kappa\in\NOinf\). Then $\Kgqka\subseteq\Kggeqka$.
\eprop

\bpropnl{\zitaa{MR3014201}{\cthm{4.12(d)}}}{P-L7*11}
 Let \(\ug\in\R\), let \(\kappa\in\NOinf\), and let \(\seqska\) be a sequence of complex \tqqa{matrices} with \taSp{} \(\seqQ{\kappa}\). Then \(\seqska\in\Kgqka\) if and only if \(\Spu{j}\in\Cgq\) for all \(j\in\mn{0}{\kappa}\).
\eprop

\section{On a coupled pair of Schur-Stieltjes-type transforms}\label{S*8}
 The main goal of this section is to prepare the elementary step of our Schur-type algorithm for the class $\SFqa$. We will be led to a situation which, roughly speaking, looks as follows: Let \(\ug\in\R\), let \(A\in\Cpq\) and let $F\colon\Cs\to\Cpq$. Then the matrix-valued functions $\STao{F}{A}\colon\Cs\to\Cpq$\index{$\STao{F}{A}$} and $\STiao{F}{A}\colon\Cs\to\Cpq$\index{$\STiao{F}{A}$} which are defined by
\beql{F8*1}
 \STaoa{F}{A}{z}
 \defeq-A\rk*{\Iq+(z-\ug)^\inv\ek*{F(z)}^\mpi A}
\eeq
and
\beql{F8*2}
 \STiaoa{F}{A}{z}
 \defeq-(z-\ug)^\inv   A\ek*{ \Iq  + A^\mpi  F(z)}^\mpi,
\eeq
 respectively, will be central objects in our further considerations.
 
 Against to the background of our later considerations, the matrix-valued functions $\STao{F}{A}$ and $\STiao{F}{A}$ are called the \notii{\taSSt{A}} of \(F\) and the \notii{\tiaSSt{A}} of \(F\).

 The generic case studied here concerns the situation where $p=q$, $A$ is a complex \tqqa{matrix} with later specified properties and $F\in\SFqa$.

 The use of the transforms introduced in \eqref{F8*2} was inspired by some considerations in the papers Chen/Hu~\zitas{MR1807884,MR2038751}. In particular, we mention~\zitaa{MR2038751}{formula~(2.3)}. Before treating more general aspects we state some relevant concrete examples for the constructions given by formulas \eqref{F8*1} and \eqref{F8*2}. First we illustrate the transformations given in \eqref{F8*1} and \eqref{F8*2} by some examples.

\bexaml{E8*1}
 Let $A\in\Cqq$. In view of \eqref{F8*1}, then:
 \benui
  \item Let $\gamma\in\Cqq$ and let $F\colon\Cs\to\Cqq$ be defined by $F(z)\defeq\gamma$. Then $\STaoa{F}{A}{z}=-A+\frac{1}{\ug-z}A\gamma^\mpi A$ for all $z\in\Cs$.
  \item Let $M\in\Cqq$, let $\tau\in\rhl$, and let $F\colon\Cs\to\Cqq$ be defined by $F(z)\defeq\frac{1+\tau-\ug}{\tau-z}M$. Then
  \[
   \STaoa{F}{A}{z}
   =\rk*{\frac{1}{1+\tau-\ug}AM ^\mpi A-A}+\frac{1}{\ug-z}\rk*{\frac{\tau-\ug}{1+\tau-\ug}AM ^\mpi A}
  \]
  for all $z\in\Cs$.
 \eenui
\eexam

 For all $\tau\in\rhl$, let $\Kronu{\tau}$ be the Dirac measure defined on $\BorK$ with unit mass at $\tau$. Furthermore, let the measure $\zmq\colon\BorK\to\Cggq$ be defined by $\zmqa{B}\defeq\Oqq$.

\bexaml{E8*2}
 Let $A\in\Cqq$. In view of \rthm{T2*1} and \rexam{E8*1}, then:
 \benui
  \item Let $\gamma\in\Cggq$ and let $F\colon\Cs\to\Cqq$ be defined by $F(z)\defeq\gamma$. Then $F\in\SFqa$ and $(\gammaF,\muF)=(\gamma,\zmq)$. If $-A\in\Cggq$, then $\STao{F}{A}\in\SFqa$ and $(\gammau{\STao{F}{A}},\muu{\STao{F}{A}})=(-A,\Kronu{\ug}A\gamma^\mpi A)$.
  \item Let $M\in\Cggq$, let $\tau\in\rhl$, and let $F\colon\Cs\to\Cqq$ be defined by $F(z)\defeq\frac{1+\tau-\ug}{\tau-z}M$. Then $F\in\SFqa$ and $(\gammaF,\muF)=(\Oqq,\Kronu{\tau}M)$. If $\frac{1}{1+\tau-\ug}AM^\mpi A\geq A$, then $\STao{F}{A}\in\SFqa$ and 
  \[
   (\gammau{\STao{F}{A}},\muu{\STao{F}{A}})
   =\rk*{\frac{1}{1+\tau-\ug}AM ^\mpi A-A,\frac{\tau-\ug}{1+\tau-\ug}AM ^\mpi A}.
  \]
 \eenui
\eexam

\bexaml{E8*4}
 Let $A,\gamma\in\Cqq$ and let $G\colon\Cs\to\Cqq$ be defined by $G(z)\defeq \gamma$. In view of \eqref{F8*2}, then $\STiaoa{G}{A}{z}=\frac{1}{\ug-z}A(\Iq+A^\mpi\gamma)^\mpi$ for all $z\in\Cs$.
\eexam

\bexaml{E8*5}
 Let $A\in\Cqq$, let $\gamma\in\Cggq$ with $A(\Iq+A^\mpi\gamma)^\mpi\in\Cggq$ and let $G\colon\Cs\to\Cqq$ be defined by $G(z)\defeq\gamma$. In view of \rthm{T2*1} and \rexam{E8*4}, then $G\in\SFqa$ and $(\gammaG,\muG)=(\gamma,\zmq)$, and, furthermore, $\STiao{G}{A}\in\SFqa$ and
 \[
  (\gammau{\STiao{G}{A}},\muu{\STiao{G}{A}})
  =\rk*{\Oqq,\Oqq,\Kronu{\ug}A(\Iq+A^\mpi\gamma)^\mpi}.
 \]
\eexam

A central theme of this paper is to choose, for a given function $F\in\SFqa$, special matrices $A\in\Cqq$ such that the function $\STao{F}{A}$ and $\STiao{F}{A}$, respectively, belongs to $\SFqa$. The following result provides a first contribution to this topic.

\bpropl{P8*6}
 Let $F\in\SFqa$ and let $A\in\Cqq$ be such that $-A\in\Cggq$. Then $\STao{F}{A}\in\SFqa$.
\eprop
\bproof
 Let $G\colon\Cs\to\Cqq$ be defined by
 \beql{P8*6.1}
  G(z)
  \defeq-A.
 \eeq
 In view of $-A\in\Cggq$ and \eqref{P8*6.1}, we see from \rthmp{T2*1}{T2*1.b} then
 \beql{P8*6.2}
  G
  \in\SFqa.
 \eeq
 Let $H\colon\Cs\to\Cqq$ be defined by
 \beql{P8*6.1h}
  H(z)
  \defeq-(z-\ug)^\inv\ek*{F(z)}^\mpi.
 \eeq
 We see from \rprop{P2*9} then
 \beql{P8*6.2h}
  H
  \in\SFqa.
 \eeq
 Taking \eqref{P8*6.2}, \eqref{P8*6.2h} and $A^\ad=A$ into account, we infer from~\zitaa{142}{\crem{3.11}} then
 \beql{P8*6.3}
  G+AHA
  \in\SFqa.
 \eeq
 Because of \eqref{F8*1}, \eqref{P8*6.1} and \eqref{P8*6.1h}, we have
 \beql{P8*6.4}
  \STao{F}{A}
  =G+AHA.
 \eeq
 Using \eqref{P8*6.4} and \eqref{P8*6.3}, we get $\STao{F}{A}\in\SFqa$.
\eproof

 Furthermore, we will show that under appropriate conditions the equations
\begin{align}\label{F8*8}
 \STiao{(\STao{F}{A})}{A}
 &=F&
 &\text{and}&
 \STao{(\STiao{G}{A})}{A}
 &=G
\end{align}
 hold true. The formulas in \eqref{F8*8} show that the functions $\STao{F}{A}$ and $\STiao{G}{A}$ form indeed a coupled pair of transformations. Furthermore it will be clear now the choice of our terminologies ``\taSSt{A}'' and ``\tiaSSt{A}''. If all Moore-Penrose inverses in \eqref{F8*1} and \eqref{F8*2} would be indeed inverse matrices, then the equations in \eqref{F8*8} could be confirmed by straightforward direct computations. Unfortunately, this is not the case in more general situations which are of interest for us. So we have to look for a convenient way to prove the equations in \eqref{F8*8} for situations which will be relevant for us.
 
 Now we verify that in important cases the formulas \eqref{F8*1} and \eqref{F8*2} can be rewritten as linear fractional transformations with appropriately chosen generating matrix-valued functions. The role of these generating functions will be played by the matix polynomials $\mHTu{A}$ and $\mHTiu{A}$ which are studied in \rapp{S*C}. In the sequel we use the terminology for linear fractional transformations of matrices which was introduced in \rapp{S*B}.

\bleml{L8*7}
 Let \(\ug\in\R\), let \(F\colon\Cs\to\Cpq\) be a matrix-valued function, and let \(A\in\Cpq\) be such that \(\ran{F(z)}\subseteq\ran{A}\) and \(\nul{F(z)}\subseteq\nul{A}\) for all \(z\in\Cs\). Then \(F(z)\in\dblftruu{-(z-\ug)A^\mpi}{\Iq-A^\mpi A}\) and \(\STaoa{F}{A}{z}=\lftrooua{p}{q}{\mHTua{A}{z}}{F(z)}\) for all \(z\in\Cs\).
\elem
\bproof
 Let \(z\in\Cs\). The matrix \(Y\defeq(z-\ug)F(z)\) fulfills \(\ran{Y}\subseteq\ran{A}\) and \(\nul{Y}\subseteq\nul{A}\), which in view of \rlem{LA*5} implies \(\ran{Y}=\ran{A}\) und \(\nul{Y}=\nul{A}\). Hence, \(YY^\mpi=AA^\mpi\) and \(Y^\mpi Y=A^\mpi A\). Thus, we obtain
 \[\begin{split}
  (-Y^\mpi A+\Iq-A^\mpi A)&(-A^\mpi Y+\Iq-A^\mpi A)\\
  &=Y^\mpi AA^\mpi Y-Y^\mpi A(\Iq-A^\mpi A)-(\Iq-A^\mpi A)A^\mpi Y+(\Iq-A^\mpi A)^2\\
  &=Y^\mpi YY^\mpi Y+\Iq-A^\mpi A
  =Y^\mpi Y+\Iq-Y^\mpi Y
  =\Iq
 \end{split}\]
 and
 \[\begin{split}
   (Y+A)(-Y^\mpi A+\Iq-A^\mpi A)
   &=-YY^\mpi A+Y-YA^\mpi A-AY^\mpi A+A(\Iq-A^\mpi A)\\
   &=-AA^\mpi A+Y-YY^\mpi Y-AY^\mpi A
   =-A(\Iq+Y^\mpi A).
 \end{split}\]
 In particular, \(\det\ek{-(z-\ug)A^\mpi F(z)+\Iq-A^\mpi A}\neq0\) and \(\ek{-(z-\ug)A^\mpi F(z)+\Iq-A^\mpi A}^\inv=-Y^\mpi A+\Iq-A^\mpi A\). In view of \eqref{F8*1}, \eqref{FC*1}, and \eqref{F*lft}, the proof is complete.
\eproof

 The follwing application of \rlem{L8*7} is important for our considerations in \rsec{S*11}.

\bpropl{L8*9}
 Let \(\ug\in\R\), let $\kappa\in\NOinf$, let $\seqska\in\Kggequa{\kappa}$, and let $F\in\SFqas{\kappa}$. Then $F(z)\in\dblftruu{-(z-\ug)\su{0}^\mpi}{\Iq-\su{0}^\mpi\su{0}}$ and $\lftrooua{q}{q}{\mHTua{\su{0}}{z}}{F(z)}=\STao{F}{\su{0}}(z)$ for all \(z\in\Cs\).
\eprop
\bproof
 In view of \rpropp{P5*7}{P5*7.a} we have
 \begin{align}\label{L8*9.1}
  \Ran{F(z)}&=\ran{\su{0}}&
  &\text{and}&
  \Nul{F(z)}&=\nul{\su{0}}.
 \end{align}
 Taking \eqref{L8*9.1} into account, we infer from \rlem{L8*7} the assertions.
\eproof

 Now we are going to consider the following situation which will turn out to be typical for larger parts of our future considerations. Let $A\in\Cggq$ and $G\in\SFqa$ be such that
\beql{F8*13}
 \ran{\gammau{G}}+\Ran{\muuA{G}{\rhl}}
 \subseteq\ran{A}.
\eeq
 Then our aim is to investigate the function $\STiao{G}{A}$ given by \eqref{F8*2}. We begin by rewriting formula~\eqref{F8*2} as linear fractional transformation. In the sequel we will often use the fact that for \(G\in\SFqa\) the matrix \(\gammaG\) given via \rthmp{T2*1}{T2*1.a} is \tnnH{}.

\bleml{L8*10}
 Let \(\ug\in\R\), let \(A\in\Cggq\), and let \(G\in\SFqa\) be such that \eqref{F8*13} holds. For all \(z\in\Cs\), then \(G(z)\in\dblftruu{(z-\ug)A^\mpi}{(z-\ug)\Iq}\) and \(\STiaoa{G}{A}{z}=\lftrooua{q}{q}{\mHTiua{A}{z}}{G(z)}\).
\elem
\bproof
 Let \(z\in\Cs\). We chose an arbitrary
 \begin{equation}\label{L8*10.0}
  v
  \in\Nul{\Iq+A^\mpi G(z)}. 
 \end{equation}
 According to \rprop{P2*7} and \eqref{F8*13}, we obtain
 \beql{L8*10.1}
  \Ran{G(z)}
  \subseteq\ran{A}.
 \eeq
 Consequently, \rremp{RA*3}{RA*3.b} implies
 \beql{L8*10.2}
  AA^\mpi G(z)
  =G(z).
 \eeq
 Setting \(F\defeq G+A\), the application of \eqref{L8*10.2} and \eqref{L8*10.0} implies
 \beql{Nr.NV}
  \ek*{ F(z)} v
  =A\ek*{\Iq+A^\mpi G(z)} v
  =\Ouu{q}{1}.
 \eeq
 Since the matrix \(\gammau{G}\) is \tnnH{}, we have \(\gammau{G}+A\in\Cggq\). From \rrem{R0832} we get then that \(F\) belongs to \(\SFqa\) and that \(\gammau{F}=\gammau{G}+A\) and \(\muu{F}=\muu{G}\) hold. Using \rprop{P2*7}, we conclude then that
 \beql{Nr.NW}
  \Nul{F(z)}
  =\nul{\gammau{F}}\cap \Nul{\muuA{F}{\rhl }}.
 \eeq
 Furthermore, because of \(\gammau{F}\geq\gammau{G}\geq\Oqq\), we have \(\nul{\gammau{F}}\subseteq\nul{\gammau{G}}\). From this, \eqref{Nr.NW}, and \rprop{P2*7} we get then
 \begin{equation}\label{Nr.NN}
  \begin{split}
  \Nul{F(z)}
  &=\nul{\gammau{F}}\cap \Nul{\muuA{F}{\rhl }}\\
  &\subseteq\nul{\gammau{G}}\cap \Nul{\muuA{G}{\rhl }}
  =\Nul{G(z)}.
 \end{split}
 \end{equation}
 From \eqref{Nr.NV} and \eqref{Nr.NN} we infer \(v\in\nul{G(z)}\). Thus, we have
 \[
  v
  =v+A^\mpi\cdot\Ouu{q}{1}
  =\ek*{\Iq+A^\mpi G(z)} v
  =\Ouu{q}{1}.
 \]
 Combining this with \eqref{L8*10.0} we get \(\nul{\Iq+A^\mpi G(z)}=\set{\Ouu{q}{1}}\). In view of \eqref{F8*2}, \eqref{FC*2}, and \eqref{F*lft} this completes the proof.
\eproof

 The following result plays an essential role in our considerations in \rsec{S*11}.

\bpropl{L8*12}
 Let \(\ug\in\R\), let \(A\in\Cggq\) and let \(G\in\SFdqaa{A}\). Then \(G(z)\in\dblftruu{(z-\ug)A^\mpi}{(z-\ug)\Iq}\) and \(\STiaoa{G}{A}{z}=\lftrooua{q}{q}{\mHTiua{A}{z}}{G(z)}\) for all \(z\in\Cs\).
\eprop
\bproof
 First observe that \(G\in\SFdqa\). In view of \rrem{R3*5}, then \(G\in\SFqa\) and \(\gammau{G}=\Oqq\). From \rrem{R-C4*7} we obtain furthermore \(\ran{\muua{G}{\rhl}}\subseteq\ran{A}\). Thus, we get \(\ran{\gammau{G}}+\ran{\muua{G}{\rhl}}\subseteq\ran{A}\). The application of \rlem{L8*10} yields for all \(z\in\Cs\) then
 \[
  \det\ek*{(z-\ug)A^\mpi G(z)+(z-\ug)\Iq}
  \neq0
 \]
 and furthermore \(\STiaoa{G}{A}{z}=\lftrooua{q}{q}{\mHTiua{A}{z}}{G(z)}\).
\eproof

 Now we formulate the first main result of this section. Assuming the situation of \rlem{L8*10}, we will obtain useful insights into the structure of the \tiaSSto{A}{F}.

\bpropl{P8*14}
 Let \(A\in\Cggq\), let \(\ug\in\R\), let \(G\in\SFqa\) be such that \eqref{F8*13} holds, and let \(u_0\defeq A(\gammaG+A)^\mpi A\). Then \(\STiao{G}{A}\colon\Cs\to\Cqq\) given by \eqref{F8*2} belongs to \(\SFuqaa{0}{\seq{u_j}{j}{0}{0}}\) and fulfills \(\ran{\STiaoa{G}{A}{z}}=\ran{A}\) and \(\nul{\STiaoa{G}{A}{z}}=\nul{A}\) for all \(z\in\Cs\).
\eprop
\bproof
 Since the matrices \(\gammau{G}\) and \(A\) are \tnnH{}, using \rremp{RA*1}{RA*1.a} we conclude that \(u_0\) is \tH{}. Let \(z\in\Cs\). According to \rlem{L8*10}, we get \(\det[\Iq+A^\mpi G(z)]\neq0\) and taking \eqref{FC*2} into account furthermore
 \beql{P8*14.1}
  \STiaoa{G}{A}{z}
  =-(z-\ug)^\inv A\ek*{ \Iq+A^\mpi G(z)}^\inv.
 \eeq
 By virtue of \rprop{P2*7} and \eqref{F8*13} we obtain \eqref{L8*10.1}. Consequently, \rremp{RA*3}{RA*3.b} implies \eqref{L8*10.2}. From \(\gammau{G}\in\Cggq\) and \(A\in\Cggq\) we get \(\gammau{G}+A\in\Cggq\). Thus, \rrem{R0832} yields that \(F\defeq G+A\) belongs to \(\SFqa\) and that \(\gammau{F}=\gammau{G}+A\) and \begin{equation}\label{P8*14.2}
  \muu{F}
  =\muu{G}
 \end{equation}
 hold. Using \rprop{P2*7}, we conclude then \eqref{Nr.NW}. Moreover, because of \(\gammau{F}=\gammau{G}+A\geq A\geq\Oqq\), we get
 \beql{P8*14.3}
  \ran{A}
  \subseteq\ran{\gammau{F}}
 \eeq
 and \(\nul{\gammau{F}}\subseteq\nul{A}\). Because of \eqref{Nr.NW}, we obtain then
 \[
  \Nul{F(z)}
  =\nul{\gammau{F}}\cap \Nul{\muuA{F}{\rhl }}
  \subseteq\nul{\gammau{F}}
  \subseteq\nul{A}.
 \]
 Consequently, \rremp{RA*3}{RA*3.a} implies \(AF^\mpi(z)F(z)=A\). Taking additionally into account \eqref{L8*10.2}, we get
 \[
  AF^\mpi(z)A\ek*{\Iq+A^\mpi G(z)}
  =AF^\mpi(z)\ek*{ A+G(z)}
  =AF^\mpi(z)F(z)
  =A.
 \]
 Hence, in view of \eqref{P8*14.1}, we conclude
 \begin{equation}\label{P8*14.4}
 \begin{split}
  \STiaoa{G}{A}{z}
  &=-(z-\ug)^\inv AF^\mpi(z)A\ek*{\Iq+A^\mpi G(z)}\ek*{ \Iq+A^\mpi G(z)}^\inv\\
  &=-(z-\ug)^\inv AF^\mpi(z)A.
 \end{split}
 \end{equation}
 Let \(H\colon\Cs\to\Cqq\) be defined by \(H(w)\defeq-(w-\ug)^\inv F^\mpi(w)\). Because of \eqref{P8*14.4}, we have then
 \[
  \STiaoa{G}{A}{z}
  =A\ek*{-(z-\ug)^\inv F^\mpi(z)} A
  =A^\ad H(z)A.
 \]
 Since \(F\) belongs to \(\SFqa\), \rprop{P2*9} yields \(H\in\SFqa\). From~\zitaa{142}{\crem{3.11}} we obtain then \(A^\ad HA\in\SFqa\), \ie{}, \(\STiao{G}{A}\in\SFqa\). Thanks to \(F\in\SFqa\), \rprop{P2*5} provides us \(\lim_{y\to\infp}F(\iu y)=\gammaF\), and \rprop{P2*7} yields
 \begin{align}\label{Nr.RG}
  \Ran{F(\iu y)}&=\ran{\gammau{F}}+\Ran{\muuA{F}{\rhl }}&\text{for each }y&\in(0,\infp).
 \end{align}
 In view of \eqref{P8*14.2}, \eqref{F8*13} and \eqref{P8*14.3}, we get
 \[
  \Ran{\muuA{F}{\rhl }}
  =\Ran{\muuA{G}{\rhl }}
  \subseteq\ran{\gammau{G}}+\Ran{\muuA{G}{\rhl }}
  \subseteq\ran{A}
  \subseteq\ran{\gammau{F}}.
 \]
 Consequently, using \eqref{Nr.RG}, we have
 \[
  \Ran{F(\iu y)}
  =\ran{\gammau{F}}+\Ran{\muuA{F}{\rhl }}
  =\ran{\gammau{F}}
 \]
 for each \(y\in(0,\infp)\). Hence, \(\rank F(\iu y)=\rank\gammau{F}\) hold for each \(y\in(0,\infp)\). Taking into account \rlem{LA*10}, we obtain then \(\lim_{y\to\infp}F^\mpi(\iu y)=\gammau{F}^\mpi\). Using this and \eqref{P8*14.4}, we get
 \[\begin{split}
  \Oqq
  &=-u_0+u_0
  =-A(\gammau{G}+A)^\mpi A+u_0
  =(-1)\cdot A\gammau{F}^\mpi A+u_0\\
  &=\rk*{\lim_{y\to\infp}\iu y\ek*{-(\iu y-\ug)^\inv}} A\ek*{\lim_{y\to\infp}F^\mpi(\iu y)} A+u_0\\
  &=\lim_{y\to\infp}\rk*{\iu y\ek*{-(\iu y-\ug)^\inv} AF^\mpi(\iu y)A+u_0}\\
  &=\lim_{y\to\infp}\iu y\ek*{-(\iu y-\ug)^\inv AF^\mpi(\iu y)A+(\iu y)^\inv u_0}\\
  &=\lim_{y\to\infp}\iu y\ek*{ \STiao{G}{A}(\iu y)+(\iu y)^\inv u_0}.
 \end{split}\]
 By virtue of \rcor{C-T6*4}, then the function \(\STiao{G}{A}\) belongs to \(\SFuqaa{0}{\seq{u_j}{j}{0}{0}}\). Because of \eqref{P8*14.1} we have \(\ran{\STiaoa{G}{A}{z}}=\ran{A}\), whereas \eqref{P8*14.4} yields \(\nul{A}\subseteq\nul{\STiaoa{G}{A}{z}}\). Hence, \rlem{LA*5} implies \(\nul{A}=\nul{\STiaoa{G}{A}{z}}\).
\eproof

 Now we indicate some generic situations in which the formulas in \eqref{F8*8} are satisfied. We start with the first formula in \eqref{F8*8}.

\bpropl{P8*17}
 Let \(\ug\in\R\), let \(F\in\SFqa\) and let \(A\in\Cqq\) be such that \(\ran{\gammau{F}}+\ran{\muua{F}{\rhl}}\subseteq\ran{A}\) and \(\nul{\gammau{F}}\cap\nul{\muua{F}{\rhl}}\subseteq\nul{A}\) hold. Then \(\STiao{\rk{\STao{F}{A}}}{A}=F\).
\eprop
\bproof
 Let \(z\in\Cs\). In view of \rprop{P2*7} the matrix \(X\defeq F(z)\) fulfills \(\ran{X}\subseteq\ran{A}\) and \(\nul{X}\subseteq\nul{A}\), which in view of \rlem{LA*5} implies \(\ran{X}=\ran{A}\) and \(\nul{X}=\nul{A}\). Hence, \rrem{RA*3} yields \(XX^\mpi=AA^\mpi\) and \(X^\mpi X=A^\mpi A\). With \(Y\defeq \STaoa{F}{A}{z}\) and \eqref{F8*1} we obtain then
  \[\begin{split}
   A^\mpi Y
   &=A^\mpi\rk*{-A\ek*{\Iq+(z-\ug)^\inv X^\mpi A}}\\
   &=-\frac{1}{z-\ug}A^\mpi AX^\mpi A-A^\mpi A
   =-\frac{1}{z-\ug}X^\mpi XX^\mpi A-A^\mpi A
   =-\frac{1}{z-\ug}X^\mpi A-A^\mpi A.
 \end{split}\]
 Furthermore, we have
 \[\begin{split}
  &\rk*{-\frac{1}{z-\ug}X^\mpi A+\Iq-A^\mpi A}\ek*{-(z-\ug)A^\mpi X+\Iq-A^\mpi A}\\
  &=X^\mpi AA^\mpi X-\frac{1}{z-\ug}X^\mpi A(\Iq-A^\mpi A)-(z-\ug)(\Iq-A^\mpi A)A^\mpi X+(\Iq-A^\mpi A)^2\\
  &=X^\mpi XX^\mpi X+\Iq-A^\mpi A
  =X^\mpi X+\Iq-X^\mpi X
  =\Iq.
 \end{split}\]
 In particular, \(\det\ek{-(z-\ug)^\inv X^\mpi A+\Iq-A^\mpi A}\neq0\) and \(\ek{-(z-\ug)^\inv X^\mpi A+\Iq-A^\mpi A}^\inv=-(z-\ug)A^\mpi X+\Iq-A^\mpi A\). Finally, we get
 \[\begin{split}
   -\frac{1}{z-\ug}A\rk{\Iq+A^\mpi Y}^\mpi
   &=-\frac{1}{z-\ug}A\rk*{\Iq-\frac{1}{z-\ug}X^\mpi A-A^\mpi A}^\mpi\\
   &=-\frac{1}{z-\ug}A\rk*{-\frac{1}{z-\ug}X^\mpi A+\Iq-A^\mpi A}^\inv\\
   &=-\frac{1}{z-\ug}A\ek*{-(z-\ug)A^\mpi X+\Iq-A^\mpi A}\\
   &=AA^\mpi X-\frac{1}{z-\ug}A(\Iq-A^\mpi A)
   =XX^\mpi X
   =X.
    \end{split}\]
  In view of \eqref{F8*1} and \eqref{F8*2}, the proof is complete.
\eproof

\bcorl{C8*18}
 Let $\kappa\in\NOinf $, let $\seqska \in\Kggeqka$, and let $F\in\SFqas{\kappa}$. Then \(\STiao{\rk{\STao{F}{\su{0}}}}{\su{0}}=F\).
\ecor
\bproof
 According to \rpart{P5*7.c} of \rprop{P5*7}, the function $F$ belongs to $\SFqa$ and \eqref{P5*7.A} holds true. Thus, the application of \rprop{P8*17} yields the assertion.
\eproof

 Now we turn our attention to the second formula in \eqref{F8*8}.

\bpropl{P8*19}
 Let \(\ug\in\R\), let \(A\in\Cggq\), and let \(G\in\SFqa\) be such that \eqref{F8*13} holds. Then \(\STao{(\STiao{G}{A})}{A}=G\).
\eprop
\bproof
 Let \(z\in\Cs\). According to \eqref{F8*1}, we have
 \beql{P8*19.2}
  \STao{(\STiao{G}{A})}{A}(z)
  =-A\rk*{\Iq+(z-\ug)^\inv\ek*{ \STiaoa{G}{A}{z}}^\mpi A}.
 \eeq
 For each \(w\in\Cs\), from \rlem{L8*10} we know that \(\det[\Iq+A^\mpi G(w)]\neq0\) and that
 \beql{P8*19.1}
  \STiao{G}{A}(w)
  =-(w-\ug)^\inv A\ek*{\Iq+A^\mpi G(w)}^\inv.
 \eeq
 From \rprop{P8*14} we know that \(\STiao{G}{A}\in\SFuqa{0}\). Because of this and \rcor{C0838}, we have \([\STiaoa{G}{A}{z}]^\ad=\STiao{G}{A}(\ko z)\). From \eqref{P8*19.1} we see that \(\ran{\STiao{G}{A}(\ko z)}=\ran{A}\) and, in view of \(A^\ad=A\), consequently \(\nul{\STiaoa{G}{A}{z}}=\ek{\ran{A}}^\oc=\nul{A}\). Hence, \rremp{RA*3}{RA*3.b} shows that
 \beql{P8*19.4}
  A\ek*{ \STiaoa{G}{A}{z}}^\mpi \STiaoa{G}{A}{z}
  =A.
 \eeq
 \rprop{P2*7} and \eqref{F8*13} provide us \eqref{L8*10.1} and, in view of \rremp{RA*3}{RA*3.b}, consequently \eqref{L8*10.2}. Using \eqref{P8*19.2}, \eqref{P8*19.1}, \eqref{P8*19.4}, and \eqref{L8*10.2}, we obtain
 \[
 \begin{split}
  &\STao{(\STiao{G}{A})}{A}(z)
  =-A+A\ek*{ \STiaoa{G}{A}{z}}^\mpi\ek*{-(z-\ug)^\inv A}\\
  &=-A+A\ek*{ \STiaoa{G}{A}{z}}^\mpi\rk*{-(z-\ug)^\inv A\ek*{\Iq+A^\mpi G(z)}^\inv}\ek*{\Iq+A^\mpi G(z)}\\
  &=-A+A\ek*{ \STiaoa{G}{A}{z}}^\mpi \STiaoa{G}{A}{z}\ek*{\Iq+A^\mpi G(z)}\\
  &=-A+A\ek*{\Iq+A^\mpi G(z)}
  =AA^\mpi G(z)
  =G(z).\qedhere
 \end{split}
 \]
\eproof

\bcorl{C8*21}
 Let \(\ug\in\R\), let $\kappa\in\NOinf$, let \(\seqska\in\Kggeqka\), and let \(G\in\SFqas{\kappa}\). Then \(\STao{\rk{\STiao{G}{\su{0}}}}{\su{0}}=G\).
\ecor
\bproof
 \rlem{L7*6} yields $\su{0}\in\Cggq$. Thus, taking \rpropp{P5*7}{P5*7.c} into account, the application of \rprop{P8*19} yields the assertion.
\eproof

\section{On the \taSSt{\su{0}} for the class $\SFqas{\kappa}$}\label{S*9}
 The central topic of this section can be described as follows. Let $m\in\NO$ and let $\seqs{m}\in\Kggequa{m}$. Then \rthm{T-P5*4} tells us that the class $\SFqas{m}$ is \tne{}. If $F\in\SFqas{m}$, then our interest is concentrated on the \taSSt{\su{0}} $\STao{F}{\su{0}}$ of $F$. We will obtain a complete description of this object. In the case $m=0$, we will show that $\STao{F}{\su{0}}$ belongs to $\SFdqaa{\su{0}}$ (see \rthm{T9*1}). Let us now consider the case $m\in\N$. If $\seq{t_j}{j}{0}{m-1}$ denotes the \saSchto{\seqs{m}}, then it will turn out (see \rthm{T9*7} and \rcor{C9*8}) that $\STao{F}{\su{0}}$ belongs to $\SFuqaa{m-1}{\seq{t_j}{j}{0}{m-1}}$. Our strategy to prove this is based on the application of Hamburger-Nevanlinna-type results for the class $\SFqa$, which were developed in \rsec{S*6}.

 Now we start with the detailed treatment of the case $m=0$.

\bthml{T9*1}
 Let \(\seqs{0}\in\Kggequa{0}\) and let \(F\in\SFqas{0}\). Then \(\STao{F}{\su{0}}\) belongs to \(\SFdqaa{\su{0}}\).
\ethm
\bproof
 Denote by \(\OSm{F}\) the \taSmo{F}. From \rrem{R1228} we infer
 \beql{T9*1.0}
  \su{0}
  =\OSmA{F}{\rhl}.
 \eeq
 According to \rprop{P0844}, the function \(G\colon\Cs  \to  \Cqq \) defined by \(G (z) \defeq  - (z-\ug)^\inv  F^\mpi (z)\) belongs to \(\SFqa\) and using \eqref{T9*1.0} we get \(\gammaG  =\ek{\OSma{F}{\rhl}}^\mpi=s_0^\mpi \). In view of \rlem{L7*6} we have \(\su{0}\in\Cggq\). Thus, \(s_0^\ad = s_0\). From this an \eqref{F8*1} we infer
 \beql{T9*1.1}
  \STao{F}{\su{0}} 
  = -s_0 + s_0^\ad Gs_0.
 \eeq
 Furthermore, \(H\defeq  s_0^\ad Gs_0\) belongs to \(\SFqa\) and \(\gammau{H} = s_0^\ad \gammaG  s_0\) by virtue of~\zitaa{142}{\crem{3.11}}. In view of \(\su{0}^\ad=\su{0}\) then \(\gammau{H}=\su{0}\). Because of \(\gammau{H} + (-s_0) =\Oqq\in \Cggq \) and \rrem{R0832}, then \(\STao{F}{\su{0}} \in\SFqa\) and \(\gammau{\STao{F}{\su{0}}} =\Oqq\). Hence, in view of \rprop{P2*5}, we obtain
 \[
  \lim_{y\to\infp}\STao{F}{\su{0}}(\iu y)
  =\gammau{\STao{F}{\su{0}}}
  =\Oqq.
 \]
 Thus, \eqref{F3*4} yields that \(\STao{F}{\su{0}}\) belongs to \(\SFdqa\). From \eqref{T9*1.1} we see moreover \(\STao{F}{\su{0}}\su{0}^\mpi\su{0}=\STao{F}{\su{0}}\). This finally shows \(\STao{F}{\su{0}}\in\SFdqaa{\su{0}}\).
\eproof

\bcorl{C9*2}
 Let \(\ug\in\R\), let $\kappa\in\NOinf$, let $\seqska\in\Kggeqka$, and let $F\in\SFqas{\kappa}$. Then \(\STao{F}{\su{0}}\in\SFdqaa{\su{0}}\).
\ecor
\bproof
 From \rrem{R5*3} we get $F\in\SFqas{0}$, which, in view of \rthm{T-P5*4}, implies \(\seqs{0}\in\Kggequa{0}\). Thus, the application of \rthm{T9*1} completes the proof.
\eproof

 Now we turn our attention to the case $m\in\N$.

\bthml{T9*7}
 Let \(\ug\in\R\), let \(m\in\N\), let \(\seqs{m}\in\Kggequa{m}\) with \saScht{} \(\seq{\tu{j}}{j}{0}{m-1}\), and let \(F\in\SFqas{m}\). Then \(\STao{F}{\su{0}}  \in\SFuqaa{m-1}{(\tu{j})_{j=0}^{m-1}}\).
\ethm
\begin{proof}
 The basic strategy of our proof is to apply \rcor{C-T6*4}. From \rthm{T-P7*2} we obtain \(\seq{\tu{j}}{j}{0}{m-1}\in\Kggequa{m-1}\). \rlem{L7*6} yields then
 \begin{align}\label{T9*7.1}
  \tu{j}&\in\CHq&\text{for all }j&\in\mn{0}{m-1}.
 \end{align}
 Furthermore, \(F\in\SFuqa{0}\) and the \taSm{} \(\OSm{F}\) of \(F\) belongs to \(\MggqKsg{m}\). \rrem{R1228} implies
 \begin{equation}\label{T9*7.2}
  s_0
  =\OSmA{F}{\rhl}.
 \end{equation}
 By virtue of \rcor{C9*2}, we have \(\STao{F}{\su{0}}\in\SFdqaa{\su{0}}\) and hence
 \beql{T9*7.3}
  \STao{F}{\su{0}} 
  \in\SFqa.
 \eeq
 In view of \(F\in\SFuqa{0}\), the application of \rpropss{P0844}{P2*5} and \eqref{T9*7.2} yields
 \beql{T9*7.4}
  \lim_{y\to\infp} (\iu y -\ug)^\inv  F^\mpi  (\iu y)
  = -s_0^\mpi.
 \eeq
 From \rcor{C6*3} we conclude that
 \[
  \lim_{y\to\infp} (\iu y)^{m+1} \ek*{ F(\iu y)+ \sum_{j=0}^m (\iu y)^{-(j+1)} s_j}
  = \Oqq 
 \]
 and, consequently, in view of \(\lim_{y\to\infp}(\iu y)^\inv=0\), then
 \begin{equation}\label{T9*7.5}
  \begin{split}
  \Oqq 
  &=\rk*{1-\ug \ek*{ \lim_{y\to\infp} (\iu y)^\inv }}\rk*{\lim_{y\to\infp}(\iu y)^{m+1} \ek*{ F(\iu y) + \sum_{j=0}^m (\iu y)^{-(j+1)}s_j}}\\
  &= \lim_{y\to\infp} (\iu y)^m (\iu y -\ug) \ek*{ F (\iu y) + \sum_{j=0}^m (\iu y)^{-(j+1)} s_j}.
 \end{split}
 \end{equation}
 Because of \rlem{L3*22} and \eqref{T9*7.2}, for each \(y\in (0,\infp)\), we have
 \[
  \Nul{F(\iu y)}
  \subseteq\Nul{\OSmA{F}{\rhl}}
  =\nul{\su{0}}
 \]
 and hence, in view of \rremp{RA*3}{RA*3.b}, then
 \begin{equation}\label{T9*7.6}
  s_0 F^\mpi  (\iu y) F (\iu y) 
  = s_0.
 \end{equation}
 Using \eqref{F8*1} and \eqref{T9*7.6}, for each \(y\in(0,\infp)\), we get
 \begin{equation}\label{T9*7.7}
  \begin{split}
  \STao{F}{\su{0}}  (\iu y)
  &=-\su{0}\ek*{\Iq+(\iu y-\ug)^\inv F^\mpi(z)\su{0}}\\
  &= -s_0 F^\mpi  (\iu y) F(\iu y) - (\iu y - \ug)^\inv  s_0 F^\mpi  (\iu y) s_0\\
  &= (\iu y-\ug)^\inv  s_0 F^\mpi  (\iu y) \ek*{ -s_0 - (\iu y-\ug) F (\iu y)}.
 \end{split}
 \end{equation}
 For each \(y\in (0,\infp)\) and each \(j\in\mn{0}{m-1}\), from \eqref{T9*7.6} and \eqref{F*sta} we get
 \[
 s_0 F^\mpi  (\iu y) F (\iu y) s_0^\mpi  \tu{j} 
 = s_0 s_0^\mpi  \tu{j} 
 = \tu{j}.
 \]
 Consequently, for each \(y\in (0,\infp)\), we have
 \begin{equation}\label{T9*7.8}
 s_0 F^\mpi  (\iu y) F(\iu y)  s_0^\mpi  \sum_{k=0}^{m-1} (\iu y)^{-(k+1)} \tu{k} 
 =\sum_{k=0}^{m-1} (\iu y)^{-(k+1)} \tu{k}.
 \end{equation}
 Taking into account \eqref{F-D1455}, we see that
 \begin{equation}\label{T9*7.9}
  \begin{split}
  -(\iu y-\ug) \sum_{j=0}^m (\iu y)^{-(j+1)}\su{j}
  &= - \ek*{ (\iu y)^{-0}s_0 + \sum_{k=1}^m  (\iu y)^{-k} (-\ug s_{k-1} + s_k) - (\iu y)^{-(m+1)} \ug s_m}\\
  &= -\sum_{j=0}^m (\iu y)^{-j} s_j^\spa  + (\iu y)^{-(m+1)} \ug s_m
 \end{split}
 \end{equation}
 and
 \[
  \sum_{j=0}^m (\iu y)^{-j} s_j^\spa  - s_0
  = \sum_{j=0}^{m-1} (\iu y)^{-(j+1)}  s_{j+1}^\spa 
 \]
 hold true for each \(y\in (0,\infp)\). Furthermore, for each \(y\in (0,\infp)\), we obtain
 \begin{equation}\label{T9*7.11}
  \begin{split}
   &-\sum_{j=0}^m \sum_{k=0}^{m-1}  (\iu y)^{-j-(k+1)}\su{j}^\spa\su{0}^\mpi\tu{k}
   =-\sum_{r=0}^{2m-1}\ek*{\sum_{\substack{k=0\\0\leq r-k\leq m}}^{m-1}(\iu y)^{-(r+1)}\su{r-k}^\spa\su{0}^\mpi\tu{k}}\\
   &=-\sum_{r=0}^{2m-1}\ek*{\sum_{\substack{k=0\\r\geq k\geq r-m}}^{m-1}(\iu y)^{-(r+1)}\su{r-k}^\spa\su{0}^\mpi\tu{k}}
   =-\sum_{r=0}^{2m-1}\sum_{k=\max\set{0,r-m}}^{\min\set{r,m-1}}(\iu y)^{-(r+1)}\su{r-k}^\spa\su{0}^\mpi\tu{k}\\
   &=-\sum_{r=0}^{m-1}\sum_{k=0}^r(\iu y)^{-(r+1)}\su{r-k}^\spa\su{0}^\mpi\tu{k}-\sum_{r=m}^{2m-1}\sum_{k=r-m}^{m-1}(\iu y)^{-(r+1)} \su{r-k}^\spa\su{0}^\mpi\tu{k}.
  \end{split}
 \end{equation}
 From \rprop{P-L7*9} we obtain \(\seqs{m}\in\Dqqu{m}\). Then~\zitaa{114arxiv}{\clem{7.11}} yields
 \beql{T9*7.12}
  \sum_{k=0}^j s_{j-k}^\spa  s_0^\mpi  \tu{k}
  = s_{j+1}^\spa 
 \eeq
 for each \(j\in\mn{0}{m-1}\). Combining \eqref{T9*7.11} and \eqref{T9*7.12}, we get
 \[
  -\sum_{j=0}^m\sum_{k=0}^{m-1}(\iu y)^{-j-(k+1)}\su{j}^\spa
  =-\sum_{j=0}^{m-1}(\iu y)^{-(j+1)}\su{j+1}^\spa-\sum_{j=m}^{2m-1}\sum_{k=j-m}^{m-1}(\iu y)^{-(j+1)}\su{j-k}^\spa\su{0}^\mpi\tu{k}.
 \]
 For each \(y\in (0,\infp)\), from \eqref{T9*7.7} and \eqref{T9*7.8} we get
 \begin{equation}\label{T9*7.13}
 \begin{split}
  &(\iu y)^m \ek*{\STao{F}{\su{0}}  (\iu y) +\sum_{k=0}^{m-1} (\iu y)^{-(k+1)} \tu{k}}\\
  &= (\iu y)^m \Biggl\{(\iu y-\ug)^\inv  s_0 F^\mpi  (\iu y)\ek*{-s_0 - (\iu y -\ug) F (\iu y)}\\
  &\hspace{107pt}+ (\iu y-\ug)^\inv  (\iu y-\ug) s_0 F^\mpi  (\iu y) F(\iu y) s_0^\mpi  \sum_{k=0}^{m-1}  (\iu y)^{-(k+1)} \tu{k}  \Biggr\}\\
  &= (\iu y)^m (\iu y-\ug)^\inv  s_0F^\mpi  (\iu y)\\
  &\qquad\times\ek*{ -s_0 - (\iu y-\ug )  F(\iu y) + (\iu y-\ug ) F(\iu y) s_0^\mpi  \sum_{k=0}^{m-1} (\iu y)^{-(k+1)} \tu{k}}.
 \end{split}
 \end{equation}
 Now we are going to derive an identity for the matrix in brackets on the right side of \eqref{T9*7.13}. For all \(z\in\ohe\), let\index{m@$M(z)$}
 \beql{T9*7.18}
  M(z)
  \defeq\Iq- s_0^\mpi  \sum_{k=0}^{m-1}z^{-(k+1)} \tu{k}.
 \eeq
 From \eqref{F-D1455} it follows \(\su{0}^\spa=\su{0}\). Combining this with \eqref{T9*7.18}, we obtain
 \begin{equation}\label{T9*7.19}
  \begin{split}
   &\rk*{\sum_{j=0}^m(\iu y)^{-j}\su{j}^\spa}\ek*{M(\iu y)}
   =\rk*{\sum_{j=0}^m(\iu y)^{-j}\su{j}^\spa}\ek*{\Iq-\su{0}^\mpi\sum_{k=0}^{m-1}(\iu y)^{-(k+1)}\tu{k}}\\
   &=\su{0}^\spa+\sum_{j=1}^m(\iu y)^{-j}\su{j}^\spa-\sum_{j=0}^m\sum_{k=0}^{m-1}(\iu y)^{-j-(k+1)}\su{j}^\spa\su{0}^\mpi\tu{k}\\
   &=\su{0}+\sum_{j=0}^{m-1}(\iu y)^{-(j+1)}\su{j+1}^\spa-\sum_{j=0}^{m-1}(\iu y)^{-(j+1)}\su{j+1}^\spa-\sum_{j=m}^{2m-1}\sum_{k=j-m}^{m-1}(\iu y)^{-(j+1)}\su{j-k}^\spa\su{0}^\mpi\tu{k}\\
   &=\su{0}-\sum_{j=m}^{2m-1}\sum_{k=j-m}^{m-1}(\iu y)^{-(j+1)}\su{j-k}^\spa\su{0}^\mpi\tu{k}.
  \end{split}
 \end{equation}
 Using \eqref{T9*7.18}, \eqref{T9*7.9}, \eqref{T9*7.19}, and \eqref{T9*7.12}, we infer then
 \begin{equation}\label{T9*7.14}
 \begin{split}
  &-s_0 - (\iu y-\ug )  F(\iu y) + (\iu y-\ug ) F(\iu y) s_0^\mpi  \sum_{k=0}^{m-1} (\iu y)^{-(k+1)} \tu{k}
  =\su{0}-\rk{\iu y-\ug}F(\iu y)M(\iu y)\\
  &=-s_0 -\rk*{(\iu y-\ug)\ek*{F (\iu y)  + \sum_{j=0}^{m} (\iu y)^{-(j+1)} s_j} -(\iu y-\ug)\sum_{j=0}^m (\iu y)^{-(j+1)} s_j}M(\iu y)\\
  &=-s_0\\
  &\qquad-\rk*{ (\iu y-\ug) \ek*{ F(\iu y) + \sum_{j=0}^{m} (\iu y)^{-(j+1)} s_j}-\sum_{j=0}^{m} (\iu y)^{-j} s_j^\spa  + (\iu y)^{-(m+1)}\ug s_m} M(\iu y)\\
  &=-\su{0}-\rk*{\rk{\iu y-\ug}\ek*{F(\iu y)+\sum_{j=0}^m(\iu y)^{-(j+1)}\su{j}}+\rk{\iu y}^{-(m+1)}\ug\su{m}}M(\iu y)\\
  &\qquad-\ek*{\sum_{j=0}^m (\iu y)^{-j}\su{j}^\spa}M(\iu y)\\
  &=-\su{0}-\rk*{\rk{\iu y-\ug}\ek*{F(\iu y)+\sum_{j=0}^m(\iu y)^{-(j+1)}\su{j}}+\rk{\iu y}^{-(m+1)}\ug\su{m}}M(\iu y)\\
  &\qquad+\su{0}-\sum_{j=m}^{2m-1}\sum_{k=j-m}^{m-1} (\iu y)^{-(j+1)}\su{j-k}^\spa\su{0}^\mpi\tu{k}\\
  &=-\sum_{j=m}^{2m-1}  \sum_{k= j-m}^{m-1}  (\iu y)^{-(j+1)} s_{j-k}^\spa  s_0^\mpi  \tu{k}\\
  &\qquad-\rk*{ (\iu y-\ug) \ek*{F(\iu y) + \sum_{j=0}^{m} (\iu y)^{-(j+1)} s_j} + (\iu y)^{-(m+1)}\ug s_m}M(\iu y)
 \end{split}
 \end{equation}
 for all \(y\in (0,\infp)\). Because of \eqref{T9*7.13}, \eqref{T9*7.14}, and \eqref{T9*7.18}, we have 
 \begin{equation}\label{T9*7.15}
  \begin{split}
  &(\iu y)^m \ek*{\STao{F}{\su{0}}  (\iu y) +\sum_{k=0}^{m-1} (\iu y)^{-(k+1)} \tu{k}}\\
  &=(\iu y)^m (\iu y-\ug)^\inv  s_0F^\mpi  (\iu y)\Biggl\{-\sum_{j=m}^{2m-1}  \sum_{k= j-m}^{m-1}  (\iu y)^{-(j+1)} s_{j-k}^\spa  s_0^\mpi  \tu{k}\\
  &\hspace{60pt}-\rk*{ (\iu y-\ug) \ek*{F(\iu y) + \sum_{j=0}^{m} (\iu y)^{-(j+1)} s_j} + (\iu y)^{-(m+1)}\ug s_m}M(\iu y)\Biggr\}\\
  &=s_0(\iu y-\ug )^\inv F^\mpi(\iu y)\Biggl\{-\sum_{j=m}^{2m-1}\ek*{(\iu y)^{m-j-1}\sum_{k=j-m}^{m-1} s_{j-k}^\spa  s_0^\mpi \tu{k}}\\
  &\hspace{60pt}-\rk*{(\iu y)^m(\iu y-\ug) \ek*{F(\iu y) + \sum_{j=0}^{m} (\iu y)^{-(j+1)} s_j} + (\iu y)^\inv\ug s_m}\\
  &\hspace{222pt}\times\ek*{\Iq  - s_0^\mpi  \sum_{k=0}^{m-1}  (\iu y)^{-(k+1)} \tu{k} }\Biggr\}
 \end{split}
 \end{equation}
 for all \(y\in (0,\infp)\). From \eqref{T9*7.4}, \eqref{T9*7.5}, and \eqref{T9*7.15} we conclude
 \begin{equation}\label{T9*7.16}
  \begin{split}
   &\Oqq
   =\su{0}(-\su{0}^\mpi)\Biggl\{-\sum_{j=m}^{2m-1}\rk*{0\cdot\sum_{k=j-m}^{m-1}\su{j-k}^\spa \su{0}^\mpi\tu{k}}\\
   &\hspace{62pt}-(\Oqq+0\cdot\ug\su{m})\rk*{\Iq-\su{0}^\mpi\sum_{k=0}^{m-1}0\cdot\tu{k}}\Biggr\}\\
   &=\su{0}\rk*{\lim_{y\to\infp}\ek*{(\iu y-\ug)^\inv F^\mpi(\iu y)}}\\
   &\qquad\times\Biggl\{-\sum_{j=m}^{2m-1}\rk*{\ek*{\lim_{y\to\infp}(\iu y)^{m-j-1}}\sum_{k=j-m}^{m-1}\su{j-k}^\spa \su{0}^\mpi\tu{k}}\\
   &\qquad\qquad-\ek*{\lim_{y\to\infp}\rk*{(\iu y)^m(\iu y-\ug)\ek*{ F(\iu y)+\sum_{j=0}^m(\iu y)^{-(j+1)}\su{j}}}+\ek*{\lim_{y\to\infp}(\iu y)^\inv}\ug\su{m}}\\
   &\hspace{213pt}\times\rk*{\Iq-\su{0}^\mpi\sum_{k=0}^{m-1}\ek*{\lim_{y\to\infp}(\iu y)^{-(k+1)}}\tu{k}}\Biggr\}\\
   &=\lim_{y\to\infp}\Biggl\{s_0(\iu y-\ug )^\inv F^\mpi(\iu y)\Biggl\{-\sum_{j=m}^{2m-1}\ek*{(\iu y)^{m-j-1}\sum_{k=j-m}^{m-1} s_{j-k}^\spa  s_0^\mpi \tu{k}}\\
   &\hspace{110pt}-\rk*{(\iu y)^m(\iu y-\ug) \ek*{F(\iu y) + \sum_{j=0}^{m} (\iu y)^{-(j+1)} s_j} + (\iu y)^\inv\ug s_m}\\
   &\hspace{252pt}\times\ek*{ \Iq  - s_0^\mpi  \sum_{k=0}^{m-1}  (\iu y)^{-(k+1)} \tu{k}  }\Biggr\}\Biggr\}\\
   &=\lim_{y\to\infp}(\iu y)^m \ek*{\STao{F}{\su{0}}  (\iu y) +\sum_{k=0}^{m-1} (\iu y)^{-(k+1)} \tu{k}}.
  \end{split}
 \end{equation}
 In view of \eqref{T9*7.1}, \eqref{T9*7.3}, and \eqref{T9*7.16}, then \rcor{C-T6*4} yields that \(\STao{F}{\su{0}} \) belongs to the class \(\SFuqaa{m-1}{\seq{\tu{j}}{j}{0}{m-1}}\).
\end{proof}

\bcorl{C9*8}
 Let \(\ug\in\R\), let $\seqsinf \in\Kggeqinfa$ with \saScht{} \(\seq{\tu{j}}{j}{0}{\infi}\), and let $F\in\SFqasinf$. Then $\STao{F}{\su{0}}$ belongs to $\SFuqaa{\infi}{\seq{\tu{j}}{j}{0}{\infi}}$.
\ecor
\bproof
 Combine \rremss{R5*3}{R7*1} and \rthm{T9*7}.
\eproof

\bpropl{P9*10}
 Let \(\ug\in\R\), let $\kappa\in\Ninf$, let $\seqska \in\Kggeqka$ with \saScht{} \(\seq{\tu{j}}{j}{0}{\kappa-1}\), and let $F\in\SFqaska$. Then $\STao{F}{\su{0}}$ belongs to $\SFuqaa{\kappa-1}{\seq{\tu{j}}{j}{0}{\kappa-1}}$.
\eprop
\bproof
 Because of $\seqska \in\Kggeqka$, we have $\seqs{m}\in\Kggequa{m}$ for all $m\in\mn{0}{\kappa}$. \rrem{R5*3} yields $F\in\bigcap_{m=0}^\kappa\SFqas{m}$. Thus, \rthm{T9*7} shows then that $\STao{F}{\su{0}}$ belongs to $\bigcap_{m=0}^{\kappa-1}\SFuqaa{m}{\seq{\tu{j}}{j}{0}{m}}$. Consequently, from \rrem{R5*3} then $\STao{F}{\su{0}}\in\SFuqaa{\kappa-1}{\seq{\tu{j}}{j}{0}{\kappa-1}}$ follows.
\eproof

\section{On the \tiaSSt{\su{0}} for special subclasses of $\SFqa$}\label{S*10}
 Against to the background of \rpropss{P8*17}{P8*19}, the considerations of \rsec{S*9} lead us to the study of two inverse problems which will be formulated now. Let $m\in\NO$, let $\seqs{m}\in\Kggequa{m}$, and let $F\in\SFqas{m}$. In the case $m=0$, it was shown in \rthm{T9*1} that $\STao{F}{\su{0}}$ belongs to $\SFdqaa{\su{0}}$. Now we start with a function $G\in\SFdqaa{\su{0}}$ and will show in \rprop{P10*1} that $\STiao{G}{\su{0}}\in\SFqas{0}$. Finally, we investigate the case $m\in\N$. Let $\seq{\tu{j}}{j}{0}{m-1}$ be the \saSchto{\seqska}, then we know from \rthm{T9*7} that $\STao{F}{\su{0}}\in\SFuqaa{m-1}{\seq{\tu{j}}{j}{0}{m-1}}$. So we are going to verify now that, for a function $G\in\SFuqaa{m-1}{\seq{\tu{j}}{j}{0}{m-1}}$, its \tiaSSt{\su{0}} $\STiao{G}{\su{0}}$ belongs to $\SFqas{m}$ (see \rthm{T10*9}).    

 First we turn our attention to a detailed treatment of the case $m=0$. An application of \rprop{P8*14} will provide us quickly the desired result.

\bpropl{P10*1}
 Let \(\ug\in\R\), let \(\seqs{0}\in\Kggequa{0}\), and let \(G\in\SFdqaa{\su{0}}\). Then \(\STiao{G}{\su{0}}\in\SFqas{0}\).
\eprop
\bproof
 Because of \(G\in\SFdqaa{\su{0}}\), we have \(G\in\SFdqa\). In view of \rrem{R3*5}, then \(G\in\SFqa\) and \(\gammau{G}=\Oqq\). From \rrem{R-C4*7} we obtain furthermore \(\ran{\muua{G}{\rhl}}\subseteq\ran{\su{0}}\). Thus, we get \(\ran{\gammau{G}}+\ran{\muua{G}{\rhl}}\subseteq\ran{\su{0}}\) and \(\su{0}=\su{0}(\gammau{G}+\su{0})^\mpi\su{0}\). The application of \rprop{P8*14} with \(A\defeq\su{0}\) provides us then \(\STiao{G}{\su{0}}\in\SFqas{0}\).
\eproof

\bcorl{C10*2}
 Let $\seqs{0}\in\Kggequa{0}$. For each $F\in\SFqas{0}$, let\index{s@$\HTua{\su{0}}{F}$}
 \[
  \HTua{\su{0}}{F}
  \defeq \STao{F}{\su{0}}.
 \]
 Then $\HTu{\su{0}}$ generates a bijective correspondence between $\SFqas{0}$ and $\SFdqaa{\su{0}}$. The inverse mapping $\HTu{\su{0}}^\inv$ is given for each $G\in\SFdqaa{\su{0}}$ by
 \(
  \HTu{\su{0}}^\inv(G)
  \defeq \STiao{G}{\su{0}}.
 \)
\ecor
\bproof
 Taking \rcor{C8*18} and \rcor{C8*21} into account, applying \rthm{T9*1} and \rprop{P10*1} completes the proof.
\eproof

 Now we turn our attention to the case $m\in\N$. Similar as in the proof of \rthm{T9*7}, we will use Hamburger-Nevanlinna-type results from \rsec{S*6}.

\bthml{T10*9}
 Let \(\ug\in\R\), let \(m\in\N\), let \(\seqs{m}\in\Kggeqka\) with \saScht{} \(\seq{\tu{j}}{j}{0}{m-1}\), and let \(G\in\SFuqaa{m-1}{\seq{\tu{j}}{j}{0}{m-1}}\). Then \(\STiao{G}{\su{0}}\) belongs to \(\SFqas{m}\).
\ethm
\bproof
 The basic strategy of our proof is to apply \rcor{C-T6*4}. \rlem{L7*6} yields
 \beql{T10*9.0}
  \su{0}
  \in\Cggq.
 \eeq
 From \rthm{T-P7*2} we get \(\seq{\tu{j}}{j}{0}{m-1}\in\Kggequa{m-1}\). According to \rpropp{P5*7}{P5*7.c} and \eqref{F*sta} we have then \(G\in\SFqa\) and \(\ran{\gammau{G}}+\ran{\muua{G}{\rhl}}=\ran{\tu{0}}\), whereas \eqref{F*sta} yields \(\ran{\tu{0}}\subseteq\ran{\su{0}}\). Thus
 \beql{T10*9.3}
  \ran{\gammau{G}}+\Ran{\muuA{G}{\rhl}}
  \subseteq\ran{\su{0}}.
 \eeq
 Consequently, in view of \eqref{T10*9.0} and \eqref{T10*9.3}, we infer from \rprop{P8*14}, that \(\STiao{G}{\su{0}}\in\SFuqaa{0}{\seq{u_j}{j}{0}{0}}\), where \(u_0\defeq\su{0}(\gammau{G}+\su{0})^\mpi\su{0}\). In particular,
 \begin{equation}\label{T10*9.1}
  \STiao{G}{\su{0}}
  \in\SFqa.
 \end{equation}
 The application of \rcor{C6*3} provides us
 \begin{align}\label{T10*9.7-8}
  \lim_{y\to\infp}(\iu y)^\ell\ek*{ G(\iu y)+\sum_{k=0}^{\ell-1}(\iu y)^{-(k+1)}\tu{k}}
  &=\Oqq
  &\text{for each }\ell&\in\mn{1}{m}.
 \end{align}
 In particular, using \eqref{T10*9.7-8} for \(\ell=1\), we get
 \begin{equation}\label{T10*9.9}
 \begin{split}
  \Iq
  &=\Iq+\su{0}^\mpi(0\cdot\Oqq-0\cdot\tu{0})\\
  &=\Iq+\su{0}^\mpi\Biggl\{\ek*{\lim_{y\to\infp}(\iu y)^\inv}\rk*{\lim_{y\to\infp}\iu y\ek*{ G(\iu y)+\sum_{j=0}^0(\iu y)^{-(j+1)}\tu{j}}}\\
  &\hspace{255pt}-\ek*{\lim_{y\to\infp}(\iu y)^\inv}\tu{0}\Biggr\}\\
  &=\lim_{y\to\infp}\ek*{\Iq+\su{0}^\mpi G(\iu y)}
  =\lim_{y\to\infp}N(\iu y),
 \end{split}
 \end{equation}
 where\index{n@$N(z)$}
 \beql{T10*9.7}
  N(z)
  \defeq\Iq+\su{0}^\mpi G(z)
 \eeq
 for all \(z\in\ohe\). From \eqref{T10*9.3}, \eqref{T10*9.7}, \eqref{FC*2}, \eqref{F*lft}, and \rlem{L8*10}, for each \(y\in(0,\infp)\), we obtain
 \beql{T10*9.14}
  \det N(\iu y)
  \neq0
 \eeq
 and
 \beql{T10*9.15}
  \STiao{G}{\su{0}}(\iu y)
  =-(\iu y-\ug)^\inv\su{0}\ek*{ N(\iu y)}^\inv.
 \eeq
 Because of \eqref{T10*9.14} and \eqref{T10*9.9}, we conclude that
 \beql{T10*9.11}
  \lim_{y\to\infp}\ek*{ N(\iu y)}^\inv
  =\Iq^\inv
  =\Iq.
 \eeq
 Using \eqref{F-D1455}, we obtain that
 \begin{equation}\label{T10*9.16}
  \begin{split}
   &(\iu y-\ug)\sum_{j=0}^m(\iu y)^{-(j+1)}\su{j}\\
   &=(\iu y)^0\su{0}+\sum_{j=1}^m(\iu y)^{-j}\su{j}-\ug\ek*{\sum_{j=0}^{m-1}(\iu y)^{-(j+1)}\su{j}}-\ug(\iu y)^{-(m+1)}\su{m}\\
   &=(\iu y)^{-0}\su{0}+\sum_{k=1}^m(\iu y)^{-k}(-\ug\su{k-1}+\su{k})-(\iu y)^{-(m+1)}\ug\su{m}\\
   &=\sum_{j=0}^m(\iu y)^{-j}\su{j}^\spa-(\iu y)^{-(m+1)}\ug\su{m}.
  \end{split}
 \end{equation}
 and, taking \(\su{0}^\spa=\su{0}\) into account,
 \begin{equation}\label{T10*9.17}
  -\su{0}+\sum_{j=0}^m(\iu y)^{-j}\su{j}^\spa
  =\sum_{j=0}^{m-1}(\iu y)^{-(j+1)}\su{j+1}^\spa
 \end{equation}
 are fulfilled for each \(y\in(0,\infp)\). Furthermore, for each \(y\in (0,\infp)\), we obtain
 \begin{equation}\label{T10*9.18}
  \begin{split}
   &-\sum_{j=0}^m\sum_{k=0}^{m-1}(\iu y)^{m-j-k}\su{j}^\spa\su{0}^\mpi\tu{k}
   =-\sum_{r=0}^{2m-1}\ek*{\sum_{\substack{k=0\\0\leq r-k\leq m}}^{m-1}(\iu y)^{m-r}\su{r-k}^\spa\su{0}^\mpi\tu{k}}\\
   &=-\sum_{r=0}^{2m-1}\ek*{\sum_{\substack{k=0\\r\geq k\geq r-m}}^{m-1}(\iu y)^{m-r}\su{r-k}^\spa\su{0}^\mpi\tu{k}}
   =-\sum_{r=0}^{2m-1}\sum_{k=\max\set{0,r-m}}^{\min\set{r,m-1}}(\iu y)^{m-r}\su{r-k}^\spa\su{0}^\mpi\tu{k}\\
   &=\sum_{r=0}^{m-1}\ek*{(\iu y)^{m-r}\rk*{-\sum_{k=0}^r\su{r-k}^\spa\su{0}^\mpi\tu{k}}}-\sum_{r=m}^{2m-1}\ek*{(\iu y)^{m-r}\sum_{k=r-m}^{m-1}\su{r-k}^\spa\su{0}^\mpi\tu{k}}.
  \end{split}
 \end{equation}
 From \rprop{P-L7*9} we get \(\seqs{m}\in\Dqqu{m}\). Hence,~\zitaa{114arxiv}{\clem{7.11}} provides us \(\sum_{l=0}^j\su{j-l}^\spa\su{0}^\mpi\tu{l}=\su{j+1}^\spa\) for all \(j\in\mn{0}{m-1}\). Consequently, for each \(y\in(0,\infp)\), we have
 \beql{T10*9.24}
  \sum_{j=0}^{m-1}(\iu y)^{m-j}\rk*{\su{j+1}^\spa-\sum_{l=0}^j\su{j-l}^\spa\su{0}^\mpi\tu{l}}
  =\Oqq.
 \eeq
 Using \eqref{T10*9.15} and \eqref{T10*9.14}, for each \(y\in(0,\infp)\), we infer
 \begin{equation}
 \begin{split}\label{T10*9.2}
  &(\iu y)^{m+1}\ek*{ \STiao{G}{\su{0}}(\iu y)+\sum_{j=0}^m(\iu y)^{-(j+1)}\su{j}}\\
  &=(\iu y)^{m+1}\rk*{-(\iu y-\ug)^\inv\su{0}\ek*{ N(\iu y)}^\inv+\sum_{j=0}^m(\iu y)^{-(j+1)}\su{j}}\\
  &=(\iu y-\ug)^\inv(\iu y)^{m+1}\rk*{-\su{0}+(\iu y-\ug)\ek*{\sum_{j=0}^m(\iu y)^{-(j+1)}\su{j}} N(\iu y)}\ek*{ N(\iu y)}^\inv.
 \end{split}
 \end{equation}
 From \eqref{T10*9.7}, \eqref{T10*9.16}, and \eqref{T10*9.17} we get
 \begin{equation}
 \begin{split}\label{T10*9.4}
  &-\su{0}+(\iu y-\ug)\ek*{\sum_{j=0}^m(\iu y)^{-(j+1)}\su{j}} N(\iu y)\\
  &=-\su{0}+(\iu y-\ug)\ek*{\sum_{j=0}^m(\iu y)^{-(j+1)}\su{j}}\ek*{\Iq+\su{0}^\mpi G(\iu y)}\\
  &=-\su{0}+\ek*{\sum_{j=0}^m(\iu y)^{-j}\su{j}^\spa-(\iu y)^{-(m+1)}\ug\su{m}}\\
  &\hspace{75pt}\times\rk*{\Iq+\su{0}^\mpi\ek*{ G(\iu y)+\sum_{k=0}^{m-1}(\iu y)^{-(k+1)}\tu{k}-\sum_{k=0}^{m-1}(\iu y)^{-(k+1)}\tu{k}}}\\
  &=-\su{0}+\sum_{j=0}^m(\iu y)^{-j}\su{j}^\spa-(\iu y)^{-(m+1)}\ug\su{m}\\
  &\qquad+\ek*{\sum_{j=0}^m(\iu y)^{-j}\su{j}^\spa-(\iu y)^{-(m+1)}\ug\su{m}}\su{0}^\mpi\\
  &\hspace{109pt}\times\rk*{\ek*{ G(\iu y)+\sum_{k=0}^{m-1}(\iu y)^{-(k+1)}\tu{k}}-\sum_{k=0}^{m-1}(\iu y)^{-(k+1)}\tu{k}}\\
  &=\sum_{j=0}^{m-1}(\iu y)^{-(j+1)}\su{j+1}^\spa-(\iu y)^{-(m+1)}\ug\su{m}\\
  &\qquad+\ek*{\sum_{j=0}^m(\iu y)^{-j}\su{j}^\spa-(\iu y)^{-(m+1)}\ug\su{m}}\su{0}^\mpi\ek*{ G(\iu y)+\sum_{k=0}^{m-1}(\iu y)^{-(k+1)}\tu{k}}\\
  &\qquad+\ek*{\sum_{j=0}^m(\iu y)^{-j}\su{j}^\spa-(\iu y)^{-(m+1)}\ug\su{m}}\su{0}^\mpi\ek*{-\sum_{k=0}^{m-1}(\iu y)^{-(k+1)}\tu{k}}
 \end{split}
 \end{equation}
 for all \(y\in(0,\infp)\). In view of \eqref{T10*9.18}, we get
 \begin{equation}\label{T10*9.8}
  \begin{split}
   &-\rk{\iu y}^{m+1}\ek*{\sum_{j=0}^m\rk{\iu y}^{-j}\su{j}^\spa-\rk{\iu y}^{-\rk{m+1}}\ug\su{m}}\su{0}^\mpi\ek*{\sum_{k=0}^{m-1}\rk{\iu y}^{-(k+1)}\tu{k}}\\
   &=-\ek*{\sum_{j=0}^m\rk{\iu y}^{m+1-j}\su{j}^\spa-\ug\su{m}}\su{0}^\mpi\ek*{\sum_{k=0}^{m-1}\rk{\iu y}^{-(k+1)}\tu{k}}\\
   &=-\sum_{j=0}^m\sum_{k=0}^{m-1}\rk{\iu y}^{m-j-k}\su{j}^\spa\su{0}^\mpi\tu{k}+\ug\su{m}\su{0}^\mpi\sum_{k=0}^{m-1}\rk{\iu y}^{-(k+1)}\tu{k}\\
   &=\sum_{j=0}^{m-1}\ek*{\rk{\iu y}^{m-j}\rk*{-\sum_{\ell=0}^j\su{j-\ell}^\spa\su{0}^\mpi\tu{\ell}}}-\sum_{j=m}^{2m-1}\ek*{\rk{\iu y}^{m-j}\sum_{\ell=j-m}^{m-1}\su{j-\ell}^\spa\su{0}^\mpi\tu{\ell}}\\
   &\qquad+\ug\su{m}\su{0}^\mpi\tu{\ell}\sum_{k=0}^{m-1}\rk{\iu y}^{-(k+1)}\tu{k}.
  \end{split}
 \end{equation}
 Because of \eqref{T10*9.4}, \eqref{T10*9.8}, and \eqref{T10*9.24}, we obtain
 \begin{equation}
 \begin{split}\label{T10*9.5}
  &(\iu y)^{m+1}\rk*{-\su{0}+(\iu y-\ug)\ek*{\sum_{j=0}^m(\iu y)^{-(j+1)}\su{j}} N(\iu y)}\\
  &=\sum_{j=0}^{m-1}(\iu y)^{m-j}\su{j+1}^\spa-\ug\su{m}\\
  &\qquad+\iu y\ek*{\sum_{j=0}^m(\iu y)^{-j}\su{j}^\spa-(\iu y)^{-(m+1)}\ug\su{m}}\su{0}^\mpi\rk*{(\iu y)^m\ek*{ G(\iu y)+\sum_{k=0}^{m-1}(\iu y)^{-(k+1)}\tu{k}}}\\
  &\qquad-(\iu y)^{m+1}\ek*{\sum_{j=0}^m(\iu y)^{-j}\su{j}^\spa-(\iu y)^{-(m+1)}\ug\su{m}}\su{0}^\mpi\ek*{\sum_{k=0}^{m-1}(\iu y)^{-(k+1)}\tu{k}}\\
  &=\sum_{j=0}^{m-1}(\iu y)^{m-j}\su{j+1}^\spa\\
  &\qquad+\iu y\ek*{\sum_{j=0}^m(\iu y)^{-j}\su{j}^\spa-(\iu y)^{-(m+1)}\ug\su{m}}\su{0}^\mpi\rk*{(\iu y)^m\ek*{ G(\iu y)+\sum_{k=0}^{m-1}(\iu y)^{-(k+1)}\tu{k}}}\\
  &\qquad+\sum_{j=0}^{m-1}\ek*{(\iu y)^{m-j}\rk*{-\sum_{l=0}^j\su{j-l}^\spa\su{0}^\mpi\tu{l}}}-\sum_{j=m}^{2m-1}\ek*{(\iu y)^{m-j}\sum_{l=j-m}^{m-1}\su{j-l}^\spa\su{0}^\mpi\tu{l}}\\
  &\qquad+\ug\su{m}\su{0}^\mpi\sum_{k=0}^{m-1}(\iu y)^{-(k+1)}\tu{k}-\ug\su{m}\\
  &=\sum_{j=0}^{m-1}(\iu y)^{m-j}\rk*{\su{j+1}^\spa-\sum_{l=0}^j\su{j-l}^\spa\su{0}^\mpi\tu{l}}\\
  &\qquad+\iu y\ek*{\sum_{j=0}^m(\iu y)^{-j}\su{j}^\spa-(\iu y)^{-(m+1)}\ug\su{m}}\su{0}^\mpi\rk*{(\iu y)^m\ek*{ G(\iu y)+\sum_{k=0}^{m-1}(\iu y)^{-(k+1)}\tu{k}}}\\
  &\qquad-\sum_{j=m}^{2m-1}\ek*{(\iu y)^{m-j}\sum_{l=j-m}^{m-1}\su{j-l}^\spa\su{0}^\mpi\tu{l}}+\ug\su{m}\ek*{\su{0}^\mpi\sum_{k=0}^{m-1}(\iu y)^{-(k+1)}\tu{k}-\Iq}\\
  &=\iu y\ek*{\sum_{j=0}^m(\iu y)^{-j}\su{j}^\spa-(\iu y)^{-(m+1)}\ug\su{m}}\su{0}^\mpi\rk*{(\iu y)^m\ek*{ G(\iu y)+\sum_{k=0}^{m-1}(\iu y)^{-(k+1)}\tu{k}}}\\
  &\qquad-\sum_{j=m}^{2m-1}\ek*{(\iu y)^{m-j}\sum_{l=j-m}^{m-1}\su{j-l}^\spa\su{0}^\mpi\tu{l}}+\ug\su{m}\ek*{\su{0}^\mpi\sum_{k=0}^{m-1}(\iu y)^{-(k+1)}\tu{k}-\Iq}
 \end{split}
 \end{equation}
 for all \(y\in(0,\infp)\). In view of \eqref{T10*9.2} and \eqref{T10*9.5}, we obtain
 \begin{equation}\label{T10*9.29}
 \begin{split}
  &(\iu y)^{m+1}\ek*{ \STiao{G}{\su{0}}(\iu y)+\sum_{j=0}^m(\iu y)^{-(j+1)}\su{j}}\\
  &=\Biggl\{(\iu y-\ug)^\inv\iu y\ek*{\sum_{j=0}^m(\iu y)^{-j}\su{j}^\spa-(\iu y)^{-(m+1)}\ug\su{m}}\su{0}^\mpi\\
  &\hspace{125pt}\times\rk*{(\iu y)^m\ek*{ G(\iu y)+\sum_{k=0}^{m-1}(\iu y)^{-(k+1)}\tu{k}}}\\
  &\qquad-\sum_{j=m}^{2m-1}\ek*{(\iu y-\ug)^\inv(\iu y)^{m-j}\sum_{l=j-m}^{m-1}\su{j-l}^\spa\su{0}^\mpi\tu{l}}\\
  &\hspace{101pt}+(\iu y-\ug)^\inv\ug\su{m}\ek*{\su{0}^\mpi\sum_{k=0}^{m-1}(\iu y)^{-(k+1)}\tu{k}-\Iq}\Biggr\}\ek*{ N(\iu y)}^\inv
 \end{split}
 \end{equation}
 for all \(y\in(0,\infp)\). Using \eqref{T10*9.7-8}, \eqref{T10*9.11}, and \eqref{T10*9.29}, we get
 \begin{equation}\label{T10*9.6}
  \begin{split}
   &\Oqq
   =\Oqq\cdot\Iq\\
   &=\Biggl\{1\cdot\rk*{\su{0}^\spa+\sum_{j=1}^m0\cdot\su{j}^\spa-0\cdot\ug\su{m}}\su{0}^\mpi\cdot\Oqq\\
   &\qquad-\sum_{j=m}^{2m-1}\rk*{0\cdot\sum_{\ell=j-m}^{m-1}\su{j-\ell}^\spa\su{0}^\mpi\tu{\ell}}+0\cdot\ug\su{m}\rk*{\su{0}^\mpi\sum_{k=0}^{m-1}0\cdot\tu{k}-\Iq}\Biggr\}\cdot\Iq\\
   &=\Biggl\{\ek*{\lim_{y\to\infp}\rk{y-\ug}^\inv y}\rk*{\su{0}^\spa+\sum_{j=1}^m\ek*{\lim_{y\to\infp}(\iu y)^{-j}}\su{j}^\spa-\ek*{\lim_{y\to\infp}(\iu y)^{-(m+1)}}\ug\su{m}}\su{0}^\mpi\\
   &\qquad\qquad\times\rk*{\lim_{y\to\infp}(\iu y)^m\ek*{G(\iu y)+\sum_{k=0}^{m-1}(\iu y)^{-(k+1)}\tu{k}}}\\
   &\qquad-\sum_{j=m}^{2m-1}\ek*{\rk*{\lim_{y\to\infp}\ek{(y-\ug)^\inv(\iu y)^{m-j}}}\sum_{\ell=j-m}^{m-1}\su{j-\ell}^\spa\su{0}^\mpi\tu{\ell}}\\
   &\qquad+\ek*{\lim_{y\to\infp}(y-\ug)^\inv}\ug\su{m}\rk*{\su{0}^\mpi\sum_{k=0}^{m-1}\ek*{\lim_{y\to\infp}(\iu y)^{-(k+1)}}\tu{k}-\Iq}\Biggr\}\rk*{\lim_{y\to\infp}\ek*{N(\iu y)}^\inv}\\
   &=\lim_{y\to\infp}(\iu y)^{m+1}\ek*{ \STiao{G}{\su{0}}(\iu y)+\sum_{j=0}^m(\iu y)^{-(j+1)}\su{j}}.
  \end{split}
 \end{equation}
 Since \(\seqs{m}\) is, in view of \rlem{L7*6}, a sequence of \tH{} matrices and because of \eqref{T10*9.1} and \eqref{T10*9.6}, the application of \rcor{C-T6*4} yields \(\STiao{G}{\su{0}}\in\SFqas{m}\).
\eproof

\bcorl{C10*10}
 Let \(\ug\in\R\), let $\seqsinf\in\Kggeqinfa$ with \saScht{} \(\seq{\tu{j}}{j}{0}{\infi}\), and let $G\in\SFuqaa{\infi}{\seq{\tu{j}}{j}{0}{\infi}}$. Then the function $\STiao{G}{\su{0}}$ belongs to $\SFqasinf$.
\ecor
\bproof
 Combine \rremss{R5*3}{R7*1} with \rthm{T10*9}.
\eproof

\bcorl{C10*12}
 Let \(\ug\in\R\), let $\kappa\in\Ninf$, let $\seqs{\kappa}\in\Kggequa{\kappa}$ with \saScht{} $\seq{\tu{j}}{j}{0}{\kappa-1}$. Denote by $\HTu{\su{0}}$ the bijective mapping defined in \rcor{C10*2}. Then $\HTu{\su{0}}$ generates a bijective correspondence between the sets $\SFqas{\kappa}$ and $\SFuqaa{\kappa-1}{\seq{\tu{j}}{j}{0}{\kappa-1}}$.  The inverse mapping $\HTu{\su{0}}^\inv$ is given for $G\in\SFuqaa{\kappa-1}{\seq{\tu{j}}{j}{0}{\kappa-1}}$ by $\HTu{\su{0}}^\inv(G)=\STiao{G}{\su{0}}$.
\ecor
\bproof
 We will consider the cases $\kappa\in\N$ and $\kappa=\infty$. In any of these cases, we can apply \rcorss{C8*18}{C8*21} to verify the shape of the inverse mapping. If $\kappa$ belongs to \(\N\), then the assertion follows from \rthmss{T9*7}{T10*9}. Finally, if $\kappa=\infty$, then one has to apply \rcorss{C9*8}{C10*10}.
\eproof

 We mention that the investigations in \rsecss{S*9}{S*10} were influenced to some extent by considerations in Chen/Hu~\zita{MR2038751}. In particular,~\zitaa{MR2038751}{\clem{2.3}} played an essential role for the choice of our strategy. This concerns first the development of a Schur-type algorithm for sequences of complex matrices and then the construction of an interrelated Schur-type algorithm for functions belonging to special subclasses of $\SFqa$. In~\zitas{MR1807884,MR2038751} Chen and Hu treat the truncated matricial Stieltjes moment problem (\(\alpha=0\) in our setting). They introduce a transformation \(\Gamma_m\) via~\zitaa{MR1807884}{formula~(9)} (see also~\zitaa{MR2038751}{formula~(3.2)}) which maps a sequence of length \(m+1\) of complex square matrices to a sequence of length \(m\) of complex square matrices. This transformation essentially coincides for \(\alpha=0\) and sequences from \(\Dqqu{m}\) with the first \taScht{}ation given via \eqref{F*sta}. To describe the respective solution sets, Chen and Hu then reduce the length of the given sequence of prescribed moments in each step by~2, using the transformation \(\Gamma_{m-1}\Gamma_m\) (see~\zitaa{MR1807884}{formula~(12)} and~\zitaa{MR2038751}{formula~(3.7)}). To prove~\zitaa{MR2038751}{\cthm{3.1}}, they use \clem{2.6} from their paper~\zita{MR1624548} on the truncated matricial Hamburger moment problem.

\section{A \Schur{}-Stieltjes-type algorithm for the class $\SFqa$}\label{S*11}
 The results of \rsec{S*9} suggest the construction of a Schur-Nevanlinna-type algorithm for the class $\SFkaqa$ with $\kappa\in\NOinf $. The main theme of this section is to work out the details of this algorithm. In \rsec{S*9}, we fixed an $m\in\NO$ and a sequence $\seqs{m}\in\Kggequa{m}$. \rthm{T-P5*4} tells us then that $\SFqas{m}\neq\emptyset$. Let $F\in\SFqas{m}$. Then the central theme of \rsec{S*9} was to study the \taSSt{\su{0}} $\STao{F}{\su{0}}$ of $F$. The following observation shows that in the case $\kappa\in\NOinf $ the results of \rsec{S*9} can be applied to arbitrary functions belonging to $\SFkaqa$.

\bleml{L11*1}
 Let \(\ug\in\R\), let $\kappa\in\NOinf$, and let $F\in\SFkaqa$ with \taSm{} \(\OSm{F}\). Then $\OSm{F}$ belongs to $\MgguqK{\kappa}$,  the sequence $\seqska$ given by \(\su{j}\defeq\int_\rhl x^j\OSm{F}(\dif x)\) belongs to \(\Kggeqka\) and \(F\) belongs to \(\SFqaska\).
\elem
\bproof
 From \eqref{F3*11} we see that $\OSm{F}\in\MgguqK{\kappa}$. Thus, from \eqref{F5*3} we infer now $F\in\SFqaska$. Hence, \rthm{T-P5*4} yields $\seqska\in\Kggeqka$.
\eproof

 In view of \rlem{L11*1}, we introduce the following construction of \tSSt{s} for functions belonging to the class $\SFkaqa$ with $\kappa\in\NOinf$.

\bdefnl{D11*2}
 Let \(\ug\in\R\), let \(\kappa\in\NOinf\), and let \(F\in\SFkaqa\) with \taSm{} \(\OSm{F}\). Then the function \(\STao{F}{\OSm{F}(\rhl)}\) is called the \notii{\tSSto{F}}.
\edefn

\begin{rem}\label{R1546}
 Let \(\ug\in\R\), let \(\kappa\in\NOinf\), let \(\seqska\in\Kggeqka\), and let \(F\in\SFqas{\kappa}\). Then \(F\in\SFkaqa\) and \(\STao{F}{\su{0}}\) is the \tSSto{F}.
\end{rem}

\begin{prop}\label{P11*4}
 Let \(\ug\in\R\), let \(\kappa\in\NOinf\), let \(\seqska\in\Kggeqka\), and let \(F\in \SFqas{\kappa}\). Then the \tSSto{F} belongs to \(\SFdqaa{\su{0}}\).
\end{prop}
\begin{proof}
 In view of \rrem{R1546}, the application of \rcor{C9*2} yields the assertion.
\end{proof}

\bleml{L11*6}
 Let \(\ug\in\R\), let \(m\in\N\), and let \(F\in\SFuqa{0}\). Furthermore, let \(G\) be the \tSSto{F}. Then \(G\) belongs to \(\SFuqa{m-1}\) if and only if \(F\) belongs to \(\SFuqa{m}\).
\elem
\bproof
 Let \(\OSm{F}\) be the \taSmo{F}. In view of \rdefn{D11*2}, we have
 \beql{L11*6.A}
  G
  =\STao{F}{\OSm{F}(\rhl)}.
 \eeq
 First suppose that
 \beql{L11*6.B}
  G
  \in\SFuqa{m-1}.
 \eeq
 In view of \eqref{L11*6.B} and \eqref{F5*3}, then \(G\in\SFuqa{0}\) and the \taSm{} \(\OSm{G}\) of \(G\) satisfies \(\OSm{G}\in\MggquK{m-1}\). Hence, the matrices
 \begin{align}\label{L11*6.1}
   \tu{j}&\defeq\int_\rhl x^j\OSma{G}{\dif x},&
   j&\in\mn{0}{m-1}
 \end{align}
 are well defined and, in view of \eqref{L11*6.B} and \eqref{L11*6.1}, we have
 \beql{L11*6.2}
   G
   \in\SFuqAA{m-1}{\seq{t_j}{j}{0}{m-1}}.
 \eeq
 In view of \eqref{L11*6.2}, then \(\SFuqaa{m-1}{\seq{t_j}{j}{0}{m-1}}\neq\emptyset\). Thus, \rthm{T-P5*4} implies \(\seq{\tu{j}}{j}{0}{m-1}\in\Kggequa{m-1}\). The matrix \(A\defeq\OSma{F}{\rhl}\) is obviously \tnnH{}. Let \(\seqs{m}\) be the \siaSchto{\seq{t_j}{j}{0}{m-1}}{A}. From \rdefn{114.D1712} then we get \(\su{0}=\OSma{F}{\rhl}\). Using \(\seq{\tu{j}}{j}{0}{m-1}\in\Kggequa{m-1}\) and \(A\in\Cggq\), we conclude from~\zitaa{114arxiv}{\cprop{10.15}} that \(\seqs{m}\in\Kggequa{m}\) and, in particular, \(\seqs{0}\in\Kggequa{0}\). In view of \(F\in\SFuqa{0}\) and \(\su{0}=\OSma{F}{\rhl}\), we infer from \rremp{R1530}{R1530.a} that \(F\in\SFqas{0}\). Because of \(\su{0}=\OSma{F}{\rhl}\), \(F\in\SFqas{0}\), \(\seqs{0}\in\Kggequa{0}\), and \rcor{C8*18}, we obtain
 \beql{L11*6.3}
  \STiao{G}{\su{0}}
  =\STiao{\rk*{\STao{F}{\OSma{F}{\rhl}}}}{\su{0}}
  =\STiao{\rk*{\STao{F}{\su{0}}}}{\su{0}}
  =F.
 \eeq
 Taking \eqref{F8*1} and \eqref{L11*6.A} into account, we get
 \begin{equation}\label{L11*6.3A}
  \nul{A}
  \subseteq\Nul{-A\rk*{\Iq+(\iu-\ug)^\inv\ek*{F(\iu)}^\mpi A}}
  =\Nul{\STaoa{F}{A}{\iu}}
  =\Nul{G (\iu)}.
 \end{equation}
 Furthermore, in view of \eqref{L11*6.B}, we infer from \rlem{L3*22} that \(\nul{G (\iu)}=\nul{\OSma{G}{\rhl}}\). Because of \eqref{L11*6.2}, we see from \rrem{R1228} that
 \begin{equation}\label{L11*6.3C}
  \tu{0}
  =\OSmA{G}{\rhl}.
 \end{equation}
 The combination of \eqref{L11*6.3A}, \(\nul{G (\iu)}=\nul{\OSma{G}{\rhl}}\), and \eqref{L11*6.3C} yields \(\nul{A}\subseteq\nul{\tu{0}}\). From \eqref{L11*6.A} and \eqref{F8*1}, we conclude that
 \begin{equation}\label{L11*6.4A}
  \Ran{G (\iu)}
  =\Ran{\STaoa{F}{\su{0}}{\iu}}
  =\Ran{-A\rk*{\Iq+(\iu-\ug)^\inv\ek*{F(\iu)}^\mpi A}}
  \subseteq\ran{A}.
 \end{equation}
 Taking \eqref{L11*6.A} into account, we infer from \rlem{L3*22} that \(\ran{G (\iu)}=\ran{\OSma{G}{\rhl}}\). In combination with \eqref{L11*6.3C} and \eqref{L11*6.4A}, this implies \(\ran{\tu{0}}\subseteq\ran{A}\). By virtue of \eqref{L11*6.3C}, \eqref{L11*6.3}, and~\zitaa{MR2988005}{\clem{B.2}}, for each \(j\in\mn{0}{m-1}\), we conclude
 \[
   \nul{t_0}
   =\Nul{\OSmA{G }{\rhl}}
   \subseteq\Nul{\int_\rhl x^j\OSma{G}{\dif x}}
   =\nul{t_j}
 \]
 and
 \[
   \ran{t_j}
   =\Ran{\int_\rhl x^j\OSma{G}{\dif x}}
   \subseteq\Ran{\OSmA{G }{\rhl}}
   =\ran{t_0}.
 \]
 Hence in view of \rdefn{D7*8}, the sequence \(\seq{t_j}{j}{0}{m-1}\) belongs to \(\Dqqu{m-1}\). Taking additionally into account \(\nul{A}\subseteq\nul{\tu{0}}\) and \(\ran{\tu{0}}\subseteq\ran{A}\), we obtain from~\zitaa{114arxiv}{\cprop{10.8}} that \(\seq{t_j}{j}{0}{m-1}\) is exactly the \saScht{} \(\seq{\su{j}^\sta{1}}{j}{0}{m-1}\) of \(\seqs{m}\). In view of \eqref{L11*6.1}, this implies \(G \in\SFuqaa{m-1}{\seq{\su{j}^\sta{1}}{j}{0}{m-1}}\). The application of \rthm{T10*9} provides us then that \(\STiao{G}{\su{0}}\) belongs to the class \(\SFqas{m}\). Thus, \eqref{L11*6.2} implies \(F\in\SFqas{m}\) and, consequently, \(F\in\SFuqa{m}\).
 
 Conversely, now suppose that \(F\) belongs to \(\SFuqa{m}\). For each \(j\in\mn{0}{m}\), then the matrix \(\su{j}\defeq\int_\rhl x^j\OSma{F}{\dif x}\) is well defined and \(F\) belongs to \(\SFqas{m}\), which, in view of \rthm{T-P5*4}, implies \(\seqs{m}\in\Kggequa{m}\). Thus, if \(\seq{\tu{j}}{j}{0}{m-1}\) denotes the \saSchto{\seqs{m}}, then \rthm{T9*7} and \rrem{R1546} yield that \(G\) belongs to \(\SFuqaa{m-1}{\seq{\tu{j}}{j}{0}{m-1}}\). In particular, \(G\) belongs to \(\SFuqa{m-1}\).
\eproof

 In view of \rlem{L11*6}, we introduce in recursive way the following notions, which will play a central role in our further investigations.

\bdefnl{D11*8}
 Let \(\ug\in\R\), let \(\kappa\in\NOinf \), and let \(F\in\SFuqa{\kappa}\). We will write \(F^\Saalpha{0}\) for \(F\) and call it the \notii{\tthSSto{0}{F}}. For each \(k\in\mn{1}{\kappa+1}\), we will call the \tSSt{} \(F^\Saalpha{k}\)\index{\(F^\Saalpha{k}\)} of the \tthSSt{(k-1)} \(F^\Saalpha{k-1}\) the \notii{\tthSSto{k}{F}}.
\edefn

\breml{R11*9}
 Let \(\ug\in\R\), let \(\kappa\in\NOinf \), let \(F\in\SFuqa{\kappa}\), and let \(k\in\mn{0}{\kappa+1}\). Because of \rdefn{D11*8} and \rlem{L11*6}, then \(F^\Saalpha{k}\in\SFuqa{\kappa-k}\). Furthermore, for each \(\ell\in\mn{0}{\kappa-k+1}\), we get \((F^\Saalpha{k})^\Saalpha{\ell}=F^\Saalpha{k+\ell}\).
\erem

 The content of our next considerations can be described as follows. Let $\kappa\in\Ninf$ and let $\seqska \in\Kggeqka$. Then we are going to study the Schur-Stieltjes-type algorithm, which is generated via \rdefn{D11*8}, particularly for functions which belong to the class $\SFqas{\kappa}$. For every choice of $F\in\SFqa$ and $\su{0}\in\Cqq$, we set\index{s@$\HTua{\su{0}}{F}$}\index{s@$\HTiua{\su{0}}{F}$}
\begin{align*}
 \HTua{\su{0}}{F}&\defeq \STao{F}{\su{0}}&
 &\text{and}&
 \HTiua{\su{0}}{F}&\defeq \STiao{F}{\su{0}}.
\end{align*}

\bpropl{P11*10}
 Let \(\ug\in\R\), let $\kappa\in\NOinf $, let $\seqska \in\Kggeqka$, and let $F$ belong to $\SFqas{\kappa} $. Further, let \(m\in\mn{0}{\kappa}\) and let $F^\Saalpha{m+1}$ be the \tthSSto{(m+1)}{F}. For each $\ell \in\mn{0}{\kappa}$, let $\seq{\su{j}^\sta{\ell}}{j}{0}{\kappa-\ell }$ be the \sthaSchto{\ell}{\seqska}. Then:
 \benui
  \il{P11*10.a}
  \(
   F^\Saalpha{m+1}
   =\rk{\HTu{\su{0}^\sta{m}}\circ\HTu{\su{0}^\sta{m-1}}\circ\dotsb\circ\HTu{\su{0}^\sta{1}}\circ\HTu{\su{0}^\sta{0}}}(F)
  \)
  and
  \[
   F^\Saalpha{m+1}
   \in
   \begin{cases}
    \SFuqaa{\kappa-(m+1)}{\seq{\su{j}^\sta{m+1}}{j}{0}{\kappa-(m+1)}}\incase{m<\kappa}\\
    \SFdqaa{\su{0}^\sta{m}}\incase{m=\kappa}
   \end{cases}.
  \]
  \il{P11*10.b} $F=\rk{\HTiu{\su{0}^\sta{0}}\circ\HTiu{\su{0}^\sta{1}}\circ\dotsb\circ\HTiu{\su{0}^\sta{m}}}(F^\Saalpha{m+1})$.
 \eenui
\eprop
\bproof
 \eqref{P11*10.a} In view of the construction of $F^\Saalpha{m+1}$ (see \rdefnss{D11*8}{D11*2}) and \rrem{R1546}, this follows by induction, using \rthm{T9*7} and \rcorss{C9*2}{C9*8}.
 
 \eqref{P11*10.b} In view of \rcorss{C10*2}{C10*12}, \rpart{P11*10.b} is an immediate consequence of~\eqref{P11*10.a}.
\eproof

 Our next considerations are aimed to study the inversion of the Schur-Stieltjes algorithm. The following result can be considered as an inverse statement with respect to \rprop{P11*10}.

\bpropl{P11*12}
 Let \(\ug\in\R\), let $m\in\NO$, and let $\seqs{m}\in\Kggequa{m}$. For each $\ell\in\mn{0}{m}$, let $\seq{\su{j}^\sta{\ell}}{j}{0}{m-\ell}$ be the \sthaSchto{\ell}{\seqs{m}}. Further, let $G\in\SFdqaa{\su{0}^\sta{m}}$ and let
 \[
  F
  \defeq\rk{\HTiu{\su{0}^\sta{0}}\circ\HTiu{\su{0}^\sta{1}}\circ\dotsb\circ\HTiu{\su{0}^\sta{m}}}(G).
 \]
 Then the following statements hold true:
 \benui
  \il{P11*12.a} $F\in\SFqas{m}$.
  \il{P11*12.b} The \tthSSt{(m+1)} $F^\Saalpha{m+1}$ of $F$ satisfies $F^\Saalpha{m+1}=G$.
 \eenui
\eprop
\bproof
 \eqref{P11*12.a} This follows by combining \rprop{P10*1} and \rthm{T10*9}.
 
 \eqref{P11*12.b} Taking the definition of $F$ and~\eqref{P11*12.a} into account, the application of \rcorss{C10*2}{C10*12} yields
 \[
  G=\rk{\HTu{\su{0}^\sta{m}}\circ\HTu{\su{0}^\sta{m-1}}\circ\dotsb\circ\HTu{\su{0}^\sta{1}}\circ\HTu{\su{0}^\sta{0}}}(F).
 \]
 Combining this with~\eqref{P11*12.a}, we infer from \rpart{P11*10.a} of \rprop{P11*10} that
 \(
  F^\Saalpha{m+1}
  =G
 \).
\eproof

 The combination of \rpropss{P11*10}{P11*12} gives us now a complete description of the Schur-Stieltjes algorithm in the class $\SFqas{m}$.

\bthml{T11*13}
 Let \(\ug\in\R\), let $m\in\NO$, and let $\seqs{m}\in\Kggequa{m}$. For each $\ell\in\mn{0}{m}$, let $\seq{\su{j}^\sta{\ell}}{j}{0}{m-\ell}$ be the \sthaSchto{\ell}{\seqs{m}}. Let\index{s@$\HTiU{\seqs{m}}$}
 \beql{F11*5}
  \HTiU{\seqs{m}}
  \defeq\HTiu{\su{0}^\sta{0}}\circ\HTiu{\su{0}^\sta{1}}\circ\dotsb\circ\HTiu{\su{0}^\sta{m}}
 \eeq
 and let\index{s@$\HTU{\seqs{m}}$}
 \beql{F11*6}
  \HTU{\seqs{m}}
  \defeq\HTu{\su{0}^\sta{m}}\circ\HTu{\su{0}^\sta{m-1}}\circ\dotsb\circ\HTu{\su{0}^\sta{1}}\circ\HTu{\su{0}^\sta{0}}.
 \eeq
 Then the following statements hold true:
 \benui
  \il{T11*13.a} The mapping $\HTiU{\seqs{m}}$ generates a bijective correspondence between the sets $\SFdqaa{\su{0}^\sta{m}}$ and $\SFqas{m}$.
  \il{T11*13.b} Let $\HTiU{\seqs{m}}^\inv$ be the inverse mapping of $\HTiU{\seqs{m}}$. For each matrix-valued function $F\in\SFqas{m}$, then
  \[
   \HTiU{\seqs{m}}^\inv(F)
   =\HTUa{\seqs{m}}{F}
   =F^\Saalpha{m+1},
  \]
  where $F^\Saalpha{m+1}$ is the \tthSSto{(m+1)}{F}.
 \eenui
\ethm
\bproof
 Combine \rpropss{P11*10}{P11*12}.
\eproof

 Our next considerations are aimed to rewrite the mappings introduced in \rthm{T11*13} as linear fractional transformations of matrices. The essential tool in realizing this goal will be the matrix polynomials introduced in \rapp{S*C}. More precisely, we will use finite products of such matrix polynomials.
 
 Let $\ug\in\R$, let $\kappa\in\NOinf $, let $\seqska $ be a sequence of complex \tpqa{matrices}, and let $m\in\mn{0}{\kappa}$. For all $\ell\in\mn{0}{m}$, let $\seq{\su{j}^\sta{\ell}}{j}{0}{\kappa-\ell}$ be the \sthaSchto{\ell}{\seqska}. Let the sequence \(\seq{\mHTiu{\su{0}^\sta{j}}}{j}{0}{m}\) be given via \eqref{FC*2}. Then, let\index{v@$\rmiupou{s}{m}$}
\beql{F11*9}
 \rmiupou{s}{m}
 \defeq\mHTiu{\su{0}^\sta{0}}\mHTiu{\su{0}^\sta{1}}\cdot\dotso\cdot\mHTiu{\su{0}^\sta{m-1}}\mHTiu{\su{0}^\sta{m}}
\eeq
 and let\index{v@$\rmiupjkou{s}{m}$}
 \[
 \rmiupou{s}{m}
 =
 \bMat
  \rmiupnwou{s}{m}&\rmiupneou{s}{m}\\
  \rmiupswou{s}{m}&\rmiupseou{s}{m}
 \eMat
 \]
 be the block representation of $\rmiupou{s}{m}$ with \tppa{block} $\rmiupnwou{s}{m}$.

\breml{R11*16}
 Let $\ug\in\R$, let $\kappa\in\Ninf$ and let $\seqska $ be a sequence of complex \tpqa{matrices}. For all $m\in\mn{1}{\kappa}$ and all $\ell\in\mn{0}{m-1}$, one can see then from \eqref{F11*9} and~\zitaa{114arxiv}{\crem{8.3}} that
 \[
  \rmiupou{s^\sta{\ell}}{m-\ell}
  =\rmiupou{s^\sta{\ell}}{m-(\ell+1)}\mHTiu{\su{0}^\sta{m}}
 \]
 and
 \[
  \rmiupou{s^\sta{\ell}}{m-\ell}
  =\mHTiu{\su{0}^\sta{\ell}}\rmiupou{t}{m-(\ell+1)},
 \]
 where $t_j\defeq\su{j}^\sta{\ell +1}$ for all $j\in\mn{0}{m-(\ell +1)}$.
\erem

 Let $\ug\in\R$, let $\kappa\in\NOinf$, let $\seqska $ be a sequence of complex \tpqa{matrices}, and let $m\in\mn{0}{\kappa}$. For all $\ell\in\mn{0}{m}$, let $\seq{\su{j}^\sta{\ell}}{j}{0}{\kappa-\ell}$ be the \sthaSchto{\ell}{\seqska}. Let the sequence \(\seq{\mHTu{\su{0}^\sta{j}}}{j}{0}{m}\) be given via \eqref{FC*1}. Then, let\index{w@$\rmupou{s}{m}$}
\beql{F11*12}
 \rmupou{s}{m}
 \defeq\mHTu{\su{0}^\sta{m}}\mHTu{\su{0}^\sta{m-1}}\cdot\dotso\cdot\mHTu{\su{0}^\sta{0}}
\eeq
 and let\index{w@$\rmupjkou{s}{m}$}
 \[
 \rmupou{s}{m}
 =
 \bMat
  \rmupnwou{s}{m}&\rmupneou{s}{m}\\
  \rmupswou{s}{m}&\rmupseou{s}{m}
 \eMat
 \]
 be the block representation of $\rmupou{s}{m}$ with \tppa{block} $\rmupnwou{s}{m}$.

\breml{R11*18}
 Let $\ug\in\R$, let $\kappa\in\Ninf$ and let $\seqska $ be a sequence of complex \tpqa{matrices}. For all $m\in\mn{1}{\kappa}$ and all $\ell\in\mn{0}{m-1}$, one can then see from \eqref{F11*12} and~\zitaa{114arxiv}{\crem{8.3}} that
 \[
  \rmupou{s^\sta{\ell}}{m-\ell}
  =\mHTu{\su{0}^\sta{m}}\rmupou{s^\sta{\ell}}{m-(\ell+1)}
 \]
 and
 \[
  \rmupou{s^\sta{\ell}}{m-\ell}
  =\rmupou{t}{m-(\ell+1)}\mHTu{\su{0}^\sta{\ell}},
 \]
 where $t_j\defeq\su{j}^\sta{\ell+1}$ for all $j\in\mn{0}{m-(\ell+1)}$.
\erem

\bpropl{P11*20}
 Let $\ug\in\R$, let $m\in\NO$, let $\seqs{m}\in\Kggequa{m}$, and let $\seq{\su{j}^\sta{m}}{j}{0}{0}$ be the \sthaSchto{m}{\seqs{m}}. Further, let $G\in\SFdqaa{\su{0}^\sta{m}}$. Then:
 \benui
  \il{P11*20.a} For all $z\in\Cs$,
  \[
   \det\ek*{\rmiupswoua{s}{m}{z}G(z)+\rmiupseoua{s}{m}{z}}
   \neq0.
  \]
  \il{P11*20.b} Let $\HTiU{\seqs{m}}$ be given via \eqref{F11*5}. Then
  \[
   \HTiUa{\seqs{m}}{G}
   =\lftrfua{q}{q}{\rmiupou{s}{m}}{G}.
  \]
 \eenui
\eprop
\bproof
 For each $\ell\in\mn{0}{m}$ let $\seq{\su{j}^\sta{\ell}}{j}{0}{m-\ell}$ be the \sthaSchto{\ell}{\seqs{m}}. For each $\ell\in\mn{0}{m}$, in view of \rthm{T-P7*2}, then
 \(
  \seq{\su{j}^\sta{\ell}}{j}{0}{m-\ell}
  \in\Kggequa{m-\ell}.
 \)
 Thus, for each $\ell\in\mn{0}{m}$, we have
 \beql{P11*20.1}
  \seq{\su{j}^\sta{\ell}}{j}{0}{0}
  \in\Kggequa{0}.
 \eeq
 Because of \eqref{F11*5}, we have
 \beql{P11*20.2}
  \HTiU{\seqs{m}}
  =\HTiu{\su{0}^\sta{0}}\circ\HTiu{\su{0}^\sta{1}}\circ\dotsb\circ\HTiu{\su{0}^\sta{m}}.
 \eeq
 For each $\ell\in\mn{0}{m}$ we infer from \rprop{P10*1} and \rthm{T10*9} inductively
 \[
  \rk{\HTiu{\su{0}^\sta{\ell}}\circ\HTiu{\su{0}^\sta{\ell+1}}\circ\dotsb\circ\HTiu{\su{0}^\sta{m}}}](G)
  \in\SFuqaa{m-\ell}{\seq{\su{j}^\sta{\ell}}{j}{0}{m-\ell}}.
 \]
 For each $\ell\in\mn{0}{m}$, then \rlem{L5*8} shows that
 \beql{P11*20.3}
  \rk{\HTiu{\su{0}^\sta{\ell}}\circ\HTiu{\su{0}^\sta{\ell+1}}\circ\dotsb\circ\HTiu{\su{0}^\sta{m}}}(G)
  \in\SFdqaa{\su{0}^\sta{\ell}}.
 \eeq
 Taking \eqref{P11*20.1}, \eqref{P11*20.2}, and \eqref{P11*20.3} into account, we get from \rprop{L8*12} that
 \beql{P11*20.4}
  \HTiUa{\seqs{m}}{G}
  =(\lftrfu{q}{q}{\mHTiu{\su{0}^\sta{0}}}\circ\lftrfu{q}{q}{\mHTiu{\su{0}^\sta{1}}}\circ\dotsb\circ\lftrfu{q}{q}{\mHTiu{\su{0}^\sta{m}}})(G).
 \eeq
 In view of \eqref{P11*20.4} and \eqref{F11*9}, now \rprop{PB*1} yields~\eqref{P11*20.a} and~\eqref{P11*20.b}.
\eproof

\bpropl{P11*21}
 Let $\ug\in\R$, let $m\in\NO$, let $\seqs{m}\in\Kggequa{m}$, and let $F$ belong to $\SFqas{m}$. Then the following statements hold true:
 \benui
  \il{P11*21.a} For all $z\in\Cs$,
  \[
   \det\ek*{\rmupswoua{s}{m}{z}F(z)+\rmupseoua{s}{m}{z}}
   \neq0.
  \]
  \il{P11*21.b} Let $\HTU{\seqs{m}}$ be given by \eqref{F11*6}. Then
  \[
   \HTUa{\seqs{m}}{F}
   =\lftrfua{q}{q}{\rmupou{s}{m}}{F}.
  \]
 \eenui
\eprop
\bproof
 For each $\ell\in\mn{0}{m}$ let $\seq{\su{j}^\sta{\ell}}{j}{0}{m-\ell}$ be the \sthaSchto{\ell}{\seqs{m}}. For each $\ell\in\mn{0}{m}$ in view of \rthm{T-P7*2}, then
 \beql{P11*21.0}
  \seq{\su{j}^\sta{\ell}}{j}{0}{m-\ell}
  \in\Kggequa{m-\ell}.
 \eeq
 Because of \eqref{F11*6}, we have
 \beql{P11*21.2}
  \HTU{\seqs{m}}
  =\HTu{\su{0}^\sta{m}}\circ\HTu{\su{0}^\sta{m-1}}\circ\dotsb\circ\HTu{\su{0}^\sta{0}}.
 \eeq
 For $m=0$, the assertions of~\eqref{P11*21.a} and~\eqref{P11*21.b} are an immediate consequence of \eqref{P11*21.2}, \eqref{F11*12}, and \rprop{L8*9}. Now let $m\in\N$. For each $\ell\in\mn{0}{m-1}$, we infer from \rthm{T9*7} inductively that
 \beql{P11*21.3}
  \rk{\HTu{\su{0}^\sta{\ell}}\circ\HTu{\su{0}^\sta{\ell-1}}\circ\dotsb\circ\HTu{\su{0}^\sta{0}}}(F)
  \in\SFuqAA{m-\ell}{\seq{\su{j}^\sta{\ell}}{j}{0}{m-\ell}}.
 \eeq
 Taking \eqref{P11*21.0}, \eqref{P11*21.2}, and \eqref{P11*21.3} into account, \rprop{L8*9} yields
 \beql{P11*21.4}
  \HTUa{\seqs{m}}{F}
  =(\lftrfu{q}{q}{\mHTu{\su{0}^\sta{m}}}\circ\lftrfu{q}{q}{\mHTu{\su{0}^\sta{m-1}}}\circ\dotsb\circ\lftrfu{q}{q}{\mHTu{\su{0}^\sta{0}}})(F).
 \eeq
 In view of \eqref{P11*21.4} and \eqref{F11*12}, now \rprop{PB*1} yields~\eqref{P11*21.a} and~\eqref{P11*21.b}.
\eproof

\section{Description of the sets $\SFqas{\kappa}$}\label{S*12}
 Let \(\ug\in\R\), let $m\in\NO$, and let $\seqs{m}\in\Kggequa{m}$. Then we know from \rthmp{T1*3+2}{T1*3+2.a} that the solution set $\MggqKsg{m}$ of \rprob{\mproblem{\rhl}{m}{=}} is \tne{}. Furthermore, \rthm{T-C5*6} tells us that $\SFqas{m}$ is exactly the set of \taSt{s} of all the measures belonging to $\MggqKsg{m}$. On the basis of our Schur-Stieltjes-type algorithm, which was introduced in \rsec{S*11}, we have obtained important insights into the structure of the set $\SFqas{m}$. Using \rthm{T7*5}, we will now rewrite the descriptions of $\SFqas{m}$, which were obtained in \rsec{S*11}, in a form which is better adapted to the original data sequence $\seqs{m}$. This is our first main theorem.

\bthml{T12*1}
 Let \(\ug\in\R\), let $m\in\NO$, let $\seqs{m}\in\Kggequa{m}$ with \taSp{} \(\seqQ{m}\), and let $\rmiupou{s}{m}$ be defined by \eqref{F11*9}. Then:
 \benui
  \il{T12*1.a} $\SFqas{m}=\lftrfuA{q}{q}{\rmiupou{s}{m}}{\SFdqaa{\Spu{m}}}$.
  \il{T12*1.b} For each $F\in\SFqas{m}$, there is a unique $G\in\SFdqaa{\Spu{m}}$ which satisfies
  \[
   \lftrfuA{q}{q}{\rmiupou{s}{m}}{G}
   =F,
  \]
  namely $G=F^\Saalpha{m+1}$, where $F^\Saalpha{m+1}$ stands for the \tthSSto{(m+1)}{F}.
 \eenui
\ethm
\bproof
 Let $\seq{\su{j}^\sta{m}}{j}{0}{0}$ be the \sthaSchto{m}{\seqs{m}}. In view of \rthm{T7*5} and \rdefn{D7*3}, then $\su{0}^\sta{m}=\Spu{m}$. Hence, the application of \rthm{T11*13} and \rprop{P11*20} completes the proof.
\eproof

 It should be mentioned that a result similar to \rthm{T12*1} for \(\ug=0\) is contained in~\zitaa{MR1802495}{\cthm{4.1(d)}}.

\bcorl{C12*2}
 Let \(\ug\in\R\), let $m\in\NO$, let $\seqs{m}\in\Kgqua{m}$, and let $\rmiupou{s}{m}$ be defined by \eqref{F11*9}. Then
 \[
  \SFqaS{m}
  =\lftrfuA{q}{q}{\rmiupou{s}{m}}{\SFdqa}.
 \]
\ecor
\bproof
 From \rprop{P-R7*10} we get
 \beql{F12*1}
  \seqs{m}
  \in\Kggequa{m}.
 \eeq
 Because of the assumption $\seqs{m}\in\Kgqua{m}$, we infer from \rprop{P-L7*11} that the matrix $\Spu{m}$ defined in \rdefn{D7*3} is \tpH{}. In particular, $\det\Spu{m}\neq0$. Consequently, \rrem{R4*1} implies $\SFdqaa{\Spu{m}}=\SFdqa$. Thus, taking \eqref{F12*1} into account, the application of \rthmp{T12*1}{T12*1.a} completes the proof.
\eproof

 Let \(\ug\in\R\), let $m\in\NO$, and let $\seqs{m}\in\Kggqua{m}$ with \taSp{} \(\seqQ{m}\). Then $\seqs{m}$ is called \notii{\tcd{}}, if $\Spu{m}=\Oqq$. Observe that the set $\Kggdqua{m}$\index{k@$\Kggdqua{m}$} of all \tcd{} sequences belonging to $\Kggqua{m}$ is a subclass of $\Kggequa{m}$ (see~\zitaa{MR3014201}{\cprop{5.9}}). This class is connected to the case of a unique solution, which was already discussed in~\zitaa{MR1802495}{\cthm{4.1(a)--(c)}}.

\bthml{T12*3}
 Let \(\ug\in\R\), let $m\in\NO$, and let $\seqs{m}\in\Kggequa{m}$. Then:
 \benui
  \item The following statements are equivalent:
  \baeqii{0}
   \il{T12*3.i} The set $\SFqas{m}$ consists of exactly one element.
   \il{T12*3.ii} The sequence $\seqs{m}\) belongs to \(\Kggdqua{m}\).
  \eaeqii
  \il{T12*3.b} If~\ref{T12*3.i} holds true, then $\det\rmiupseoua{s}{m}{z}\neq0$ for all $z\in\Cs$ and
  \[
   \SFqaS{m}
   =\set*{\rmiupneou{s}{m}(\rmiupseou{s}{m})^\inv}.
  \]
 \eenui
\ethm
\bproof
 \bimp{T12*3.i}{T12*3.ii}
  Use \rthm{T3*2}, \rrem{R3*18}, \eqref{F5*3}, \rthm{T-C5*6}, and~\zitaa{MR2735313}{\cthmss{6.1}{6.3}}.
 \eimp
 
 \bimp{T12*3.ii}{T12*3.i}
  Let $\seqQ{m}$ be the \taSpo{\seqs{m}}. Then~\ref{T12*3.ii} means $\Spu{m}=\Oqq$. Thus, \rpropp{P4*9}{P4*9.a}, implies that the set $\SFdqaa{\Spu{m}}$ consists of exactly one element, namely the constant function defined on $\Cs$ with value $\Oqq$. Using \rthmp{T12*1}{T12*1.a}, then~\ref{T12*3.i} follows, and we see moreover that~\eqref{T12*3.b} holds true.
 \eimp
\eproof

By using \rprop{P4*9}, we are able to derive an alternate description of the set $\SFqas{m}$ if $m\in\NO$ and $\seqs{m}\in\Kggequa{m}$ are arbitrarily given. This gives reformulations of our main result stated in \rthm{T12*1}.

\bthml{T12*7}
 Let $m\in\NO$ and let $\seqs{m}\in\Kggequa{m}$ with \taSp{} $\seqQ{m}$. Suppose that $r\defeq\rank\Spu{m}$ fulfills $r\geq1$. Let $u_1,u_2,\dotsc,u_r$ be an orthonormal basis of $\ran{\Spu{m}}$ and let $U\defeq\mat{ u_1,u_2,\dotsc,u_r}$. Then:
 \benui
  \il{T12*7.a} $\SFqas{m}=\lftrfua{q}{q}{\rmiupou{s}{m}}{\setaa{UfU^\ad}{f\in\SFdua{r}}}$.
  \il{T12*7.b} Let $F\in\SFqas{m}$. Then there is a unique $f\in\SFdua{r}$ such that
  \[
   \lftrfua{q}{q}{\rmiupou{s}{m}}{UfU^\ad}
   =F,
  \]
  namely $f\defeq U^\ad F^\Saalpha{m+1}U$, where $F^\Saalpha{m+1}$ is the \tthSSto{(m+1)}{F}.
 \eenui
\ethm
\bproof
 Since $\seqs{m}\in\Kggequa{m}$ implies $\seqs{m}\in\Kggqua{m}$ we get from \rprop{P-L7*7} that $\Spu{m}\in\Cggq$. In particular, $\Spu{m}^\ad=\Spu{m}$. Thus, \rpropp{P4*9}{P4*9.b1} yields
 \beql{F12*3}
  \SFdqaa{\Spu{m}}
  =\setaa*{UfU^\ad}{f\in\SFdua{r}},
 \eeq
 whereas \rrem{RA*6} implies $U^\ad U=\Iu{r}$.
 
 \eqref{T12*7.a} Because of \eqref{F12*3}, we infer~\eqref{T12*7.a} from \rthmp{T12*1}{T12*1.a}.
 
 \eqref{T12*7.b} Using \eqref{F12*3}, $U^\ad U=\Iu{r}$, and~\eqref{T12*7.a}, we conclude~\eqref{T12*7.b} from \rthmp{T12*1}{T12*1.b}.
\eproof

\appendix
\section{Some considerations on \tnnH{} measures}\label{A1608}
 In this appendix, we summarize some facts on integration with respect to \tnnH{} measures. For a detailed treatment of this subject we refer to Rosenberg~\zita{MR0163346} (see also Kats~\zita{MR0080280}). Let \((\Omega,\mathfrak{A})\) be a measurable space. If \(\mu = \matauuo{\mu_{jk}}{j,k}{1}{q}\) is a \tnnH{} measure on \((\Omega,\mathfrak{A})\), then each entry function \(\mu_{jk}\) is a complex measure on \((\Omega,\mathfrak{A})\). In particular, \(\mu_{11},\mu_{22},\dotsc,\mu_{qq}\) are finite \tnn{} real-valued measures. For each \(H\in\Cggq\), the inequality \(H\leq(\tr H)\Iq\) holds true. Hence, each entry function \(\mu_{jk}\) is absolutely continuous with respect to the so-called trace measure \(\tau\defeq\sum_{j=1}^q\mu_{jj}\)\index{t@$\tau$} of \(\mu\), \ie{}, for each \(M\in\mathfrak{A}\), which satisfies \(\tau(M)=0\), it follows \(\mu(M)=\Oqq\). The Radon-Nikodym derivatives \symba{\dif\mu_{jk}/\dif\tau}{d} are thus well defined up to sets of zero \(\tau\)\nobreakdash-measure. Obviously, the matrix-valued function \(\mu_\tau'\defeq\matauuo{\dif\mu_{jk}/\dif\tau}{j,k}{1}{q}\)\index{$\mu_\tau'$} is \(\mathfrak{A}\)\nobreakdash-measurable and integrable with respect to \(\tau\). The matrix-valued function \(\mu_\tau'\) is said to be the trace derivative of \(\mu\). If \(\nu\) is a \tnn{} real-valued measure on \(\mathfrak{A}\), then let the class of all \(\mathfrak{A}\)\nobreakdash-measurable \tpqa{matrix-valued} functions \(\Phi=\matau{\phi_{jk}}{\substack{j=1,\dotsc,p\\ k=1,\dotsc,q}}\) on \(\Omega\) such that each \(\phi_{jk}\) is integrable with respect to \(\nu\) be denoted by \symba{\pqLaaaC{p}{q}{\Omega}{\mathfrak{A}}{\nu}}{l}. An ordered pair \((\Phi,\Psi)\) consisting of an \(\mathfrak{A}\)\nobreakdash-measurable \tpqa{matrix-valued} function \(\Phi\) on \(\Omega\) and an \(\mathfrak{A}\)\nobreakdash-measurable \taaa{r}{q}{matrix-valued} function \(\Psi\) on \(\Omega\) is said to be integrable with respect to a \tnnH{} measure \(\mu\) on \((\Omega,\mathfrak{A})\) if \(\Phi\mu_\tau'\Psi^\ad\) belongs to \(\pqLaaaC{p}{r}{\Omega}{\mathfrak{A}}{\tau}\), where \(\tau\) is the trace measure of \(\mu\). In this case, the integral of \((\Phi,\Psi)\) with respect to \(\mu\) is defined by\index{$\int_\Omega\Phi\dif\mu\Psi^\ad$}
 \[
  \int_\Omega\Phi\dif\mu\Psi^\ad
  \defeq\int_\Omega\Phi\mu_\tau'\Psi^\ad\dif\tau
 \]
 and for any \(M\in\mathfrak{A}\), the pair \((1_M\Phi,1_M\Psi)\) is integrable with respect to \(\mu\), where \symba{1_M}{1} is the indicator function of the set \(M\). Then the integral of \((\Phi,\Psi)\) with respect to \(\mu\) over \(M\) is defined by \(\int_M\Phi\dif\mu\Psi^\ad\defeq\int_\Omega\rk{1_M\Phi}\mu_\tau'\rk{1_M\Psi}^\ad\dif\tau\)\index{$\int_M\Phi\dif\mu\Psi^\ad$}. An \(\mathfrak{A}\)\nobreakdash-measurable complex-valued function \(f\) on \(\Omega\) is said to be integrable with respect to a \tnnH{} measure \(\mu\) on \((\Omega,\mathfrak{A})\) if the pair \((f\Iq,\Iq)\) is integrable with respect to \(\mu\). In this case the integral of \(f\) with respect to \(\mu\) is defined by\index{$\int_\Omega f\dif\mu$}
 \[
  \int_\Omega f\dif\mu
  \defeq\int_\Omega(f\Iq)\dif\mu\Iq^\ad
 \]
 and for any \(M\in\mathfrak{A}\), the function \(1_M f\) is integrable with respect to \(\mu\). Then the integral of \(f\) with respect to \(\mu\) over \(M\) is defined by \(\int_M f\dif\mu\defeq\int_M(f\Iq)\dif\mu\Iq^\ad\)\index{$\int_M f\dif\mu$}. We denote by \symba{\Loaaaa{1}{\Omega}{\mathfrak{A}}{\mu}{\C}}{l} the set of all \(\mathfrak{A}\)\nobreakdash-measurable complex-valued functions \(f\) on \(\Omega\) which are integrable with respect to a \tnnH{} measure \(\mu\) on \((\Omega,\mathfrak{A})\). 
 

 We consider a \tnnH{} measure \(\mu\) on \((\Omega,\mathfrak{A})\) with trace measure \(\tau\) and an \(\mathfrak{A}\)\nobreakdash-measurable complex-valued function \(f\) on \(\Omega\):
 
\breml{R1643}
 The function \(f\) belongs to \(\LaaaC{\Omega}{\mathfrak{A}}{\mu}\) if and only if \(f\in\LaaaC{\Omega}{\mathfrak{A}}{\tau}\).
\erem

\breml{R1141}
 Let \(u\in\Cq\). Then \(\nu\defeq u^\ad\mu u\) is a finite measure on \((\Omega,\mathfrak{A})\) which is absolutely continuous with respect to \(\tau\). If \(f\) belongs to \(\Loaaaa{1}{\Omega}{\mathfrak{A}}{\mu}{\C}\), then \(\int_\Omega\abs{f}\dif\nu<\infp\) and \(\int_\Omega f\dif\nu=u^\ad\rk{\int_\Omega f\dif\mu}u\).
\erem

\breml{R1558}
 The function \(f\) belongs to \(\Loaaaa{1}{\Omega}{\mathfrak{A}}{\mu}{\C}\) if and only if \(\int_\Omega\abs{f}\dif(u^\ad\mu u)<\infp\) for all \(u\in\Cq\).
\erem

\breml{R1647}
 If \(f\in\Loaaaa{1}{\Omega}{\mathfrak{A}}{\mu}{\C}\), then the functions \(\re f\) and \(\im f\) both belong to \(\Loaaaa{1}{\Omega}{\mathfrak{A}}{\mu}{\C}\) and \(\int_\Omega\re f\dif\mu=\re\int_\Omega f\dif\mu\) and \(\int_\Omega\im f\dif\mu=\im\int_\Omega f\dif\mu\).
\erem

\bpropnl{Lebesgue's dominated convergence for \tnnH{} measures~\zitaa{142}{\cprop{A.6}}}{P1102}
 Let \(\mu\) be a \tnnH{} \tqqa{measure} on a measurable space \((\Omega,\mathfrak{A})\) with trace measure \(\tau\). For all \(n\in\N\), let \(f_n\colon\Omega\to\C\) be an \(\mathfrak{A}\)\nobreakdash-measurable function. Let \(f\colon\Omega\to\C\) be an \(\mathfrak{A}\)\nobreakdash-measurable function and let \(g\in\LaaaC{\Omega}{\mathfrak{A}}{\mu}\) be such that \(\lim_{n\to\infi}f_n(\omega) = f(\omega)\) for  \(\tau\)\nobreakdash-a.\,a.\ \(\omega\in\Omega\) and that \(\abs{f_n(\omega)}\leq\abs{g(\omega)}\) for all \(n\in\N\) and \(\tau\)\nobreakdash-a.\,a.\ \(\omega\in\Omega\). Then \(f\in\LaaaC{\Omega}{\mathfrak{A}}{\mu}\), \(f_n\in\LaaaC{\Omega}{\mathfrak{A}}{\mu}\) for all \(n\in\N\), and \( \lim_{n\to\infi}\int_\Omega f_n\dif\mu=\int_\Omega f\dif\mu\).
\eprop

\section{Some Results On Moore-Penrose Inverses of Matrices}\label{S*A}
For the convenience of the reader, we state some well-known and some special results on Moore-Penrose inverses of matrices (see \eg{}, Rao/Mitra~\zita{MR0338013} or~\zita{MR1152328}{\csec{1}}). If $A\in\Cpq$, then (by definition) the Moore-Penrose inverse $A^\mpi$ of $A$ is the unique matrix $A^\mpi\in\Cqp$\index{$A^\mpi$} which satisfies the four equations
\begin{align*}
 AA^\mpi A&=A,&A^\mpi AA^\mpi&=A^\mpi,&(AA^\mpi)^\ad&=AA^\mpi,&
 &\text{and}&
 (A^\mpi A)^\ad&=A^\mpi A.
\end{align*}

\breml{RA*1}
 Let $A\in\Cpq$. Then one can easily check that:
 \benui
  \il{RA*1.a} $(A^\mpi)^\mpi=A$, $(A^\mpi)^\ad=(A^\ad)^\mpi$, and $\Ip-AA^\mpi\in\Cggq$.
  \il{RA*1.b} $\ran{A^\mpi}=\ran{A^\ad}$, $\rank(A^\mpi)=\rank A$, and $\nul{A^\mpi}=\nul{A^\ad}$.
 \eenui
\erem

\begin{prop}[see, \eg{}~\zita{MR1152328}{\cthm{1.1.1}}]\label{PA*2}
 If $A\in\Cpq$ then a matrix $G\in\Cqp$ is the Moore-Penrose inverse of $A$ if and only if $AG=P_A$ and $GA=P_G$, where $P_A$ and $P_G$ are respectively, the matrices associated with the orthogonal projection in $\Cp$ onto $\ran{A}$ and the orthogonal projection in $\Cq$ onto $\ran{G}$.
\end{prop}

\breml{RA*3}
 Let $A\in\Cpq$. Then it is readily checked that:
 \benui
  \il{RA*3.a} Let $r\in\N$ and let $B\in\Coo{r}{q}$. Then $\nul{A}\subseteq\nul{B}$ if and only if $BA^\mpi A=B$. Furthermore, $\nul{B}=\nul{A}$ if and only if $B^\mpi B=A^\mpi A$.
  \il{RA*3.b} Let $s\in\N$ and let $C\in\Coo{p}{s}$. Then $\ran{C}\subseteq\ran{A}$ if and only if $AA^\mpi C=C$. Furthermore, $\ran{C}=\ran{A}$ if and only if $CC^\mpi=AA^\mpi$.
 \eenui
\erem

\blemnl{\zitaa{MR3380267}{\clem{A.5}}}{LA*5}
 Let $A,X\in\Cpq$. Then the following statements are equivalent:
\baeqi{0}
  \il{LA*5.i} $\nul{A}\subseteq\nul{X}$ and $\ran{A}\subseteq\ran{X}$.
  \il{LA*5.ii} $\nul{A}=\nul{X}$ and $\ran{A}=\ran{X}$.
  \il{LA*5.iii} $A^\mpi A=X^\mpi X$ and $AA^\mpi=XX^\mpi$.
  \il{LA*5.iv} $\nul{A^\mpi}=\nul{X^\mpi}$ and $\ran{A^\mpi}=\ran{X^\mpi}$.
 \eaeqi
\elem

\breml{RA*6}
 Let $A\in\Cpq\setminus\set{\Opq}$. Let $r\defeq\rank A$, let $u_1,u_2,\dotsc,u_r$ be an orthonormal basis of $\ran{A^\ad}$, and let $U\defeq\mat{u_1,u_2,\dotsc,u_r}$. Then $U^\ad U=\Iu{r}$, and, in view of \rprop{PA*2} and \rremp{RA*1}{RA*1.b}, furthermore $UU^\ad=A^\mpi A$.
\erem

At the end of this section, we give a slight generalization of a result due to S.~L.~Campbell and C.~D.~Meyer~Jr. This result can be proved by an obvious modification of the proof given in~\zitaa{MR1105324}{\cthm{10.4.1}}.

\begin{lem}\label{LA*10}
 Suppose that $\seq{A_n}{n}{1}{\infi}$ is a sequence of complex \tpqa{matrices} which converges to a complex \tpqa{matrix} $A$. Then $\seq{A_n^\mpi}{n}{1}{\infi}$ is convergent if and only if there is a positive integer $m$ such that $\rank A_n=\rank A$ for each integer $n$ with $n\geq m$. In this case, $\seq{A_n^\mpi}{n}{1}{\infi}$ converges to $A^\mpi$.
\end{lem}

\section{On Linear Fractional Transformations of Matrices}\label{S*B}
 In this appendix, we summarize some basic facts on linear fractional transformations of matrices which are needed in the paper. This material is mostly taken from~\zita{MR566141} and~\zitaa{MR1152328}{\csec{1.6}}.

Let $a\in\Coo{p}{p}$, $b\in\Coo{p}{q}$, $c\in \Cqp$, $d\in\Cqq$, and let 
\[ 
 E
 \defeq
 \bMat
  a&b\\
  c&d
 \eMat.
\]
If the set\index{q@$\dblftruu{c}{d}$}
\[
 \dblftruu{c}{d}
 \defeq\setaa*{x\in\Cpq}{\det(cx+d)\neq0}
\]
is \tne{}, then the linear fractional transformation $\lftroou{p}{q}{E}\colon\dblftruu{c}{d}\to\Coo{p}{q}$\index{s@$\lftroou{p}{q}{E}$} is defined by
\beql{F*lft}
 \lftrooua{p}{q}{E}{x}
 \defeq(ax+b)(cx+d)^\inv.
\eeq

 The following well-known result shows that the composition of two linear fractional transformations is again a mapping of this type.

\begin{prop}[see, \eg{}~\zitaa{MR1152328}{\cprop{1.6.3}}]\label{PB*1}
 Let $a_1,a_2\in\Cpp$, let $b_1,b_2\in\Cpq$, let $c_1,c_2\in\Cqp$, and let $d_1,d_2\in\Cqq$ be such that
 \[
  \rank\mat{ c_1,d_1}
  =\rank\mat{ c_2,d_2}
  =q.
 \]
 Furthermore, let $E_1\defeq\bsma a_1&b_1\\ c_1&d_1\esma$, $E_2\defeq\bsma a_2&b_2\\ c_2&d_2\esma$, $E\defeq E_2E_1$, and $E=\bsma a&b\\c&d\esma$ be the block representation of $E$ with \tppa{block} $a$. Then $\mathcal{Q}\defeq\setaa{x\in\dblftruu{c_1}{d_1}}{\lftrooua{p}{q}{E_1}{x}\in\dblftruu{c_2}{d_2}}$ is a nonempty subset of the set $\dblftruu{c}{d}$ and $\lftrooua{p}{q}{E_2}{\lftrooua{p}{q}{E_1}{x}}=\lftrooua{p}{q}{E}{x}$ holds true for all $x\in\mathcal{Q}$.
\end{prop}

We make the following convention: If a \tne{} subset $\mathcal{G}$ of $\C$ and a matrix-valued function $V\colon\mathcal{G}\to\Coo{2q}{2q}$ with \tqqa{block} partition $V=\bsma v_{11}&v_{12}\\v_{21}&v_{22}\esma$ and a matrix-valued function $F\colon\mathcal{G}\to\Cqq$ with $\det[v_{21}(z)F(z)+v_{22}(z)]\neq0$ for all $z\in\mathcal{G}$ are given, then we will use the notation $\lftrfua{q}{q}{V}{F}$\index{s@$\lftrfua{q}{q}{V}{F}$} for the function $\lftrfua{q}{q}{V}{F}\colon\mathcal{G}\to\Cqq$ defined by $\lftrfuaa{q}{q}{V}{F}{z}\defeq\lftrooua{q}{q}{V(z)}{F(z)}$ for all $z\in\mathcal{G}$.

\section{The Matrix Polynomials $\mHTiu{A}$ and $\mHTu{A}$}\label{S*C}
 In this appendix, we study special linear \taaa{(p+q)}{(p+q)}{matrix} polynomials which are intensively used in \rsec{S*8}. Let \(\ug\in\C\) and let $A\in\Cpq$. Then we define the mappings $\mHTu{A}\colon\C\to\Coo{(p+q)}{(p+q)}$\index{w@$\mHTu{A}$} and $\mHTiu{A}\colon\C\to\Coo{(p+q)}{(p+q)}$\index{v@$ \mHTiu{A}$} by
\beql{FC*1}
 \mHTua{A}{z}
 \defeq
 \begin{pmat}[{|}]
  (z-\ug)\Ip&A\cr\-
  -(z-\ug)A^\mpi&\Iq-A^\mpi A\cr
 \end{pmat}
\eeq
and
\beql{FC*2}
 \mHTiua{A}{z}
 \defeq
 \begin{pmat}[{|}]
  \Opp&-A\cr\-
  (z-\ug)A^\mpi&(z-\ug)\Iq\cr
 \end{pmat}.
\eeq

 The use of the matrix polynomial $ \mHTiu{A}$ was inspired by some constructions in the paper~\zita{MR2038751}. In particular, we mention~\zitaa{MR2038751}{formula~(2.3)}. In their constructions Chen and Hu used Drazin inverses instead of Moore-Penrose inverses of matrices. Since both types of generalized inverses coincide for \tH{} matrices (see~\zitaa{MR3014199}{\cprop{A.2}}) we can conclude that in the generic case the matrix polynomials $ \mHTiu{A}$ coincide for \(\ug=0\) with the objects used in~\zita{MR2038751}.

\breml{RC*1}
 Let $A\in\Cpq$ and let $\ug,z\in\C$. Then one can easily see that
 \[
  \mHTiua{A}{z}\mHTua{A}{z}
  =(z-\ug)\diag\mat{AA^\mpi,\Iq}
 \]
 and
 \[
  \mHTua{A}{z}\mHTiua{A}{z}
  =(z-\ug)\diag\mat{AA^\mpi,\Iq}.
 \]
\erem

Now we are going to study the linear fractional transformation generated by the matrix $\mHTua{A}{z}$ for arbitrarily given $z\in\C\setminus\set{\ug}$.

\bleml{LC*2}
 Let \(\ug\in\C\), let $A\in\Cpq$, and let $z\in\C\setminus\set{\ug}$. Then:
 \benui
  \il{LC*2.a} The matrix $-(z-\ug)^\inv A$ belongs to $\dblftruu{-(z-\ug)A^\mpi}{\Iq-A^\mpi A}$. In particular, $\dblftruu{-(z-\ug)A^\mpi}{\Iq-A^\mpi A}\neq\emptyset$.
  \il{LC*2.b} Let $X\in\Cpq$ be such that $\ran{A}\subseteq\ran{X}$ and $\nul{A}\subseteq\nul{X}$. Then $X\in\dblftruu{-(z-\ug)A^\mpi}{\Iq-A^\mpi A}$ and $\ek{-(z-\ug)A^\mpi X +\Iq-A^\mpi A}^\inv=-(z-\ug)^\inv X^\mpi A+\Iq-A^\mpi A$.
 \eenui
\elem
\bproof
 First observe, that \(w\defeq z-\ug\neq0\).
 
 \eqref{LC*2.a} This follows from $\rk{-wA^\mpi}\rk{-w^\inv A}+\rk{\Iq-A^\mpi A}=\Iq$.
 
 \eqref{LC*2.b} In view of \rlem{LA*5}, we have $AA^\mpi=XX^\mpi$ and $A^\mpi A=X^\mpi X$. Hence,
 \begin{equation*}\begin{split}
  &(-wA^\mpi X +\Iq-A^\mpi A)(-w^\inv X^\mpi A+\Iq-A^\mpi A)\\
  &=A^\mpi XX^\mpi A-wA^\mpi X+wA^\mpi XA^\mpi A\\
  &\qquad-w^\inv X^\mpi A+\Iq-A^\mpi A+w^\inv A^\mpi AX^\mpi A-A^\mpi A+A^\mpi AA^\mpi A\\
  &=A^\mpi AA^\mpi A-wA^\mpi X+wA^\mpi XX^\mpi X\\
  &\qquad-w^\inv X^\mpi A+\Iq-A^\mpi A+w^\inv X^\mpi XX^\mpi A-A^\mpi A+A^\mpi A\\
  &=A^\mpi A-wA^\mpi X+wA^\mpi X-w^\inv X^\mpi A+\Iq-A^\mpi A+w^\inv X^\mpi A
  =\Iq.
 \end{split}\end{equation*}
 This completes the proof of~\eqref{LC*2.b}. 
\eproof

\bleml{LC*3}
 Let \(\ug\in\C\), let $A\in\Cpq$, let $\mHTu{A}$ be defined via \eqref{FC*1}, and let $X\in\Cpq$ be such that the inclusions $\ran{A}\subseteq\ran{X}$ and $\nul{A}\subseteq\nul{X}$ are satisfied. Furthermore, let $z\in\C\setminus\set{\ug}$. Then the matrix $X$ belongs to $\dblftruu{-(z-\ug)A^\mpi}{\Iq-A^\mpi A}$. Furthermore,
  \[
  \lftrooua{p}{q}{\mHTua{A}{z}}{X}
  =-A\ek*{\Iq+(z-\ug)^\inv X^\mpi A}
  \]
  and
  \begin{align*}
   \Ran{\lftrooua{p}{q}{\mHTua{A}{z}}{X}}&\subseteq\ran{A},&
   \nul{A}&\subseteq\Nul{\lftrooua{p}{q}{\mHTua{A}{z}}{X}}.
  \end{align*}
\elem
\bproof
 First observer that \(w\defeq z-\ug\neq0\). In view of \rlemp{LC*2}{LC*2.b}, we have $X\in\dblftruu{-wA^\mpi}{\Iq-A^\mpi A}$ and
 \beql{LC*3.1}
  \rk{-wA^\mpi X +\Iq-A^\mpi A}^\inv
  =-w^\inv X^\mpi A+\Iq-A^\mpi A.
 \eeq
 Because of the choice of $X$, \rpartss{RA*3.a}{RA*3.b} of \rrem{RA*3} yield $XA^\mpi A=X$ and $XX^\mpi A=A$, and, in view of \eqref{LC*3.1}, then
 \begin{equation}
 \begin{split}\label{LC*3.4}
  \lftrooua{p}{q}{\mHTua{A}{z}}{X}
  &=\rk{wX+A}(-wA^\mpi X+\Iq-A^\mpi A)^\inv\\
  &=\rk{wX+A}(-w^\inv X^\mpi A+\Iq-A^\mpi A)\\
  &=-XX^\mpi A-w^\inv AX^\mpi A+\rk{wX+A}(\Iq-A^\mpi A)\\
  &=-A-w^\inv AX^\mpi A
  =-A(\Iq+w^\inv X^\mpi A).
 \end{split}
 \end{equation}
 The remaining assertions are immediate consequences of \eqref{LC*3.4}.
\eproof

\breml{RC*4}
 Let \(\ug\in\C\), let $A\in\Cpq$, let $z\in\C\setminus\set{\ug}$, and let $X\in\dblftruu{(z-\ug)A^\mpi}{(z-\ug)\Iq}$. Then from \eqref{FC*2} we see that
 \[
  \lftrooua{p}{q}{\mHTiua{A}{z}}{X}
  =-(z-\ug)^\inv A\rk{\Iq+A^\mpi X}^\mpi
 \]
 and, in view of $\det[(z-\ug)A^\mpi X+(z-\ug)\Iq]\neq0$, thus $\ran{\lftrooua{p}{q}{\mHTiua{A}{z}}{X}}=\ran{A}$.
\erem

\bibliography{141k}
\bibliographystyle{abbrv}

\vfill\noindent
\begin{minipage}{0.5\textwidth}
 Universit\"at Leipzig\\
Fakult\"at f\"ur Mathematik und Informatik\\
PF~10~09~20\\
D-04009~Leipzig
\end{minipage}
\begin{minipage}{0.49\textwidth}
 \begin{flushright}
  \texttt{
   fritzsche@math.uni-leipzig.de\\
   kirstein@math.uni-leipzig.de\\
   maedler@math.uni-leipzig.de
  } 
 \end{flushright}
\end{minipage}

\end{document}